\documentclass{amsart}

\usepackage{amsmath} 
\usepackage{amssymb}

\newtheorem{theorem}{Theorem}[section] 
\newtheorem{clx}{Claim}[theorem]
\newtheorem{lemma}[theorem]{Lemma} 
\newtheorem{proposition}[theorem]{Proposition} 
\newtheorem{observation}[theorem]{Observation} 
\newtheorem{corollary}[theorem]{Corollary} 

\theoremstyle{definition}
\newtheorem{definition}[theorem]{Definition}

\newtheorem{problem}[theorem]{Problem}

\theoremstyle{remark}
\newtheorem{remark}[theorem]{Remark}

\numberwithin{equation}{section}
\newcommand{\forces}{\Vdash}

\newcommand{\bV}{{\bf V}} 
\newcommand{\lesdot}{\mathrel{\mathord{<}\!\!\raise 
0.8 pt\hbox{$\scriptstyle\circ$}}}

\newcommand{\conc}{{}^\frown\!}
\newcommand{\lh}{{\rm lh}\/}
\newcommand{\rest}{{\restriction}}
\newcommand{\Dom}{{\rm Dom}}

\newcommand{\cD}{{\mathcal D}}

\newcommand{\cH}{{\mathcal H}}

\newcommand{\cI}{{\mathcal I}}

\newcommand{\bbP}{{\mathbb P}}
\newcommand{\cP}{{\mathcal P}}

\newcommand{\bbQ}{{\mathbb Q}}

\newcommand{\cT}{{\mathcal T}}
\newcommand{\cU}{{\mathcal U}}

\newcommand{\cX}{{\mathcal X}}
\newcommand{\cY}{{\mathcal Y}}

\newcommand{\cf}{{\rm cf}\/} 
 
\newcommand{\fil}{{\rm fil}\/}

\newcommand{\otp}{{\rm otp}\/}

\newcommand{\st}{{\bf st}} 

\newcommand{\vare}{\varepsilon}

\newcommand{\bqz}{{{\mathbb Q}^0_\lambda}}

\newcommand{\bqd}{{{\mathbb Q}^*_\lambda}}

\newcommand{\Agame}{{\Game^{\rm rcA}_{\bar{\mu}}}}
\newcommand{\agame}{{\Game^{{\rm tree}{\bf A}}_{\bar{\mu}}}}
\newcommand{\tagame}{{\Game^{{\rm rc}{\bf 2a}}_{\bar{\mu}}}}
\newcommand{\bagame}{{\Game^{{\rm rc}{\bf 2b}}_{\bar{\mu},\cU}}}

\newcommand{\pos}{{\rm pos}}
\newcommand{\vcom}{{\rm COM}^\bV}
\newcommand{\vinc}{{\rm INC}^\bV}
\newcommand{\fcom}{{\rm COM}^{\bV^{\bbP_\gamma}}}
\newcommand{\finc}{{\rm INC}^{\bV^{\bbP_\gamma}}}

\newcommand{\rk}{{\rm rk}}
\newcommand{\vtl}{\vartriangleleft}

\newcount\skewfactor
\def\mathunderaccent#1#2 {\let\theaccent#1\skewfactor#2
\mathpalette\putaccentunder}
\def\putaccentunder#1#2{\oalign{$#1#2$\crcr\hidewidth
\vbox to.2ex{\hbox{$#1\skew\skewfactor\theaccent{}$}\vss}\hidewidth}}
\def\name{\mathunderaccent\tilde-3 }

\begin{document}

\title{Reasonable ultrafilters, again}

\author{Andrzej Ros{\l}anowski}
\address{Department of Mathematics\\
 University of Nebraska at Omaha\\
 Omaha, NE 68182-0243, USA}
\email{roslanow@member.ams.org}
\urladdr{http://www.unomaha.edu/logic}

\author{Saharon Shelah}
\address{Einstein Institute of Mathematics\\
Edmond J. Safra Campus, Givat Ram\\
The Hebrew University of Jerusalem\\
Jerusalem, 91904, Israel\\
 and  Department of Mathematics\\
 Rutgers University\\
 New Brunswick, NJ 08854, USA}
\email{shelah@math.huji.ac.il}
\urladdr{http://shelah.logic.at}
\thanks{The first author would like to thank the Hebrew University of
  Jerusalem and the Lady Davis Fellowship Trust for awarding him a
  {\em Sch\"onbrunn Visiting Professorship\/} under which this research was
  carried out. \\
Both authors acknowledge support from the United States-Israel
Binational Science Foundation (Grant no. 2002323). This is publication
890 of the second author.}     

\date{September 2010}

\begin{abstract}
  We continue investigations of {\em reasonable ultrafilters\/} on
  uncountable cardinals defined in Shelah \cite{Sh:830}. We introduce
  stronger properties of ultrafilters and we show that those properties may
  be handled in $\lambda$--support iterations of reasonably bounding forcing
  notions. We use this to show that consistently there are reasonable
  ultrafilters on an inaccessible cardinal $\lambda$ with generating systems 
  of size less than $2^\lambda$. We also show how ultrafilters generated by
  small systems can be killed by forcing notions which have enough
  reasonable completeness to be iterated with $\lambda$--supports.
\end{abstract}

\maketitle

\section{Introduction}
{\em Reasonable ultrafilters\/} were introduced in Shelah \cite{Sh:830} in
order to suggest a line of research that would repeat in some sense the
beautiful theory created around the notion of {\em P--points on
$\omega$}. Most of the generalizations of P--points to uncountable
cardinals in the literature go into the direction of normal ultrafilters
and large cardinals (see, e.g., Gitik \cite{Gi81}), but one may be
interested in the opposite direction. If one wants to keep away from {\em
normal ultrafilters on $\lambda$}, one may declare interest in
ultrafilters which do not include some clubs and even demand that quotients 
by a closed unbounded subset of $\lambda$ do not extend the club filter of
$\lambda$. Such ultrafilters are called {\em weakly reasonable
ultrafilters}, see \ref{1.5}, \ref{easyob}. But if we are interested in
generalizing P--points, we have to consider also properties that would
correspond to {\em any countable family of members of the ultrafilter has a
pseudo-intersection in the ultrafilter}. The choice of the right property
in the declared context of {\em very non-normal ultrafilters\/} is not
clear, and one of the goals of the present paper is to show that the {\em
very reasonable ultrafilters\/} suggested in Shelah \cite{Sh:830} (see
Definition \ref{1.3} here) are very reasonable indeed, that is we may prove
interesting theorems on them.     

In the first section we recall some of the concepts and results presented in
Shelah \cite{Sh:830} and we introduce strong properties of generating
systems (super and strong reasonability, see Definitions \ref{superplus},
\ref{supersemi}) and we show that there may exist super reasonable systems
which generate ultrafilters (Propositions \ref{1.5x}, \ref{forcesuper}).  

In the next section we recall from \cite{RoSh:860} some properties of
forcing notions relevant for $\lambda$--support iterations. We also improve
in some sense a result of \cite{RoSh:860} and we show a preservation theorem
for {\em the nice double {\bf a}--bounding property} (Theorem
\ref{presdouble}). 

Then in the third section we show that super reasonable families generating
ultrafilters will be still at least strongly reasonable and will continue to
generate ultrafilters after forcing with $\lambda$--support iterations of
{\bf A}--bounding forcing notions. Therefore, for an inaccessible cardinal
$\lambda$, it is consistent that $2^\lambda=\lambda^{++}$ and there is a
very reasonable ultrafilter generated by a system of size $\lambda^+$
(Corollary \ref{p.4A}). It should be stressed that ``generating an
ultrafilter'' has the specific meaning stated in Definition \ref{1.3}(3). In
particular, ``having a small generating system'' {\em does not\/} imply
``having small ultrafilter base''.

The fourth section shows that some technical inconveniences of the proofs
from the third sections reflect the delicate nature of our concepts, not
necessarily our lack of knowledge. We give an example of a nicely double {\bf
  a}--bounding forcing notion which kills ultrafilters generated by systems
from the ground model. Then we show that for an inaccessible cardinal
$\lambda$, it is consistent that $2^\lambda=\lambda^{++}$ and there is no 
ultrafilter generated by a system of size $\lambda^+$ (see Corollary
\ref{p.4A}).

Studies of ultrafilters generated according to the schema introduced in
\cite{Sh:830} are also carried out in Ros{\l}anowski and Shelah
\cite{RoSh:889}.         
\medskip

\noindent {\bf Notation:}\quad Our notation is rather standard and
compatible with that of classical textbooks (like Jech \cite{J}). In forcing
we keep the older convention that {\em a stronger condition is the larger
  one}.  

\begin{enumerate}
\item Ordinal numbers will be denoted be the lower case initial letters of
the Greek alphabet ($\alpha,\beta,\gamma,\delta\ldots$) and also by $i,j$
(with possible sub- and superscripts). Cardinal numbers will be called
$\kappa,\lambda,\mu$ (with possible sub- and superscripts). {\bf $\lambda$
  is always assumed to be regular, sometimes even strongly inaccessible}.

By $\chi$ we will denote a {\em sufficiently large\/} regular cardinal; 
$\cH(\chi)$ is the family of all sets hereditarily of size less than
$\chi$. Moreover, we fix a well ordering $<^*_\chi$ of $\cH(\chi)$. 

\item A sequence is a function with the domain being a set of ordinals. For
  two sequences $\eta,\nu$ we write $\nu\vartriangleleft\eta$ whenever $\nu$
  is a proper initial segment of $\eta$, and $\nu \trianglelefteq\eta$ when
  either $\nu\vartriangleleft\eta$ or $\nu=\eta$.  The length of a sequence
  $\eta$ is the order type of its domain and it is denoted by $\lh(\eta)$. 

\item We will consider several games of two players. One player will be
called {\em Generic\/} or {\em Complete\/} or just {\em COM\/}, and we
will refer to this player as ``she''. Her opponent will be called {\em
Antigeneric\/} or {\em Incomplete} or just {\em INC\/} and will be
referred to as ``he''. 

\item For a forcing notion $\bbP$, all $\bbP$--names for objects in the
  extension via $\bbP$ will be denoted with a tilde below (e.g.,
  $\name{\tau}$, $\name{X}$). The canonical $\bbP$--name for the generic
  filter in $\bbP$ is called $\name{G}_\bbP$. The weakest element of $\bbP$
  will be denoted by $\emptyset_\bbP$ (and we will always assume that there
  is one, and that there is no other condition equivalent to it). We will
  also assume that all forcing notions under consideration are atomless.

  By ``$\lambda$--support iterations'' we mean iterations in which domains
  of conditions are of size $\leq\lambda$. However, we will pretend that
  conditions in a $\lambda$--support iteration
  $\bar{\bbQ}=\langle\bbP_\zeta, \name{\bbQ}_\zeta:\zeta<\zeta^*\rangle$ are
  total functions on $\zeta^*$ and for $p\in\lim(\bar{\bbQ})$ and
  $\alpha\in\zeta^*\setminus\Dom(p)$ we will let
  $p(\alpha)=\name{\emptyset}_{\name{\bbQ}_\alpha}$.

\item For a filter $D$ on $\lambda$, the family of all $D$--positive subsets
  of $\lambda$ is called $D^+$. (So $A\in D^+$ if and only if $A\subseteq
  \lambda$ and $A\cap B\neq\emptyset$ for all $B\in D$.)

  The club filter of $\lambda$ is denoted by $\cD_\lambda$.
\end{enumerate}

\section{More reasonable ultrafilters on $\lambda$}
Here we recall some basic definitions and results from \cite{Sh:830}, and
then we introduce even stronger properties of ultrafilters and/or generating
systems. We also show that assumptions like
$\diamondsuit_{S^{\lambda^+}_\lambda}$ imply the existence of such objects.  

As explained in the introduction, we are interested in ultrafilters (on an
uncountable cardinal $\lambda$) which are far from being normal. Weakly
reasonable ultrafilters defined below do not contain some clubs even if we
look at their quotients by a club. 

\begin{definition}
[{\cite[Def. 1.4]{Sh:830}}]
\label{1.5}
We say that a uniform ultrafilter $D$ on $\lambda$ is {\em weakly
  reasonable\/} if for every function $f\in {}^\lambda\lambda$ there is a
club $C$ of $\lambda$ such that 
\[\bigcup\{[\delta,\delta + f(\delta)):\delta \in C\}\notin D.\]
\end{definition}

\begin{observation}
[{\cite[Obs. 1.5]{Sh:830}}]
\label{easyob}
Let $D$ be a uniform ultrafilter on $\lambda$. Then the following conditions
are equivalent: 
\begin{enumerate}
\item[(A)] $D$ is weakly reasonable,
\item[(B)] for every increasing continuous sequence $\langle\delta_\xi:\xi< 
\lambda\rangle\subseteq\lambda$ there is a club $C^*$ of $\lambda$ such that 
\[\bigcup\big\{[\delta_\xi,\delta_{\xi+1}):\xi\in C^*\big\}\notin D.\] 
\end{enumerate}
\end{observation}

We want to investigate ultrafilters on $\lambda$ which are generated by
systems defining ``largeness in $\lambda$'' by giving a condition based on
``largeness in intervals below $\lambda$''. The family $\bqz$ introduced
below is a natural generalization of the approach used in \cite[Sections 5,
6]{RoSh:470}. The directness of $G^*$ is an easy way to guarantee that
$\fil(G^*)$ is a filter, and $({<}\lambda^+)$--directness has the flavour of
$P$--pointness.  

\begin{definition}
[{\cite[Def. 2.5]{Sh:830}}]
\label{1.3} 
\begin{enumerate}
\item Let $\bqz$ consist of all tuples 
\[p=(C^p,\langle Z^p_\delta:\delta\in C^p\rangle,\langle d^p_\delta:\delta
\in C^p \rangle)\]  
such that  
\begin{enumerate}
\item[(i)]   $C^p$ is a club of $\lambda$ consisting of limit ordinals only,
  and for $\delta\in C^p$: 
\item[(ii)]  $Z^p_\delta=\big[\delta,\min\big(C^p\setminus(\delta+1)\big)
\big)$ and 
\item[(iii)] $d^p_\delta\subseteq\cP(Z^p_\delta)$ is a proper non-principal
  ultrafilter on $Z^p_\delta$. 
\end{enumerate}
\item For $q\in\bqz$ we let 
\[\fil(q)\stackrel{\rm def}{=}\big\{A\subseteq\lambda: (\exists\vare<
\lambda)(\forall\delta\in C^q\setminus\vare)(A\cap Z^q_\delta\in
d^q_\delta)\big\},\]  
and for a set $G^*\subseteq\bqz$ we let $\fil(G^*)\stackrel{\rm def}{=}
\bigcup\{\fil(p):p\in G^*\}$. We also define a binary relation $\leq^0$ on
$\bqz$ by 
\begin{center}
$p\leq^0 q$\quad if and only if\quad $\fil(p)\subseteq\fil(q)$.
\end{center}
\item We say that an ultrafilter $D$ on $\lambda$ is {\em reasonable\/} if
it is weakly reasonable (see \ref{1.5}) and there is a directed (with
respect to $\leq^0$) set $G^*\subseteq\bqz$ such that $D=\fil(G^*)$. The
family $G^*$ may be called {\em the generating system for $D$}.     
\item An ultrafilter $D$ on $\lambda$ is said to be {\em very reasonable} if 
it is weakly reasonable and there is a $({<}\lambda^+)$--directed (with
respect to $\leq^0$) set $G^*\subseteq\bqz$ such that $D=\fil(G^*)$. 
\end{enumerate}
\end{definition}

\begin{definition}
\label{eproduct} 
Suppose that
\begin{enumerate}
\item[(a)] $X$ is a non-empty set and $e$ is an ultrafilter on $X$,
\item[(b)] $d_x$ is an ultrafilter on a set $Z_x$ (for $x\in X$). 
\end{enumerate}
We let 
\[\bigoplus\limits^e_{x\in X} d_x=\big\{A\subseteq\bigcup\limits_{x\in X}
Z_x: \{x\in X:Z_x\cap A\in d_x\}\in e\big\}.\]
(Clearly, $\bigoplus\limits^e_{x\in X} d_x$ is an ultrafilter on
$\bigcup\limits_{x\in X} Z_x$.)
\end{definition}

\begin{proposition}
[{\cite[Prop. 2.9]{Sh:830}}]
\label{1.3C}
Let $p,q\in\bqz$. Then the following are equivalent:
\begin{enumerate}
\item[(a)]  $p\leq^0 q$,
\item[(b)] there is $\vare<\lambda$ such that  
\[\big(\forall\alpha\in C^q\setminus\vare\big)\big(\forall A\in d^q_\alpha 
\big)\big(\exists \beta\in C^p\big)\big(A\cap Z^p_\beta\in d^p_\beta
\big),\]    
\item[(c)] there is $\vare<\lambda$ such that\\ 
if $\alpha\in C^q\setminus\vare$, $\beta_0=\sup\big(C^p\cap(\alpha+1)\big)$, 
$\beta_1=\min\big( C^p\setminus\min(C^q\setminus(\alpha+1))\big)$,\\
then there is an ultrafilter $e$ on $[\beta_0,\beta_1)\cap C^p$ such that 
\[d^q_\alpha=\big\{A\cap Z^q_\alpha:A\in\bigoplus\limits^e\{d^p_\beta:
\beta\in [\beta_0,\beta_1)\cap C^p\}\big\}.\] 
\end{enumerate}
\end{proposition}

\begin{observation}
[Compare {\cite[Prop. 2.3(4)]{Sh:830}}]
\label{forult}
If $p\in\bqz$, $A\subseteq\lambda$, then there is $q\in\bqz$ such that
$p\leq^0 q$ and either $A\in\fil(q)$ or $\lambda\setminus A\in\fil(q)$.
\end{observation}

\begin{definition}
[{\cite[Def. 2.10]{Sh:830}}]
\label{restrictable}
Let $p\in\bqz$. Suppose that $X\in [C^p]^{\textstyle\lambda}$ and
$C\subseteq C^p$ is a club of $\lambda$ such that
\begin{quotation}
if $\alpha<\beta$ are successive elements of $C$,\\
then $|[\alpha,\beta)\cap X|=1$.
\end{quotation}
(In this situation we say that {\em $p$ is restrictable to $\langle X,C
\rangle$}.) We define {\em the restriction of $p$ to $\langle X,C
\rangle$} as an element $q=p\rest\langle X,C\rangle\in \bqz$ such that 
$C^q=C$, and if $\alpha<\beta$ are successive elements of $C$, $x\in
[\alpha,\beta)\cap X$, then $Z^q_\alpha=[\alpha,\beta)$ and $d^q_\alpha=
\{A\subseteq Z^q_\alpha: A\cap Z^p_x\in d^p_x\}$. 
\end{definition}
 
\begin{proposition}
[{\cite[Prop. 2.11]{Sh:830}}]
\label{2.10}
\begin{enumerate}
\item  If $G^*\subseteq\bqz$ is $\leq^0$--directed and $|G^*|\leq\lambda$,
  then $G^*$ has a $\leq^0$--upper bound. (Hence, in particular, $\fil(G^*)$
  is not an ultrafilter.) 
\item Assume that $G^*\subseteq\bqz$ is $\leq^0$--directed and
  $\leq^0$--downward closed, $p\in G^*$, $X\in [C^p]^{\textstyle\lambda}$
  and $C\subseteq C^p$ is a club of $\lambda$ such that $p$ is restrictable
  to $\langle X,C\rangle$. If $\bigcup\limits_{x\in X}Z_x^p\in\fil(G^*)$,
  then $p\rest\langle X,C\rangle\in G^*$. 
\end{enumerate}
\end{proposition}

The following definition is used here to simplify our notation in
\ref{superplus} only. However, these concepts play a more central role in
\cite{RoSh:889}.  

\begin{definition}
\label{genfil} 
\begin{enumerate}
\item Let $\bqd$ be the family of all sets $r$ such that  
\begin{enumerate}
\item[(a)] members of $r$ are triples $(\alpha,Z,d)$ such that
$\alpha<\lambda$, $Z\subseteq [\alpha,\lambda)$, $\aleph_0\leq |Z|<\lambda$
  and $d$ is a non-principal ultrafilter on $Z$, and  
\item[(b)] $\big(\forall\xi<\lambda\big)\big(|\{(\alpha,Z,d)\in r:\alpha=
  \xi\}|<\lambda\big)$, and $|r|=\lambda$.
\end{enumerate}
For $r\in\bqd$ we define
\[\fil^*(r)=\big\{A\subseteq\lambda:\big(\exists\vare<\lambda\big)\big(\forall 
(\alpha,Z,d)\in r\big)\big(\vare\leq\alpha\ \Rightarrow\ A \cap Z\in d\big)
\big\},\] 
and we define a binary relation $\leq^*$ on $\bqd$ by 

$r_1\leq^* r_2$ if and only if ($r_1,r_2\in\bqd$ and)
$\fil^*(r_1)\subseteq \fil^*(r_2)$.
\item For a set $G_*\subseteq\bqd$ we let $\fil^*(G_*)=\bigcup\big\{
  \fil^*(r): r\in G_*\}$.   
\item We say that an $r\in\bqd$ is {\em strongly disjoint\/} if and only if 
\begin{itemize}
\item $\big(\forall\xi<\lambda\big)\big(|\{(\alpha,Z,d)\in r:\alpha=\xi\}|   
<2\big)$, and 
\item $\big(\forall (\alpha_1,Z_1,d_1),(\alpha_2,Z_2,d_2)\in r\big)\big(
\alpha_1<\alpha_2\ \Rightarrow\ Z_1\subseteq\alpha_2\big)$.
\end{itemize}
\item For $p\in\bqz$ we let $\#(p)=\{(\alpha,Z^p_\alpha,d^p_\alpha):
  \alpha\in C^p\}$. 
\end{enumerate}
\end{definition}

\begin{observation}
\label{qstar}
\begin{enumerate}
\item If $p\in\bqz$ then $\#(p)\in\bqd$ is strongly disjoint and
  $\fil(p)=\fil^*(\#(p))$.  Also, if $r\in\bqd$ is strongly disjoint, then
  $\fil^*(r)=\fil(p)$ for some $p\in\bqz$.
\item Let $r,s\in\bqd$. Then $r\leq^* s$\quad if and only if\quad there is
  $\vare<\lambda$ such that   
\[\big(\forall(\alpha,Z,d)\in s\big)\big(\forall A\in
d\big)\big(\alpha>\vare\ \Rightarrow\ \big(\exists (\alpha',Z',d')\in r\big)
\big(A\cap Z'\in d'\big)\big).\]    
\end{enumerate}
\end{observation}

The various definitions of super reasonable ultrafilters introduced in
Definition \ref{superplus} below are motivated by the proof of ``the Sacks
forcing preserves $P$--points''. In that proof, a fusion sequence is
constructed so that at a stage $n<\omega$ of the construction one deals with
{\em finitely many\/} nodes in a condition (the nodes that are declared to
be kept). We would like to carry out this kind of argument, e.g., for
forcing notions used in \cite[B.8.3, B.8.5]{RoSh:777}, but now we have to
deal with $<\lambda$ nodes  in a tree, and the ultrafilter we try to
preserve is not that complete. So what do we do? We deal with {\em
  finitely\/} many nodes at a time eventually taking care of everybody. One
can think that in the definition below the set $I_\alpha$ is the set of
nodes we have to keep and the finite sets $u_{\alpha,i}$ are the nodes taken
care of at a substage $i$.  

The technical aspects of \ref{superplus} are motivated by the iteration
theorems in \cite{RoSh:860} and \cite{RoSh:888}: our games here are taylored
to fit the games played on trees of conditions in $\lambda$--support
iterations, see Theorems \ref{lemult}, \ref{lemsup} later. As said earlier,
the main goal is to have a property of $G^*$ which implies the preservation
of ``$\fil(G^*)$ is an ultrafilter'' by many forcing notions. We would also
love to preserve that property itself, but we failed to achieve it. The
``super reasonability'' is what we need to preserve the ultrafilter (see
\ref{lemult}), ``strong reasonability'' is what we can prove about $G^*$ in
the extension (see \ref{lemsup}).

\begin{definition}
\label{superplus}
Let $G^*\subseteq\bqz$ and let $\bar{\mu}=\langle\mu_\alpha:\alpha<
\lambda\rangle$ be a sequence of cardinals, $2\leq\mu_\alpha\leq\lambda$ for
$\alpha<\lambda$.  
\begin{enumerate}
\item We define a game $\Game^\boxplus_{\bar{\mu}}(G^*)$ between two
players, COM and INC. A play of $\Game^\boxplus_{\bar{\mu}}(G^*)$ lasts
$\lambda$ steps and at a stage $\alpha<\lambda$ of the play the players
choose $I_\alpha,i_\alpha,\bar{u}_\alpha$ and $\langle r_{\alpha,i},
r'_{\alpha,i},(\beta_{\alpha,i},Z_{\alpha,i},d_{\alpha,i}):i<i_\alpha\rangle$
applying the following procedure.  
\begin{itemize}
\item First, INC chooses a non-empty set $I_\alpha$ of cardinality
$<\mu_\alpha$ and an enumeration $\bar{u}_\alpha=\langle u_{\alpha,i}:
i<i_\alpha\rangle$ of $[I_\alpha]^{<\omega}$ (so $i_\alpha<\mu_\alpha\cdot
\aleph_0$).       
\item Next the two players play a subgame of length $i_\alpha$. In the
$i^{\rm th}$ move of the subgame, 
\begin{enumerate}
\item COM chooses $r_{\alpha,i}\in G^*$, and then 
\item INC chooses $r'_{\alpha,i}\in G^*$ such that $r_{\alpha,i}\leq^0
  r'_{\alpha,i}$, and finally  
\item COM picks $(\beta_{\alpha,i},Z_{\alpha,i},d_{\alpha,i})\in
  \#\big(r_{\alpha,i}'\big)$ such that $\beta_{\alpha,i}>\alpha$. 
\end{enumerate}
\end{itemize}
In the end of the play COM wins if and only if
\begin{enumerate}
\item[$(\boxplus)$]  there is $r\in G^*$ such that for every $\bar{j}=
  \langle j_\alpha:\alpha<\lambda\rangle\in \prod\limits_{\alpha<\lambda}
  I_\alpha$ we have  
\[\{(\beta_{\alpha,i},Z_{\alpha,i},d_{\alpha,i}):\alpha<\lambda,\ j_\alpha
\in u_{\alpha,i} \mbox{ and }i<i_\alpha\}\leq^* \#(r).\] 
\end{enumerate}
A game $\Game^\boxminus_{\bar{\mu}}(G^*)$ is defined similarly to
$\Game^\boxplus_{\bar{\mu}}(G^*)$ except that $(\boxplus)$ is weakened to   
\begin{enumerate}
\item[$(\boxminus)$] for every $\bar{j}\in\prod\limits_{\alpha<\lambda} 
I_\alpha$ the set $\bigcup\{Z_{\alpha,i}:\alpha<\lambda,\ i<i_\alpha$ and
$j_\alpha \in u_{\alpha,i}\}$ belongs to $\fil(G^*)$.
\end{enumerate}
\item We say that the family $G^*$ is {\em $\bar{\mu}$--super reasonable\/}
  ({\em $\bar{\mu}$--super$^-$ reasonable,\/} respectively) if 
  \begin{enumerate}
\item[(i)] $G^*$ is $({<}\lambda^+)$--directed (with respect to $\leq^0$),
    and 
\item[(ii)] if $s\in G^*$, $r\in\bqz$ and for some $\alpha<\lambda$ we have
  $C^r=C^s\setminus \alpha$ and $d^r_\beta=d^s_\beta$ for $\beta\in C^r$,
  then $r\in G^*$, and  
\item[(iii)]  INC has no winning strategy in the game
  $\Game^\boxplus_{\bar{\mu}}(G^*)$ ($\Game^\boxminus_{\bar{\mu}}(G^*)$,
  respectively). 
  \end{enumerate}
\item We say that a uniform ultrafilter $D$ on $\lambda$ is
  {\em $\bar{\mu}$--super reasonable\/} ({\em $\bar{\mu}$--super$^-$
    reasonable}, respectively) if there is a $\bar{\mu}$--super reasonable  
  ($\bar{\mu}$--super$^-$ reasonable, respectively) set $G^*\subseteq\bqz$
  such that $D=\fil(G^*)$. 
\item If $\mu_\alpha=\lambda$ for all $\alpha<\lambda$, then we omit
  $\bar{\mu}$ and say just {\em super reasonable\/} or {\em super$^-$
    reasonable\/} (in reference to both ultrafilters on $\lambda$ and
  families $G^*\subseteq\bqz$). Also in this case we may write
  $\Game^\boxplus$ instead of $\Game^\boxplus_{\bar{\mu}}$. 
\end{enumerate}
\end{definition}

\begin{definition}
\label{supersemi}
Let $G^*\subseteq\bqz$ be directed with respect to $\leq^0$ and let
$\bar{\mu}=\langle\mu_\alpha:\alpha<\lambda\rangle$ be a sequence of
cardinals, $2\leq\mu_\alpha\leq\lambda$ for $\alpha<\lambda$. 
\begin{enumerate}
\item A game $\Game^\oplus_{\bar{\mu}}(G^*)$ between two players, COM and
INC is defined as follows. A play of $\Game^\oplus_{\bar{\mu}}(G^*)$ lasts
$\lambda$ steps and at a stage $\alpha<\lambda$ of the play the players
choose $I_\alpha,i_\alpha,\bar{u}_\alpha$ and $\langle r_{\alpha,i},
\delta_{\alpha,i},(\beta_{\alpha,i},Z_{\alpha,i},d_{\alpha,i}):i<i_\alpha\rangle$ 
applying the following procedure.  
\begin{itemize}
\item First, INC chooses a non-empty set $I_\alpha$ of cardinality
$<\mu_\alpha$, and then COM chooses $i_\alpha<\lambda$ and a sequence
$\bar{u}_\alpha=\langle u_{\alpha,i}:i<i_\alpha\rangle$ of non-empty finite
subsets of $I_\alpha$ such that $I_\alpha=\bigcup\limits_{i<i_\alpha}
u_{\alpha,i}$.    
\item Next the two players play a subgame of length $i_\alpha$. In the
$i^{\rm th}$ move of the subgame, 
\begin{enumerate}
\item COM chooses $r_{\alpha,i}\in G^*$ and then 
\item INC chooses $\delta_{\alpha,i}<\lambda$, and finally 
\item COM picks $(\beta_{\alpha,i},Z_{\alpha,i},d_{\alpha,i})\in
  \#\big(r_{\alpha,i}\big)$ such that $\beta_{\alpha,i}$ is above
  $\delta_{\alpha,i}$ and $\alpha$.   
\end{enumerate}
\end{itemize}
In the end of the play COM wins if and only if
\begin{enumerate}
\item[$(\oplus)$]  there is $r\in G^*$ such that for every $\bar{j}= \langle
  j_\alpha:\alpha<\lambda\rangle\in\prod\limits_{\alpha<\lambda} I_\alpha$
  we have  
\[\{(\beta_{\alpha,i},Z_{\alpha,i},d_{\alpha,i}):\alpha<\lambda,\ j_\alpha
\in u_{\alpha,i} \mbox{ and }i<i_\alpha\}\leq^* \#(r).\] 
\end{enumerate}
A game $\Game^\ominus_{\bar{\mu}}(G^*)$ is defined similarly to
$\Game^\oplus_{\bar{\mu}}(G^*)$ except that $(\oplus)$ is weakened to   
\begin{enumerate}
\item[$(\ominus)$] for every $\bar{j}\in\prod\limits_{\alpha<\lambda} 
I_\alpha$ the set $\bigcup\{Z_{\alpha,i}:\alpha<\lambda,\ i<i_\alpha$ and
$j_\alpha \in u_{\alpha,i}\}$ belongs to $\fil(G^*)$.
\end{enumerate}
\item If $G^*\subseteq\bqz$ is $({<}\lambda^+)$--directed (with respect to
$\leq^0$) and INC has no winning strategy in the game
$\Game^\oplus_{\bar{\mu}}(G^*)$, then we say that $G^*$ is {\em
$\bar{\mu}$--strongly reasonable}. Also, $G^*$ is said to be {\em
$\bar{\mu}$--strongly$^-$ reasonable\/} if it is  $({<}\lambda^+)$--directed
and INC has no winning strategy in the  game
$\Game^\ominus_{\bar{\mu}}(G^*)$.   
\item We say that a uniform ultrafilter $D$ on $\lambda$ is
  {\em $\bar{\mu}$--strongly reasonable\/} ({\em $\bar{\mu}$--strongly$^-$
    reasonable}, respectively) if there is a $\bar{\mu}$--strongly
  reasonable ($\bar{\mu}$--strongly$^-$ reasonable, respectively) set
  $G^*\subseteq\bqz$ such that $D=\fil(G^*)$. If $\mu_\alpha=\lambda$ for
  all $\alpha<\lambda$, then we omit $\bar{\mu}$ and say just {\em strongly
    reasonable\/} or {\em strongly$^-$ reasonable\/}. 
\end{enumerate}
\end{definition}

\begin{observation}
Assume that $2\leq\mu_\alpha\leq\kappa_\alpha\leq\lambda$ for
$\alpha<\lambda$ and $\bar{\mu}=\langle\mu_\alpha:\alpha<\lambda \rangle$,
$\bar{\kappa}=\langle\kappa_\alpha:\alpha<\lambda \rangle$. Then for a
family $G^*\subseteq\bqz$ and/or a uniform ultrafilter $D$ on $\lambda$ the 
following implications hold. 
\smallskip

%\hspace{-1cm}
\begin{tabular}{ccccc}
\small $\bar{\kappa}$--super reasonable $\Rightarrow$ &
\small $\bar{\mu}$--super reasonable $\Rightarrow$ &
\small $\bar{\mu}$--strongly reasonable\\
$\Downarrow$ & $\Downarrow$ & $\Downarrow$\\
\small $\bar{\kappa}$--super$^-$ reasonable $\Rightarrow$ &
\small $\bar{\mu}$--super$^-$ reasonable $\Rightarrow$ &
\small $\bar{\mu}$--strongly$^-$ reasonable 
\end{tabular}
\end{observation}

\begin{proposition}
\label{cd1.5x}
Assume that $2\leq\mu_\alpha\leq\lambda$ for $\alpha<\lambda$ and
$\bar{\mu}=\langle\mu_\alpha:\alpha<\lambda \rangle$. If a uniform
ultrafilter $D$ on $\lambda$ is $\bar{\mu}$--strongly$^-$ reasonable, then
it is very reasonable.  
\end{proposition}

\begin{proof} 
Pick a $\bar{\mu}$--strongly$^-$ reasonably family $G^*\subseteq\bqz$ such
that $D=\fil(G^*)$. Then $G^*$ is $({<}\lambda^+)$--directed and the proof 
will be completed once we show that $D$ is weakly reasonable. 

Let $f\in{}^\lambda\lambda$. We will argue that for some club $C=\{
\gamma_\alpha:\alpha<\gamma\}\subseteq\lambda$ we have
$\bigcup\{[\delta,\delta+f(\delta)):\delta\in C\}\notin D$, where
$\gamma_\alpha$ are given by the arguments below.    

We consider the following strategy $\st(f)$ for INC in
$\Game^\ominus_{\bar{\mu}} (G^*)$. The strategy $\st(f)$ instructs 
INC to construct on the side an increasing continuous sequence
$\langle\gamma_\alpha:\alpha<\lambda\rangle\subseteq \lambda$ so that at a
stage $\alpha<\lambda$ of the play, when
\[\big\langle I_\xi,i_\xi,\bar{u}_\xi,\langle r_{\xi,i},\delta_{\xi,i},
(\beta_{\xi,i},Z_{\xi,i},d_{\xi,i}):i<i_\xi\rangle:\xi<\alpha\big\rangle\]  
is the result of the play so far, then
\begin{itemize}
\item if $\alpha$ is limit, then $\gamma_\alpha=\sup(\gamma_\xi:\xi<
  \alpha)$,  
\item if $\alpha$ is not limit, then $\gamma_\alpha=\sup\big(\bigcup\{
Z_{\xi,i}: i<i_\xi,\ \xi<\alpha\}\big)+1$.
\end{itemize}
Now (at the stage $\alpha$) $\st(f)$ instructs INC to choose
$I_\alpha=\{0\}$ and then (after COM picks $i_\alpha,\bar{u}_\alpha$) he is
instructed to play in the subgame of this stage as follows. At stage
$i<i_\alpha$, after COM has picked $r_{\alpha,i}$, INC lets
\[\delta_{\alpha,i}=\gamma_\alpha+f(\gamma_\alpha)+\sup\big(\bigcup\{ 
Z_{\alpha,j}:j<i\}\big)+890.\] 
(After this COM chooses $(\beta_{\alpha,i},Z_{\alpha,i},d_{\alpha,i})\in
  \#\big(r_{\alpha,i}\big)$ with $\beta_{\alpha,i}>\delta_{\alpha,i}$.) 

The strategy $\st(f)$ cannot be the winning one for INC, so there is a play   
\[\big\langle I_\alpha,i_\alpha,\bar{u}_\alpha,\langle r_{\alpha,i},
\delta_{\alpha,i},(\beta_{\alpha,i},Z_{\alpha,i}, d_{\alpha,i}):
i<i_\alpha\rangle:\alpha<\lambda\big\rangle\]  
of $\Game^\ominus_{\bar{\mu}}(G^*)$ in which INC follows $\st(f)$ but 
\[A^*\stackrel{\rm def}{=}\bigcup\big\{Z_{\alpha,i}:\alpha<\lambda,\
i<i_\alpha\big\}\in\fil(G^*)=D\]
(note that necessarily $u_{\alpha,i}=I_\alpha=\{0\}$). It follows from the
choice of $\gamma_\alpha,\delta_{\alpha,i}$ that for each $\alpha<\lambda$
\[[\gamma_\alpha,\gamma_\alpha+f(\gamma_\alpha))\cap\bigcup\big\{
Z_{\xi,i}:\xi<\lambda,\ i<i_\xi\big\}=\emptyset,\]
and hence also $\bigcup\big\{[\gamma_\alpha,\gamma_\alpha+
f(\gamma_\alpha)):\alpha<\lambda\big\}\cap A^*=\emptyset$. Consequently
$\bigcup\big\{[\gamma_\alpha,\gamma_\alpha+f(\gamma_\alpha)): \alpha<
\lambda\big\}\notin D$ and one can easily finish the proof.  
\end{proof} 

\begin{proposition}
\label{1.5x} 
Assume $\lambda=\lambda^{<\lambda}$ and
$\diamondsuit_{S^{\lambda^+}_\lambda}$ holds. There exists a sequence
$\langle r_\xi:\xi<\lambda^+\rangle\subseteq\bqz$ such that  
\begin{enumerate}
\item[(i)]   $(\forall\xi<\zeta<\lambda^+)(r_\xi\leq^0 r_\zeta)$, and 
\item[(ii)]  the family 
\[G^*\stackrel{\rm def}{=} \big\{r\in\bqz:(\exists\xi<\lambda^+)(r\leq^0
r_\xi)\big\}\]  
is super reasonable and $\fil(G^*)$ is an ultrafilter on $\lambda$. 
\end{enumerate}
\end{proposition}

\begin{proof}
The sequence $\langle r_\xi:\xi<\lambda^+\rangle$ will be constructed
inductively. At successor stages we will use \ref{forult} to make sure that
$\fil(G^*)$ is an ultrafilter. At limit stages we will use \ref{2.10}(1) to
find upper bounds to the sequence constructed so far. Moreover, at (some)
stages $\xi$ of cofinality $\lambda$ the element $r_\xi$ will be chosen so
that ``it kills'' a strategy for INC in $\Game^\boxplus(G^*)$ predicted by
the diamond sequence.

For $\alpha<\lambda$ let $X^1_\alpha$ be the set of all legal plays of
$\Game^\boxplus(\bqz)$ of the form 
\begin{enumerate}
\item[$(\odot)^1_\alpha$] \quad $\big\langle I_\gamma,i_\gamma,
  \bar{u}_\gamma,\langle r_{\gamma,i},r_{\gamma,i}',(\beta_{\gamma,i},
  Z_{\gamma,i}, d_{\gamma,i}):i< i_\gamma\rangle:\gamma<\alpha\big\rangle$   
\end{enumerate}
where each $I_\gamma$ (for $\gamma<\alpha$) is an ordinal below
$\lambda$. Also let $X^1=\bigcup\limits_{\alpha<\lambda}X^1_\alpha$. Next,
for $\alpha<\lambda$, $0<I<\lambda$ and an enumeration $\bar{u}=\langle 
u_j:j<i\rangle$ of $[I]^{<\omega}$ let $X^2_{\alpha,I,\bar{u}}$ be the set
of all legal plays of $\Game^\boxplus(\bqz)$ of the form 
\begin{enumerate}
\item[$(\odot)^2_{\alpha,I,\bar{u}}$]\quad $\bar{\sigma}\conc \langle
  (I,i,\bar{u})\rangle\conc \langle r_j,r_j',(\beta_j,Z_j,d_j):j<j^*\rangle
  \conc \langle r\rangle$,
\end{enumerate}
where $\bar{\sigma}\in X^1_\alpha$, $j^*<i$ (and $\langle
r_j,r_j',(\beta_j,Z_j,d_j):j<j^*\rangle\conc \langle r\rangle$ is a legal
partial play of the subgame of level $\alpha$; in particular $r_j,r_j',r\in
\bqz$). Also let  
\[\begin{array}{ll}
X^2=\bigcup\big\{X^2_{\alpha,I,\bar{u}}:&\alpha<\lambda\ \mbox{ and }\
0<I<\lambda\ \mbox{ and}\\
&\bar{u}=\langle u_j:j<i\rangle\mbox{ is an enumeration  of }
[I]^{<\omega}\big\}.
\end{array}\] 
Any strategy for INC in $\Game^\boxplus(\bqz)$ can be interpreted as
a function $\st$ such that 
\begin{enumerate}
\item[$(\odot)^3$]  the domain of $\st$ is $X^1\cup X^2$,   
\item[$(\odot)^4$]  if $\bar{\sigma}\in X^1_\alpha$, $\alpha<\lambda$, then
  $\st(\bar{\sigma})=(I,i,\bar{u})$ for some $I<\lambda$ and an enumeration
  $\bar{u}=\langle u_j:j<i\rangle$ of $[I]^{<\omega}$, 
\item[$(\odot)^5$] if $\bar{\sigma}\in X^2_{\alpha,I,\bar{u}}$,
  $\alpha<\lambda$, $0<I<\lambda$, $\bar{u}=\langle u_j:j<i\rangle=
  [I]^{<\omega}$, and $\bar{\sigma}=\bar{\sigma}_0\conc\langle
  (I,i,\bar{u})\rangle\conc \langle r_j,r_j',(\beta_j,Z_j,d_j):j<j^*\rangle\conc
  \langle r\rangle$, then $\st(\bar{\sigma})\in\bqz$ is such that
  $r\leq^0\st(\bar{\sigma})$.  
\end{enumerate}
Below, whenever we say {\em a strategy for INC\/} we mean a function $\st$
satisfying conditions $(\odot)^3$--$(\odot)^5$.  

Since $|\bqz|=2^{2^{<\lambda}}=\lambda^+$, we may pick a bijection $\pi_0: 
\bqz\stackrel{1-1}{\longrightarrow}\lambda^+$ and for $\xi<\lambda^+$ let
$\cX_\xi$ consist of all $\bar{\sigma}\in X^1\cup X^2$ such that
$\pi_0(r)<\xi$ for all elements $r\in\bqz$ involved in the representation of
$\bar{\sigma}$ as in $(\odot)^1,(\odot)^2$. We also let $\cY_\xi$ consist of
all pairs $(\bar{\sigma},a)$ such that 
\begin{itemize}
\item $\bar{\sigma}\in\cX_\xi$ and $a=\st(\bar{\sigma})$ for some strategy
  $\st$ of INC, and 
\item if $\bar{\sigma}\in X^2$ (and so $a\in\bqz$) then $\pi_0(a)<\xi$. 
\end{itemize}
Note that $|\cX_\xi|\leq\lambda$ and $|\cY_\xi|\leq\lambda$  (for each
$\xi<\lambda^+$). Put $\cY=\bigcup\limits_{\xi<\lambda^+}\cY_\xi$. Plainly
$|\cY|=\lambda^+$ so we may fix a bijection $\pi_1:\lambda^+ \stackrel{\rm
  onto}{\longrightarrow}\cY$. Let 
\[C=\{\xi<\lambda^+:\pi_1[\xi]=\cY_\xi\};\]
it is a club of $\lambda^+$. 

Let $\langle A_\zeta:\zeta<\lambda^+\rangle$ list all subsets of $\lambda$
and let $\langle B_\zeta:\zeta\in S^{\lambda^+}_\lambda\rangle$ be a diamond
sequence on $S^{\lambda^+}_\lambda=\{\zeta<\lambda^+:\cf(\zeta)=\lambda\}$. 
By induction on $\xi<\lambda^+$ we choose a $\leq^0$--increasing sequence
$\langle r_\xi:\xi<\lambda^+\rangle \subseteq\bqz$ applying the following
procedure. Assume $\xi<\lambda^+$ and we have constructed $\langle
r_\zeta:\zeta<\xi\rangle$.   
\medskip

\noindent{\sc Case 0}:\quad $\xi=0$.\\
We let $r_0$ be the $<^*_\chi$--first member of $\bqz$.
\smallskip

\noindent{\sc Case 1}:\quad $\xi=\zeta+1$.\\ 
Pick $r_\xi\in\bqz$ such that $r_\zeta\leq^0 r_\xi$ and either $A_\zeta\in
\fil(r_\xi)$ or $\lambda\setminus A_\zeta\in\fil(r_\xi)$ (remember
Observation \ref{forult}).
\smallskip

\noindent{\sc Case 2}:\quad $\xi$ is a limit ordinal, $\cf(\xi)<\lambda$.\\  
Pick $r_\xi\in\bqz$ such that $(\forall\zeta<\xi)(r_\zeta\leq^0 r_\xi)$
(exists by Proposition \ref{2.10}(1)).
\smallskip

\noindent{\sc Case 3}:\quad $\xi$ is a limit ordinal, $\cf(\xi)=\lambda$.\\  
Now we ask if 
\begin{enumerate}
\item[$(\odot)^6_\xi$]\quad $\xi\in C$ and $(\forall\zeta<\xi)(\pi_0(
  r_\zeta)<\xi)$ and there is a strategy $\st$ for INC in $\Game^\boxplus
  (\bqz)$ such that $\pi_1[B_\xi]=\st\cap\cY_\xi=\st\rest\cX_\xi$. 
\end{enumerate}
If the answer to $(\odot)^6_\xi$ is negative, then we choose $r_\xi\in\bqz$
as in Case 2. 

Suppose now that the answer to $(\odot)^6_\xi$ is positive (so in
particular $\xi\in C$) and $\st$ is a strategy for INC such that
$\pi_1[B_\xi]= \st\cap\cY_\xi=\st\rest\cX_\xi$. Let $\bar{\xi}=\langle
\xi_\alpha:\alpha<\lambda\rangle$ be an increasing continuous sequence
cofinal in $\xi$. Consider a play 
\[\bar{\sigma}=\big\langle I_\alpha,i_\alpha,\bar{u}_\alpha,\langle
r_{\alpha,i},r_{\alpha,i}',(\beta_{\alpha,i},Z_{\alpha,i},d_{\alpha,i}):i<i_\alpha
\rangle: \alpha<\lambda\big\rangle\]
of $\Game^\boxplus(\bqz)$ in which INC follows the strategy $\st$ and COM 
proceeds as follows. When playing $\Game^\boxplus(\bqz)$, at step
$i<i_\alpha$ of the subgame of level $\alpha<\lambda$ (of
$\Game^\boxplus(\bqz)$) COM chooses $r_{\alpha,i}=r_{\xi_\alpha}$ and then,
after INC determines $r_{\alpha,i}'$ by $\st$, she picks the
$<^*_\chi$--first $(\beta_{\alpha,i},Z_{\alpha,i},d_{\alpha,i})\in\#(
r_{\alpha,i}')$ satisfying: 
\begin{enumerate}
\item[$(\odot)^7_{\xi,\alpha,i}$]\quad $(\forall\gamma\leq\alpha)(\forall
  A\in d_{\alpha,i})(\exists\delta\in C^{r_{\xi_\gamma}})(A\cap
  Z^{r_{\xi_\gamma}}_\delta\in d^{r_{\xi_\gamma}}_\delta)$ (remember
  \ref{1.3C}) and 
\item[$(\odot)^8_{\xi,\alpha,i}$]\quad  $(\forall\gamma<\alpha)(\forall
  j<i_\gamma)(Z_{\gamma,j}\subseteq\beta_{\alpha,i})$ and $(\forall j<i)(
  Z_{\alpha,j}\subseteq\beta_{\alpha,i})$. 
\end{enumerate}
The above rules fully determine the play $\bar{\sigma}$ and it should be
clear that $\bar{\sigma}\rest \alpha\in \cX_\xi$ for each $\alpha<\lambda$.
Note that $\bar{\sigma}$ depends on $B_\xi$ and $\bar{\xi}$ only (and not on
$\st$, provided it is as required by $(\odot)^6_\xi$). 

By the demands $(\odot)^8_{\xi,\alpha,i}$, we may choose an increasing
continuous sequence $\langle\gamma_\alpha:\alpha<\lambda\rangle
\subseteq\lambda$ such that $\gamma_0=0$ and $(\forall\alpha<\lambda)
(\forall i<i_\alpha)(Z_{\alpha,i}\subseteq[\gamma_\alpha,\gamma_{\alpha+1})
)$.  Now, for $\alpha<\lambda$ choose an ultrafilter $e_\alpha$ on
$i_\alpha$ such that 
\begin{enumerate}
\item[$(\odot)^9_{\xi,\alpha}$]\quad $\big(\forall j\in I_\alpha\big)\big(\{i<
  i_\alpha: j\in u_{\alpha,i}\}\in e_\alpha\big)$ 
\end{enumerate}
and let $d_\alpha$ be an ultrafilter on $[\gamma_\alpha,\gamma_{\alpha+1})$
such that 
\begin{enumerate}
\item[$(\odot)^{10}_{\xi,\alpha}$]\quad $\bigoplus\limits^{e_\alpha}\big
  \{d_{\alpha,i}:i<i_\alpha\big\}\subseteq d_\alpha$.  
\end{enumerate}
Now let $r_\xi\in\bqz$ be such that 
\begin{itemize}
\item $C^{r_\xi}=\{\gamma_\alpha:\alpha<\lambda\}$, and 
\item if $\delta=\gamma_\alpha$, then $Z^{r_\xi}_\delta=[\gamma_\alpha,
  \gamma_{\alpha+1})$ and $d^{r_\xi}_\delta=d_\alpha$.  
\end{itemize}
One easily verifies that $r_{\xi_\alpha}\leq^0 r_\xi$ for all
$\alpha<\lambda$ (remember $(\odot)^7$ and the choice of $d_\alpha$; use
\ref{1.3C}) and so $r_\zeta\leq^0 r_\xi$ for all $\zeta<\xi$. It follows
from $(\odot)^9_{\xi,\alpha}$ and $(\odot)^{10}_{\xi,\alpha}$ that 
\begin{enumerate}
\item[$(\odot)^{11}_\xi$] for every  $\bar{j}=\langle j_\alpha: \alpha<
  \lambda\rangle\in\prod\limits_{\alpha<\lambda} I_\alpha$ we have 
\[\big\{(\beta_{\alpha,i},Z_{\alpha,i},d_{\alpha,i}): \alpha< \lambda\ \&\
j_\alpha\in u_{\alpha,i}\ \&\ i< i_{\alpha,i}\big\}\leq^*\#(r_\xi).\]    
\end{enumerate}
After the construction of $\langle r_\xi:\xi<\lambda^+\rangle$ is carried
out we let 
\[G^*=\{r\in\bqz:(\exists \xi<\lambda^+)(r\leq^0 r_\xi)\}.\]
Plainly, $G^*$ satisfies demands (i) and (ii) of \ref{superplus}(2) and
$\fil(G^*)$ is an ultrafilter on $\lambda$ (remember Case 1 of the
construction). We should argue that INC has no winning strategy in
$\Game^\boxplus(G^*)$. To this end suppose that $\st^\boxplus$ is a strategy
of INC in $\Game^\boxplus(G^*)$. Pick $\xi\in  S^{\lambda^+}_\lambda\cap C$
such that $(\forall\zeta<\xi)(\pi_0(r_\zeta)<\xi)$ and $\pi_1[B_\xi]=
\st^\boxplus\cap\cY_\xi=\st^\boxplus\rest\cX_\xi$. Then when choosing
$r_\xi$ we gave a positive answer to $(\odot)^6_\xi$ and we constructed a
play $\bar{\sigma}$ of $\Game^\boxplus(\bqz)$. In that play, INC follows
$\st^\boxplus$ and COM chooses members of $G^*$, so it is a play of 
$\Game^\boxplus(G^*)$ . Now the condition $(\odot)^{11}_\xi$ means that
$r_\xi$ witnesses that COM wins the play $\bar{\sigma}$ and consequently
$\st^\boxplus$ is not a winning strategy for INC. 
\end{proof}

\begin{proposition}
\label{forcesuper}
Let $\bqz=(\bqz,\leq^0)$. 
\begin{enumerate}
\item $\bqz$ is a $({<}\lambda^+)$--complete forcing notion of size
  $2^{2^{<\lambda}}$.
\item $\forces_{\bqz}$`` $\name{G}_{\bqz}$ is a super reasonable family and
  $\fil(\name{G}_{\bqz})$ is an ultrafilter ''. 
\end{enumerate}
\end{proposition}

\begin{proof}
(1)\quad Should be clear; see also Proposition \ref{2.10}(1). 

\noindent(2)\quad By the completeness of $\bqz$, forcing with it does not
add new subsets of $\lambda$, and by \ref{1.3C}
\[\forces_{\bqz}\mbox{``  }\fil(\name{G}_{\bqz})\mbox{ is a uniform
  ultrafilter on $\lambda$ ''.}\]
It should also be clear that $\name{G}_{\bqz}$ satisfies the demands of
\ref{superplus}(2)(i+ii) (in $\bV^\bqz$). Let us argue that 
\[\forces_{\bqz}\mbox{``  INC has no winning strategy in }\Game^\boxplus
(\name{G}_{\bqz})\mbox{ ''}\]
and to this end suppose $p\in\bqz$ and $\name{\st}$ is a $\bqz$--name such
that 
\[p\forces_{\bqz}\mbox{``  $\name{\st}$ is a strategy of INC in }\Game^\boxplus
(\name{G}_{\bqz})\mbox{ ''}.\]
We are going to construct a condition $q\in\bqz$ stronger than $p$ and a
play $\bar{\sigma}$ of $\Game^\boxplus(\bqz)$ such that 
\[q\forces_{\bqz}\mbox{``  $\bar{\sigma}$ is a play of }\Game^\boxplus
(\name{G}_{\bqz})\mbox{ in which INC follows $\name{\st}$ but COM wins''}.\] 
Let $X^1,X^2$ be defined as in the proof of \ref{1.5x} (see
$(\odot)^1_\alpha$, $(\odot)^2_{\alpha,I,\bar{u}}$ there). We may assume
that  
\[p\forces_{\bqz}\mbox{``  $\name{\st}$ is a function satisfying
  $(\odot)^3$--$(\odot)^5$ of the proof of \ref{1.5x} ''.}\]
By induction on $\alpha<\lambda$ we choose conditions $p_\alpha\in\bqz$ and
partial plays $\bar{\sigma}_\alpha\in X^1_\alpha$ so that 
\begin{enumerate}
\item[$(\boxdot)_1$] $p\leq^0 p_\alpha\leq^0 p_\beta$ and
  $\bar{\sigma}_\alpha\vtl \bar{\sigma}_\beta$ for $\alpha<\beta<\lambda$, 
\item[$(\boxdot)_2$] $p_\alpha\forces_{\bqz}$`` $\bar{\sigma}_\alpha$ is a
  partial play of $\Game^\boxplus(\name{G}_{\bqz})$ in which INC uses
  $\name{\st}$ '', 
\item[$(\boxdot)_3$] if $\bar{\sigma}_\alpha=\big\langle I_\gamma,i_\gamma,
  \bar{u}_\gamma,\langle r_{\gamma,i},r_{\gamma,i}',(\beta_{\gamma,i},
  Z_{\gamma,i}, d_{\gamma,i}):i<i_\gamma\rangle:\gamma<\alpha\big\rangle$,
  then for every $\gamma<\delta<\alpha$ and $j<i<i_\gamma$ we have 
\[r_{\gamma,i}'\leq^0 p_\alpha\quad\mbox{ and }\quad Z_{\gamma,j}
\subseteq\beta_{\gamma,i}\quad\mbox{ and }\quad Z_{\gamma,j} \subseteq
\beta_{\delta,0}.\] 
\end{enumerate}
Suppose that $\alpha=\alpha^*+1$ and we have determined $p_{\alpha^*},
\bar{\sigma}_{\alpha^*}$. Pick $p_\alpha'\geq^0 p_{\alpha^*}$ and $I_\alpha,
i_\alpha, \bar{u}_\alpha$ such that $p_\alpha'\forces\name{\st}(
\bar{\sigma}_{\alpha^*})=(I_\alpha,i_\alpha,\bar{u}_\alpha)$. Now choose
inductively $p^i_\alpha,r_{\alpha,i},r_{\alpha,i}'$ and $(\beta_{\alpha,i},
Z_{\alpha,i}, d_{\alpha,i})$ for $i<i_\alpha$ so that for each
$i<j<i_\alpha$ we have 
\begin{enumerate}
\item[$(\boxdot)_4$]
  \begin{enumerate}
  \item[(i)]   $p^0_\alpha=p_\alpha'$, $p^i_\alpha\leq^0 p^j_\alpha$,
    $p^i_\alpha=r_{\alpha,i}\leq^0 r_{\alpha,i}'\leq^0 p^{i+1}_\alpha$, and  
  \item[(ii)]  $p^{i+1}_\alpha\forces$`` $r_{\alpha,i}'$ is the answer by
    $\name{\st}$ at stage $i$ of the subgame '', 
  \item[(iii)]  $\beta_{\alpha,i}$ satisfies the demand in $(\boxdot)_3$ and
    $(\beta_{\alpha,i},Z_{\alpha,i},d_{\alpha,i})\in\#(r_{\alpha,i}')$, 
   \item[(iv)]  $(\forall A\in d_{\alpha,i})(\forall\gamma\leq\alpha^*) 
     (\exists\delta\in C^{p_\gamma})(A\cap Z^{p_\gamma}_\delta\in
     d^{p_\gamma}_\delta)$. 
  \end{enumerate}
\end{enumerate}
Then $p_{\alpha+1}$ is any $\leq^0$--upper bound to $\{p^i_\alpha:
i<i_\alpha\}$.  

The limit stages of the construction should be clear.

After the construction is carried out and we have $\bar{\sigma}_\lambda=
\bigcup\{\bar{\sigma}_\alpha:\alpha<\lambda\}$, we define $r\in\bqz$ like
$r_\xi$ in the proof of \ref{1.5x} (see $(\odot)^9_{\xi,\alpha}+
(\odot)^{10}_\xi$  there). Then $r$ is $\leq^0$--stronger then all
$p_\alpha$ (for $\alpha<\lambda$) and 
\[r\forces_{\bqz}\mbox{`` $\bar{\sigma}_\lambda$ is a play of
  $\Game^\boxplus(\name{G}_{\bqz})$ in which INC uses $\name{\st}$ but COM
  wins ''.}\]
(Note that the respective version of $(\odot)^{11}_\xi$ of the proof of
\ref{1.5x} holds. By the completeness it continues to hold in $\bV^{\bqz}$.) 
\end{proof}

\section{More on reasonably complete forcing}

\begin{definition}
\label{strcom}
Let $\bbP$ be a forcing notion.
\begin{enumerate}
\item For a condition $r\in\bbP$ let $\Game_0^\lambda(\bbP,r)$ be the
following game of two players, {\em Complete} and  {\em Incomplete}:   
\begin{quotation}
\noindent the game lasts at most $\lambda$ moves and during a play the
players attempt construct a sequence $\langle (p_i,q_i): i<\lambda\rangle$
of pairs of conditions from $\bbP$ in such a way that $(\forall
j<i<\lambda)(r\leq p_j\leq q_j\leq p_i)$ and at the stage $i<\lambda$ of the
game, first  Incomplete chooses $p_i$ and then Complete chooses $q_i$.  
\end{quotation}
Complete wins if and only if for every $i<\lambda$ there are legal moves for
both players. 
\item We say that the forcing notion $\bbP$ is {\em strategically
$({<}\lambda)$--complete\/} if Complete has a winning strategy in the game 
$\Game_0^\lambda(\bbP,p)$ for each condition $p\in\bbP$. 
\item Let $N\prec (\cH(\chi),\in,<^*_\chi)$ be a model such that
${}^{<\lambda} N\subseteq N$, $|N|=\lambda$ and $\bbP\in N$. We say that a
condition $p\in\bbP$ is {\em $(N,\bbP)$--generic in the standard
sense\/} (or just: {\em $(N,\bbP)$--generic\/}) if for every
$\bbP$--name $\name{\tau}\in N$ for an ordinal we have $p\forces$``
$\name{\tau}\in N$ ''. 
\item $\bbP$ is {\em $\lambda$--proper in the standard sense\/} (or just:
{\em $\lambda$--proper\/}) if there is $x\in \cH(\chi)$  such that for
every model $N\prec (\cH(\chi),\in,<^*_\chi)$ satisfying  
\[{}^{<\lambda} N\subseteq N,\quad |N|=\lambda\quad\mbox{ and }\quad\bbP,x  
  \in N, \]
and every condition $p\in N\cap\bbP$ there is an $(N,\bbP)$--generic
condition $q\in\bbP$ stronger than $p$.
\end{enumerate}
\end{definition}

\begin{theorem}
[See Shelah {\cite[Ch. III, Thm 4.1]{Sh:f}}, Abraham {\cite[\S 2]{Ab}} and
  Eisworth {\cite[\S 3]{Ei03}}] 
\label{lppcc}
Assume $2^\lambda=\lambda^+$, $\lambda^{<\lambda}=\lambda$. Let $\bar{\bbQ}=
\langle\bbP_i,\name{\bbQ}_i:i<\lambda^{++}\rangle$ be $\lambda$--support
iteration such that for all $i<\lambda^{++}$ we have  
\begin{itemize}
\item $\bbP_i$ is $\lambda$--proper,
\item $\forces_{\bbP_i}$`` $|\name{\bbQ}_i|\leq\lambda^+$ ''. 
\end{itemize}
Then 
\begin{enumerate}
\item for every $\delta<\lambda^{++}$,
  $\forces_{\bbP_\delta}2^\lambda=\lambda^+$, and
\item the limit $\bbP_{\lambda^{++}}$ satisfies the $\lambda^{++}$--cc.
\end{enumerate}

\end{theorem}

\begin{proposition}
[{\cite[Prop. A.1.6]{RoSh:777}}]
\label{pA.6}
Suppose $\bar{\bbQ}=\langle\bbP_i,\name{\bbQ}_i: i<\gamma\rangle$ is a
$\lambda$--support iteration and, for each $i<\gamma$, 
\[\forces_{\bbP_i}\mbox{`` $\name{\bbQ}_i$ is strategically
$({<}\lambda)$--complete ''.}\]
Then, for each $\vare\leq\gamma$ and $r\in\bbP_\vare$, there
is a winning strategy $\st(\vare,r)$ of Complete in the game
$\Game_0^\lambda(\bbP_\vare,r)$ such that, whenever $\vare_0<\vare_1\leq
\gamma$ and $r\in\bbP_{\vare_1}$, we have:
\begin{enumerate}
\item[(i)]  if $\langle (p_i,q_i):i<\lambda\rangle$ is a play of
$\Game_0^\lambda(\bbP_{\vare_0},r\rest\vare_0)$ in which Complete follows
the strategy $\st(\vare_0,r\rest\vare_0)$, then $\langle (p_i\conc r\rest
[\vare_0,\vare_1),q_i\conc r\rest [\vare_0,\vare_1)):i<\lambda\rangle$ is a
play of $\Game_0^\lambda(\bbP_{\vare_1},r)$ in which Complete uses
$\st(\vare_1,r)$;  
\item[(ii)] if $\langle (p_i,q_i):i<\lambda\rangle$ is a play of
$\Game_0^\lambda(\bbP_{\vare_1},r)$ in which Complete plays according to the
strategy $\st(\vare_1,r)$, then $\langle (p_i\rest\vare_0,q_i\rest\vare_0):i
<\lambda\rangle$ is a play of $\Game_0^\lambda(\bbP_{\vare_0},r\rest
\vare_0)$ in which Complete uses $\st(\vare_0,r\rest\vare_0)$; 
\item[(iii)] if $\langle (p_i,q_i):i<i^*\rangle$ is a partial play of
$\Game_0^\lambda(\bbP_{\vare_1},r)$ in which Complete uses $\st(\vare_1,r)$
and $p'\in\bbP_{\vare_0}$ is stronger than all $p_i\rest\vare_0$ (for
$i<i^*$), then there is $p^*\in \bbP_{\vare_1}$ such that
$p'=p^*\rest\vare_0$ and $p^*\geq p_i$ for $i<i^*$. 
\end{enumerate}
\end{proposition}

\begin{definition}
[Compare {\cite[Def. 2.2]{RoSh:860}}]
\label{dA.5}
\begin{enumerate}
\item Let $\gamma$ be an ordinal, $w\subseteq\gamma$. {\em A standard
    $(w,1)^\gamma$--tree\/} is a pair $\cT=(T,\rk)$ such that  
\begin{itemize}
\item $\rk:T\longrightarrow w\cup\{\gamma\}$, 
\item if $t\in T$ and $\rk(t)=\vare$, then $t$ is a sequence $\langle
(t)_\zeta: \zeta\in w\cap\vare\rangle$, 
\item $(T,\vtl)$ is a tree with root $\langle\rangle$ and such that every
chain in $T$ has a $\vartriangleleft$--upper bound in $T$, 
\item if $t\in T$, then there is $t'\in T$ such that $t\trianglelefteq t'$
  and $\rk(t')=\gamma$.
\end{itemize}
We will keep the convention that $\cT^x_y$ is $(T^x_y,\rk^x_y)$.
\item Let $\bar{\bbQ}=\langle\bbP_i,\name{\bbQ}_i:i<\gamma\rangle$ be a
$\lambda$--support iteration. {\em A standard tree of conditions in
$\bar{\bbQ}$\/} is a system $\bar{p}=\langle p_t:t\in T\rangle$ such that 
\begin{itemize}
\item $(T,\rk)$ is a standard $(w,1)^\gamma$--tree for some $w\subseteq
\gamma$, and 
\item $p_t\in\bbP_{\rk(t)}$ for $t\in T$, and
\item if $s,t\in T$, $s\vtl t$, then $p_s=p_t\rest\rk(s)$. 
\end{itemize}
\item Let $\bar{p}^0,\bar{p}^1$ be standard trees of conditions in
  $\bar{\bbQ}$, $\bar{p}^i=\langle p^i_t:t\in T\rangle$. We write
  $\bar{p}^0\leq \bar{p}^1$ whenever for each $t\in T$ we have $p^0_t\leq
  p^1_t$.  
\end{enumerate}
\end{definition}

Note that our standard trees and trees of conditions are a special case of
that introduced in \cite[Def. A.1.7]{RoSh:777} when $\alpha=1$. Also, the
rank function $\rk$ is essentially the function giving the level of a node,
adjusted to have values in $w\cup\{\gamma\}$ via the canonical increasing
bijection. 

\begin{proposition}
  [See {\cite[Prop. A.1.9]{RoSh:777}}]
\label{dectree}
Assume that $\bar{\bbQ}=\langle\bbP_i,\name{\bbQ}_i:i<\gamma\rangle$ is a 
$\lambda$--support iteration such that for all $i<\gamma$ we have 
\[\forces_{\bbP_i}\mbox{`` $\name{\bbQ}_i$ is strategically
$({<}\lambda)$--complete ''.}\]
Suppose that $\bar{p}=\langle p_t:t\in T\rangle$ is a standard tree of
conditions in $\bar{\bbQ}$, $|T|<\lambda$, and $\cI\subseteq\bbP_\gamma$
is open dense. Then there is a standard tree of conditions $\bar{q}=\langle
q_t:t\in T\rangle$ such that $\bar{p}\leq \bar{q}$ and $(\forall t\in
T)(\rk(t)=\gamma\ \Rightarrow\ q_t\in\cI)$, and such that conditions
$q_{t_0},q_{t_1}$ are incompatible whenever $t_0,t_1\in T$, $\rk(t_0)=
\rk(t_1)$ but $t_0\neq t_1$.
\end{proposition}

\begin{definition}
[See {\cite[Def. 3.1]{RoSh:860}}]
\label{p.1A}
Let $\bbQ$ be a forcing notion and let $\bar{\mu}=\langle\mu_\alpha:
\alpha<\lambda\rangle$ be a sequence of regular cardinals such that
$\aleph_0\leq\mu_\alpha\leq\lambda$ for all $\alpha<\lambda$.
\begin{enumerate}
\item For a condition $p\in\bbQ$ we define {\em a {\bf r}easonable A--{\bf 
c}ompleteness game\/} $\Agame(p,\bbQ)$ between two players, Generic and
Antigeneric, as follows. A play of $\Agame(p,\bbQ)$ lasts $\lambda$ steps
and during a play a sequence        
\[\Big\langle I_\alpha,\langle p^\alpha_t,q^\alpha_t:t\in I_\alpha\rangle:
\alpha<\lambda\Big\rangle\]
is constructed. Suppose that the players have arrived to a stage $\alpha<
\lambda$ of the game. Now, 
\begin{enumerate}
\item[$(\aleph)_\alpha$]  first Generic chooses a non-empty set $I_\alpha$
of cardinality $<\mu_\alpha$ and a system $\langle p^\alpha_t:t\in I_\alpha
\rangle$ of conditions from $\bbQ$,
\item[$(\beth)_\alpha$]  then Antigeneric answers by picking a system 
$\langle q^\alpha_t:t\in I_\alpha\rangle$ of conditions from $\bbQ$ such that 
$(\forall t\in I_\alpha)(p^\alpha_t\leq q^\alpha_t)$. 
\end{enumerate}
At the end, Generic wins the play 
\[\Big\langle I_\alpha,\langle p^\alpha_t,q^\alpha_t:t\in I_\alpha\rangle:
\alpha<\lambda\Big\rangle\]  
of $\Agame(p,\bbQ)$ \quad if and only if 
\begin{enumerate}
\item[$(\circledast)^{\rm rc}_{\rm A}$] there is a condition $p^*\in\bbQ$
stronger than $p$ and such that  
\[p^*\forces_{\bbQ}\mbox{`` }\big(\forall\alpha<\lambda\big)\big(\exists
t\in I_\alpha\big)\big(q^\alpha_t\in\name{G}_{\bbQ}\big)\mbox{ ''}.\]
\end{enumerate}
\item We say that a forcing notion $\bbQ$ is {\em reasonably A--bounding 
over $\bar{\mu}$\/} if    
\begin{enumerate}
\item[(a)] $\bbQ$ is strategically $({<}\lambda)$--complete,  and 
\item[(b)] for any $p\in\bbQ$, Generic has a winning strategy in the game
$\Agame(p,\bbQ)$. 
\end{enumerate}
\end{enumerate}
\end{definition}

\begin{definition}
[See {\cite[Def. 3.2]{RoSh:860}}]
\label{betterA}
Let $\bar{\bbQ}=\langle\bbP_\xi,\name{\bbQ}_\xi:\xi<\gamma\rangle$ be a
$\lambda$--support iteration and let $\bar{\mu}=\langle\mu_\alpha:
\alpha<\lambda\rangle$ be a sequence of regular cardinals such that
$\aleph_0\leq\mu_\alpha\leq\lambda$ for all $\alpha<\lambda$.  
  \begin{enumerate}
\item For a condition $p\in\bbP_\gamma=\lim(\bar{\bbQ})$ we define {\em a
tree A--completeness game\/} $\agame(p,\bar{\bbQ})$ between two players,
Generic and Antigeneric, as follows. A play of $\agame(p,\bar{\bbQ})$ lasts
$\lambda$ steps and in the course of a play a sequence $\langle\cT_\alpha,
\bar{p}^\alpha,\bar{q}^\alpha:  \alpha<\lambda\rangle$ is
constructed. Suppose that the players have arrived to a stage
$\alpha<\lambda$ of the game. Now, 
\begin{enumerate}
\item[$(\aleph)_\alpha$] first Generic picks a standard
$(w,1)^\gamma$--tree $\cT_\alpha$ such that $|T_\alpha|<\mu_\alpha$ and a
tree of conditions $\bar{p}^\alpha=\langle p^\alpha_t:t\in T_\alpha\rangle
\subseteq\bbP_\gamma$ (so Generic, as a part of choosing $\cT_\alpha$, picks
also $w=w_\alpha$), 
\item[$(\beth)_\alpha$] then Antigeneric answers by choosing a tree of
conditions 

$\bar{q}^\alpha=\langle q^\alpha_t:t\in T_\alpha\rangle\subseteq
\bbP_\gamma$ such that $\bar{p}^\alpha\leq\bar{q}^\alpha$.   
\end{enumerate}
At the end, Generic wins the play $\langle\cT_\alpha,\bar{p}^\alpha,
\bar{q}^\alpha:  \alpha<\lambda\rangle$ of $\agame(p,\bar{\bbQ})$ if and
only if 
\begin{enumerate}
\item[$(\circledast)^{\rm tree}_{\bf A}$] there is a condition
  $p^*\in\bbP_\gamma$ stronger than $p$ and such that  
\[p^*\forces_{\bbP_\gamma}\mbox{`` }\big(\forall\alpha< \lambda\big)
\big(\exists t\in T_\alpha\big)\big(\rk_\alpha(t)=\gamma\ \&\
q^\alpha_t\in \name{G}_{\bbP_\gamma}\big)\mbox{ ''}\] 
\end{enumerate}
\item We say that $\bbP_\gamma=\lim(\bar{\bbQ})$ is {\em reasonably$^*$ 
$A(\bar{\bbQ})$--bounding over $\bar{\mu}$} if Generic has a winning
strategy in the game $\agame(p,\bar{\bbQ})$ for every $p\in\bbP_\gamma$.   
\end{enumerate}
\end{definition}

\begin{theorem}
[See {\cite[Thm 3.2]{RoSh:860}}]
\label{verB}
Assume that 
\begin{enumerate}
\item[(a)] $\lambda$ is a strongly inaccessible cardinal,
\item[(b)] $\bar{\mu}=\langle\mu_\alpha:\alpha<\lambda\rangle$, each
$\mu_\alpha$ is a regular cardinal satisfying (for $\alpha<\lambda$)
\[\aleph_0\leq\mu_\alpha\leq\lambda\qquad\mbox{ and }\qquad \big(\forall f\in
{}^\alpha \mu_\alpha\big)\big(\big|\prod_{\xi<\alpha} f(\xi)\big|<
\mu_\alpha\big),\]
\item[(c)] $\bar{\bbQ}=\langle\bbP_\xi,\name{\bbQ}_\xi:\xi<\gamma\rangle$ is a
$\lambda$--support iteration such that for every $\xi<\gamma$,  
\[\forces_{\bbP_\xi}\mbox{`` $\name{\bbQ}_\xi$ is reasonably A--bounding
  over $\bar{\mu}$ ''.}\] 
\end{enumerate}
Then $\bbP_{\gamma}=\lim(\bar{\bbQ})$ is reasonably$^*$ 
$A(\bar{\bbQ})$--bounding over $\bar{\mu}$ (and so $\bbP_\gamma$ is also
$\lambda$--proper).  
\end{theorem}

In \cite[\S 3]{RoSh:860}, in addition to A--reasonable completeness game we
considered its variant called {\em {\bf a}--reasonable completeness game.\/}
In that variant, at stage $\alpha<\lambda$ of the game the players played a
subgame to construct a sequence $\langle
p^\alpha_\xi,q^\alpha_\xi:\xi<i_\alpha \rangle$ (corresponding to $\langle
p^\alpha_t,q^\alpha_t:t\in I_\alpha \rangle$). In the following definition
we introduce a further modification of that game. In the new game, the
players will again play subgames, in some sense repeating several times the
subgames from the {\bf a}--reasonable completeness game. 

\begin{definition}
  \label{po27}
Let $\bbQ$ be a forcing notion and let $\bar{\mu}=\langle \mu_\alpha:
\alpha<\lambda\rangle$ be a sequence of cardinals such that $\aleph_0\leq
\mu_\alpha<\lambda$ for all $\alpha<\lambda$. Suppose also that $\cU$ is
a normal filter on $\lambda$.
\begin{enumerate}
\item For a condition $p\in\bbQ$ we define {\em a  reasonable double--{\bf
a}--completeness game\/} $\tagame(p,\bbQ)$ between Generic and Antigeneric
as follows. A play of $\tagame(p,\bbQ)$ lasts at most $\lambda$ steps and in
the course of the play the players try to construct a sequence  
  \begin{enumerate}
  \item[$(\boxtimes)$] \quad $\big\langle\xi_\alpha,\langle p^\alpha_\gamma,
    q^\alpha_\gamma: \gamma<\mu_\alpha\cdot \xi_\alpha\rangle:\alpha<
    \lambda\big\rangle$. 
  \end{enumerate}
(Here $\mu_\alpha$ is treated as an ordinal and $\mu_\alpha\cdot\xi_\alpha$
is the ordinal product of $\mu_\alpha$ and $\xi_\alpha$.) Suppose that the
players have arrived to a stage $\alpha<\lambda$ of the game. First,
Antigeneric picks a non-zero ordinal $\xi_\alpha<\lambda$. Then the two
players start a subgame of length $\mu_\alpha\cdot\xi_\alpha$ alternately
choosing the terms of the sequence $\langle p^\alpha_\gamma,q^\alpha_\gamma:
\gamma<\mu_\alpha\cdot\xi_\alpha\rangle$. At a stage $\gamma=\mu_\alpha
\cdot i+j$ (where $i<\xi_\alpha$, $j<\mu_\alpha$) of the subgame, first
Generic picks a condition $p^\alpha_\gamma\in\bbQ$ stronger than all
conditions $q^\alpha_\delta$ for $\delta<\gamma$ of the form $\delta=
\mu_\alpha\cdot i'+j$ (where $i'<i$), and then Antigeneric answers with a
condition $q^\alpha_\gamma$ stronger than $p^\alpha_\gamma$.  

At the end, Generic wins the play $(\boxtimes)$ of $\tagame(p,\bbQ)$ if and
only if both players had always legal moves and 
\begin{enumerate}
\item[$(\circledast)^{\rm rc}_{\bf 2a}$] there is a condition $p^*\in\bbQ$
  stronger than $p$ and such that
\[p^*\forces_{\bbQ}\mbox{`` }\big(\forall\alpha<\lambda\big)\big(\exists
j<\mu_\alpha\big)\big(\{q^\alpha_{\mu_\alpha\cdot i+j}:i< \xi_\alpha \}
\subseteq \name{G}_{\bbQ}\big)\mbox{ ''}.\]
\end{enumerate}
\item Games $\bagame(p,\bbQ)$ (for $p\in\bbQ$) are defined similarly, we
  only replace condition $(\circledast)^{\rm rc}_{\bf 2a}$ by 
\begin{enumerate}
\item[$(\circledast)^{\rm rc}_{\bf 2b}$] there is a condition $p^*\in\bbQ$
  stronger than $p$ and such that
\[p^*\forces_{\bbQ}\mbox{`` }\big\{\alpha<\lambda:\big(\exists j<
  \mu_\alpha\big)\big(\{q^\alpha_{\mu_\alpha\cdot i+j}:i< \xi_\alpha \} 
\subseteq \name{G}_{\bbQ}\big)\big\}\in \cU^{\bbQ}\mbox{ ''},\]
where $\cU^{\bbQ}$ is the ($\bbQ$--name for the) normal filter generated by
$\cU$ in $\bV^{\bbQ}$. 
\end{enumerate}
\item A strategy $\st$ for Generic in $\tagame(p,\bbQ)$ (or
  $\bagame(p,\bbQ)$) is said to be {\em nice\/} if for every play
  $\big\langle\xi_\alpha,\langle p^\alpha_\gamma, q^\alpha_\gamma:
  \gamma<\mu_\alpha\cdot \xi_\alpha\rangle:\alpha< \lambda \big\rangle$ in
  which she uses $\st$, for every $\alpha<\lambda$, the conditions in
  $\{p^\alpha_\gamma:\gamma<\mu_\alpha\}$ are pairwise incompatible.  (These
  are conditions played in the first ``run'' of the subgame. Note that then
  $p^\alpha_\gamma,p^\alpha_{\gamma'}$ are incompatible whenever $\gamma
  \not\equiv \gamma'\mod \mu_\alpha$.) 
\item Let ${\bf x}\in \{{\bf a},{\bf b}\}$. A forcing notion $\bbQ$ is {\em 
    nicely double {\bf x}--bounding over $\bar{\mu}$} (and $\cU$ if ${\bf
    x}={\bf b}$) if 
  \begin{enumerate}
  \item[(a)] $\bbQ$ is strategically $({<}\lambda)$--complete, and 
  \item[(b)] Generic has a nice winning strategy in the game
    $\tagame(p,\bbQ)$ ($\bagame(p,\bbQ)$ if ${\bf x}={\bf b}$) for every
    $p\in\bbQ$.    
  \end{enumerate}
\end{enumerate}
\end{definition}

\begin{remark}
\begin{enumerate}
\item Reasonable double ${\bf x}$--boundedness (for ${\bf x}\in \{{\bf
    a},{\bf b}\}$) is an iterable relative of reasonable ${\bf
    x}$--boundedness introduced in \cite[Definition 3.1, pp
  206-207]{RoSh:860}. Technical differences in the definitions of suitable
  games are to achieve the preservation of the corresponding property in
  $\lambda$--support iterations (see Theorems \ref{presdouble},
  \ref{seconddouble} below). 
\item The game $\bagame(p,\bbQ)$ is easier to win for Generic than
  $\tagame(p,\bbQ)$ (because the winning criterion is weaker). Therefore, if
  we are interested in $\lambda$--properness for $\lambda$--support
  iterations only, then \ref{seconddouble} will cover a larger class of
  forcing notions than \ref{presdouble}.
\end{enumerate}
\end{remark}

\begin{definition}
[See {\cite[Def. 6.1]{RoSh:860}}] 
\label{silver}
Suppose that $\lambda$ is inaccessible and $\bar{\kappa}=\langle 
\kappa_\alpha:\alpha<\lambda\rangle$ is a sequence of cardinals,
$1<\kappa_\alpha<\lambda$ for $\alpha<\lambda$. We define a forcing notion
$\bbP^{\bar{\kappa}}$ as follows.\\ 
{\bf A condition in $\bbP^{\bar{\kappa}}$} is a pair $p=(f^p,C^p)$ such that  
\[C^p\subseteq\lambda\mbox{ is a club of $\lambda$ and }f^p\in\prod
\{\kappa_\iota:\iota\in\lambda\setminus C^p\}.\]
{\bf The order $\leq_{\bbP^{\bar{\kappa}}}=\leq$ of $\bbP^{\bar{\kappa}}$}
is given by:\\
$p\leq_{\bbP^{\bar{\kappa}}} q$\quad if and only if\quad $C^q\subseteq C^p$
and $f^p\subseteq f^q$.  
\end{definition}

\begin{proposition}
\begin{enumerate}
\item Assume that $\bar{\kappa},\lambda$ are as in \ref{silver} above and
  let a sequence $\bar{\mu}=\langle\mu_\alpha:\alpha<\lambda\rangle$ be
  chosen so that $\prod\limits_{\beta<\alpha}\kappa_\beta\leq\mu_\alpha<\lambda$
  (for $\alpha<\lambda$). Then the forcing notion $\bbP^{\bar{\kappa}}$ is
  nicely double {\bf b}--bounding over $\bar{\mu},\cD_\lambda$. 
\item If $\kappa_\alpha=\kappa$ for all $\alpha<\lambda$ and
  $\mu_\alpha\geq\kappa^\alpha$, then $\bbP^{\kappa}$ is nicely double {\bf
    a}--bounding over $\bar{\mu}$. 
\end{enumerate}
\end{proposition}

\begin{proof}
(1)\quad A natural modification of the proof of \cite[Prop. 6.1]{RoSh:860}
works here. Note that if $\bar{\delta}=\langle\delta_\alpha:\alpha< \lambda
\rangle$ is an increasing continuous sequence constructed as there during a
play of $\Game^{{\rm rc}{\bf 2b}}_{\bar{\mu},\cD_\lambda}(p,
\bbP^{\bar{\kappa}})$, then the set $B\stackrel{\rm
  def}{=}\big\{\alpha<\lambda: \prod\limits_{\beta<\alpha}
\kappa_{\delta_\beta}\leq \mu_\alpha\big\}$ is in the filter
$\cD_\lambda$. In the game, the stages $\alpha\in\lambda\setminus B$ are
ignored and only those for $\alpha\in B$ are ``active''. Also, at each stage 
$\alpha$ we may create $\mu_\alpha$ ``not active'' steps at each run of the
subgame by picking an antichain of conditions incompatible with $p$.
\smallskip

\noindent (2)\quad Similar; we get double {\bf a}--bounding here as {\em at
  each stage\/} $\alpha<\lambda$ of the game we know that
$\prod\limits_{\beta<\alpha}\kappa_{\delta_\beta}=\kappa^\alpha\leq\mu_\alpha$
(so all steps are ``active'').  
\end{proof}

\begin{theorem}
  \label{presdouble}
Assume that 
\begin{enumerate}
\item[(a)] $\lambda$ is a strongly inaccessible cardinal,
\item[(b)] $\bar{\mu}=\langle\mu_\alpha:\alpha<\lambda\rangle$ is a sequence
of cardinals below $\lambda$ such that $(\forall\alpha<\lambda)(
\aleph_0\leq \mu_\alpha=\mu_\alpha^{|\alpha+1|})$, 
\item[(c)] $\bar{\bbQ}=\langle\bbP_\zeta,\name{\bbQ}_\zeta:\zeta<\zeta^*
\rangle$ is a $\lambda$--support iteration such that for every $\zeta
<\gamma$,   
\[\forces_{\bbP_\zeta}\mbox{`` $\name{\bbQ}_\zeta$ is nicely double {\bf
    a}--bounding over $\bar{\mu}$ ''.}\] 
\end{enumerate}
Then $\bbP_{\zeta^*}=\lim(\bar{\bbQ})$ is nicely double {\bf a}--bounding
over $\bar{\mu}$ (and so $\bbP_{\zeta^*}$ is also $\lambda$--proper).
\end{theorem}

\begin{proof}
Our arguments refine those presented in the proof of \cite[Theorem 3.2,
p. 217]{RoSh:860}, but the differences in the games involved eliminate the
use of trees of conditions. However, trees of conditions are implicitely
present here too. The tree at level $\delta$ of the argument is indexed by  
\[T_\delta=\bigcup\big\{\prod\limits_{\xi\in w_\delta\cap \zeta}\mu_\delta:
\zeta\in w_\delta\cup\{\zeta^*\}\big\}\]
and it is formed in part by conditions played in the game for various $t\in
\prod\limits_{\xi\in w_\delta}\mu_\delta=\{t\in T_\delta:\rk(t)=\zeta^*\}$;
note the coherence demand in $(\boxtimes)_7$.  

Let $p\in\bbP_{\zeta^*}$. We will describe a strategy $\st$ for Generic in
the game $\tagame(p,\bbP_{\zeta^*})$. The strategy $\st$ instructs Generic
to play the game $\tagame$ on each relevant coordinate $\zeta<\zeta^*$ using
her winning strategy $\name{\st}_\zeta$. At stage $\delta<\lambda$ Generic
will be concerned with coordinates $\zeta\in w_\delta$ for some set
$w_\delta$ of size $<\lambda$. If $\xi_\delta<\lambda$ is the ordinal put by
Antigeneric in the play of $\tagame(p,\bbP_{\zeta^*})$, then in the
simulated plays on coordinates $\zeta\in w_\delta$ Generic pretends that her
opponent put $\xi_\delta^*=\mu_\delta\cdot \xi_\delta$. The innings of the
two players, Generic and Antigeneric, in the subgame of level $\delta$ on a
coordinate $\zeta$ are 
\[\name{\bar{p}}_{\delta,\zeta}=\langle \name{p}^\gamma_{\delta,\zeta}:
\gamma<\mu_\delta\cdot \xi^*_\delta\rangle\ \mbox{ and }\ 
\name{\bar{q}}_{\delta,\zeta}=\langle \name{q}^\gamma_{\delta,\zeta}:
\gamma<\mu_\delta\cdot \xi^*_\delta\rangle,\] 
respectively. Generic's innings in the subgame of $\tagame(p,
\bbP_{\zeta^*})$ will be associated with sequences
$t\in\prod\limits_{\zeta\in w_\delta} \mu_\delta=\langle t^\delta_j: j<
\mu_\delta\rangle=\bar{t}^\delta$. The innings of the two players will be
$p_\vare^\delta,q^\delta_\vare$ (for $\vare<\mu_\delta\cdot \xi_\delta$) and
they will be related to what happens at coordinates $\zeta\in w_\delta$ as
follows. If $t=t^\delta_j$, $\zeta\in w_\delta$ and
$\beta=(t)_\zeta<\mu_\delta$, then in the subgame of $\tagame(p,
\bbP_{\zeta^*})$ of level $\delta$ at stages of the form
$\vare=\mu_\delta\cdot i+j$ we will have $p^\delta_\vare(\zeta)=
\name{p}^\gamma_{\delta,\zeta}$ and $q^\delta_\vare(\zeta)=
\name{q}^\gamma_{\delta,\zeta}$, where $\gamma=\mu_\delta\cdot
\vare+\beta$. 

To keep track of what happens at coordinates $\zeta\notin w_\delta$ Generic
will use conditions $r_\delta$. 

Let us note that the construction of $\st$ presented in detail below would
be somewhat simpler if we knew that all forcings $\name{\bbQ}_\zeta$ are
$({<}\lambda)$--complete (and not only strategically
$({<}\lambda)$--complete). Then $\name{\st}^0_\xi,  r^-_\delta$ and
$p^{\delta,*}_\vare$ could be eliminated as their role is to make sure that
some sequences of conditions (related to $r_\delta$ and/or $p^\delta_\vare$)
have upper bounds. However, many natural forcing notions tend to have
strategic completeness only (see \cite[Part B]{RoSh:777}).   

Let us formalize the ideas presented above. For each $\zeta<\zeta^*$ pick a
$\bbP_\zeta$--name $\name{\st}^0_\zeta$ such that 
\[\begin{array}{r}
\forces_{\bbP_\zeta}\mbox{`` }\name{\st}^0_\zeta\mbox{ is a winning strategy 
  for Complete in }\Game_0^\lambda\big(\name{\bbQ}_\zeta,
  \name{\emptyset}_{\name{\bbQ}_\zeta}\big)\mbox{ such that }\ \\
\mbox{ if Incomplete plays $\name{\emptyset}_{\name{\bbQ}_\zeta}$ then 
  Complete answers with $\name{\emptyset}_{\name{\bbQ}_\zeta}$ as well ''.}
  \end{array}\] 
In the course of a play of $\tagame(p,\bbP_{\zeta^*})$, at a stage
$\delta<\lambda$, Generic will be instructed to construct on the side
\begin{enumerate}
\item[$(\otimes)_\delta$] $w_\delta,\bar{t}^\delta,\xi^*_\delta$,
  $\name{\st}_\zeta$ (for $\zeta\in w_{\delta+1}\setminus w_\delta$), 
  $\name{\bar{p}}_{\delta,\zeta},\name{\bar{q}}_{\delta,\zeta}$,
  $p^{\delta,*}_\vare$ (for $\vare< \mu_\delta\cdot\xi_\delta$), and
  $r^-_\delta, r_\delta$. 
\end{enumerate}
These objects will be chosen so that if 
\[\big\langle\xi_\delta,\langle p^\delta_\gamma,q^\delta_\gamma:\gamma< 
\mu_\delta\cdot \xi_\delta\rangle:\delta<\lambda\big\rangle\] 
is a play of $\tagame(p,\bbP_{\zeta^*})$ in which Generic follows $\st$, and
the additional objects constructed at stage $\delta<\lambda$ are listed in 
$(\otimes)_\delta$, then the following conditions are satisfied (for each
$\delta<\lambda$).  
\begin{enumerate}
\item[$(\boxtimes)_1$] $r^-_\delta,r_\delta\in \bbP_{\zeta^*}$, $r_0^-(0)=
  r_0(0)=p(0)$, $w_\delta\subseteq\zeta^*$, $|w_\delta|=|\delta+1|$,
  $\bigcup\limits_{\alpha<\lambda}\Dom(r_\alpha)=\bigcup\limits_{\alpha<
    \lambda} w_\alpha$, $w_0=\{0\}$, $w_\delta\subseteq w_{\delta+1}$ and if
  $\delta$ is limit then $w_\delta=\bigcup\limits_{\alpha<\delta} w_\alpha$.
\item[$(\boxtimes)_2$] For each $\alpha<\delta<\lambda$ we have
  $(\forall\zeta\in w_{\alpha+1})(r_\alpha(\zeta)= r^-_\delta(\zeta)=
  r_\delta(\zeta))$ and $p\leq r_\alpha^-\leq r_\alpha\leq r^-_\delta\leq
  r_\delta$, and $p^{\delta,*}_\vare\in \bbP_{\zeta^*}$ (for $\vare<
  \mu_\delta\cdot \xi_\delta$). 
\item[$(\boxtimes)_3$] If $\zeta\in\zeta^*\setminus w_\delta$, then 
\[\begin{array}{ll}
r_\delta\rest\zeta\forces_{\bbP_\zeta}&\mbox{`` the sequence } \langle
r^-_\alpha(\zeta), r_\alpha(\zeta):\alpha\leq\delta\rangle\mbox{ is a legal
  partial play of }\\ 
&\quad\Game_0^\lambda\big(\name{\bbQ}_\zeta,
\name{\emptyset}_{\name{\bbQ}_\zeta}\big)\mbox{ in which Complete follows
}\name{\st}^0_\zeta\mbox{ ''}  
\end{array}\]
and if $\zeta\in w_{\delta+1}\setminus w_\delta$, then $\name{\st}_\zeta$ is
a $\bbP_\zeta$--name  for a nice winning strategy for Generic in
$\tagame(r_\delta(\zeta),\name{\bbQ}_\zeta)$. (And $\st_0$ is a nice winning
strategy of Generic in $\tagame(p(0),\bbQ_0)$.)     
\item[$(\boxtimes)_4$] $\bar{t}^\delta=\langle t^\delta_j:j<\mu_\delta
\rangle$ is an enumeration of $\prod\limits_{\zeta\in w_\delta}\mu_\delta=
{}^{w_\delta}\mu_\delta$. 
\item[$( \boxtimes)_5$] $\xi^*_\delta=\mu_\delta\cdot \xi_\delta$ (the
ordinal product) and $\name{\bar{p}}_{\delta,\zeta}=\langle
\name{p}^\gamma_{\delta,\zeta}: \gamma<\mu_\delta\cdot\xi^*_\delta\rangle$
and $\name{\bar{q}}_{\delta,\zeta}=\langle\name{q}^\gamma_{\delta,\zeta}:
\gamma<\mu_\delta\cdot\xi^*_\delta\rangle$ are $\bbP_\zeta$--names for
sequences of conditions in $\name{\bbQ}_\zeta$ of length $\mu_\delta\cdot
\xi^*_\delta$ (for $\zeta\in\bigcup\limits_{\alpha<\lambda}w_\alpha$).    
\item[$( \boxtimes)_6$] If $\zeta\in w_{\beta+1}\setminus w_\beta$,
  $\beta<\delta$ (or $\zeta=\beta=0$),  then 
\[\begin{array}{r}
\forces_{\bbP_\zeta}\mbox{`` }\langle\xi^*_\alpha,\langle
\name{p}^\gamma_{\alpha,\zeta},\name{q}^\gamma_{\alpha,\zeta}: \gamma<
\mu_\alpha\cdot \xi^*_\alpha\rangle:\alpha\leq\delta\rangle\mbox{ is a
partial play of }\ \ \\
\tagame(r_\beta(\zeta), \name{\bbQ}_\zeta)\mbox{ in which Generic uses
  $\name{\st}_\zeta$ ''.} 
  \end{array}\]
\item[$( \boxtimes)_7$] If $\vare=\mu_\delta\cdot i+j$, $i<\xi_\delta$,
  $j<\mu_\delta$, then 
\[\Dom(p^{\delta,*}_\vare)=\Dom(p^\delta_\vare)=w_\delta\cup\Dom(p)\cup  
\bigcup\limits_{\alpha<\delta}\Dom(r_\alpha)\cup \bigcup\limits_{\vare'< 
  \vare}\Dom(q^\delta_{\vare'}),\]   
and for each $\zeta\in w_\delta\cup\{\zeta^*\}$ the condition
$p^{\delta,*}_\vare\rest\zeta$ is an upper bound to 
\[\begin{array}{l}
\{p\rest \zeta\}\cup\{r_\alpha\rest\zeta:\alpha<\delta\}\cup\\
\{q^\delta_{\vare'}\rest \zeta:\vare'=\mu_\delta\cdot i'+j'<\vare\ \&\
i'<\xi_\delta\ \&\ j'<\mu_\delta\ \&\ t^\delta_{j'}\rest\zeta= t^\delta_j 
\rest\zeta\}.
\end{array}\]    
\item[$( \boxtimes)_8$] If $j<\mu_\delta$, $i<\xi_\delta$, $\zeta\in
  w_\delta$, $(t^\delta_j)_\zeta=\beta$ and $\vare=\mu_\delta\cdot i+j$, 
  $\gamma=\mu_\delta\cdot\vare+\beta$, then $p^{\delta,*}_\vare(\zeta)
  =p^\delta_\vare(\zeta)=\name{p}^\gamma_{\delta,\zeta}$ and
  $q^\delta_\vare\rest\zeta\forces_{\bbP_\zeta} q^\delta_\vare(\zeta)= 
  \name{q}^\gamma_{\delta,\zeta}$.  
\item[$( \boxtimes)_9$] If $\vare=\mu_\delta\cdot i+j$, $i<\xi_\delta$,
  $j<\mu_\delta$, $\zeta\in\zeta^*\setminus w_\delta$ and
  $t\in \prod\{\mu_\delta:\xi\in w_\delta\cap\zeta\}$, $t\trianglelefteq
  t^\delta_j$, then 
\[\begin{array}{l}
p^\delta_\vare\rest\zeta\forces_{\bbP_\zeta}\mbox{`` the sequence }\\
\quad\ \langle p^{\delta,*}_{\vare'}(\zeta),p^\delta_{\vare'}(\zeta):
{\vare'}=\mu_\delta\cdot i'+j'\leq\vare\ \&\ i'<\xi_\delta\ \&\
j'<\mu_\delta\ \&\ t\trianglelefteq t^\delta_{j'}\rangle\\   
\quad\ \mbox{ is a legal partial play of }\Game^\lambda_0(
\name{\bbQ}_\zeta,p(\zeta))\mbox{ in which Complete follows
  $\name{\st}^0_\zeta$ ''.}  
  \end{array}\]
\item[$( \boxtimes)_{10}$] $\Dom(r_\delta^-)=\Dom(r_\delta)= \bigcup\{ 
\Dom(q^\delta_\vare):\vare<\mu_\delta\cdot\xi_\delta\}$ and if $\zeta\in
\zeta^*\setminus w_\delta$,   $t\in \prod\{\mu_\delta:\xi\in
w_\delta\cap\zeta\}$, and $q\in\bbP_\zeta$, $q\geq r^-_\delta\rest \zeta$
and $q\geq q^\delta_\vare\rest\zeta$ whenever $\vare=\mu_\delta\cdot i+j$,
$i<\xi_\delta$, $j<\mu_\delta$ and $t\trianglelefteq t^\delta_j$, then  
\[\begin{array}{ll}
q\forces_{\bbP_\zeta}&\mbox{`` if the set }\\
&\qquad\{p(\zeta)\}\cup \{r_\alpha(\zeta):\alpha<\delta\}\cup\\
&\qquad\{q^\delta_\vare(\zeta):\vare=\mu_\delta\cdot i+j\ \&\ i<\xi_\delta\
\&\ j<\mu_\delta\ \&\ t\trianglelefteq t^\delta_j\}\\  
&\mbox{ has an upper bound in }  
\name{\bbQ}_\xi,\ \mbox{ then $r^-_\delta(\zeta)$ is such an upper bound,}\\
&\mbox{ otherwise }\\
&\ r^-_\delta(\zeta)\mbox{ is just an upper bound to $\{p(\zeta)\}
  \cup\{r_\alpha(\zeta):\alpha<\delta\}$ ''.}  
  \end{array}\]
\end{enumerate}
Assume that the two players arrived to stage $\delta$ of
$\tagame(p,\bbP_{\zeta^*})$ and 
\[\big\langle\xi_\alpha,\langle p^\alpha_\vare, q^\alpha_\vare:
\vare<\mu_\alpha\cdot \xi_\alpha\rangle:\alpha<\delta \big\rangle\]
is the play constructed so far, and that Generic followed $\st$ and
determined objects listed in $(\otimes)_\alpha$ (for $\alpha<\delta$) with
properties $(\boxtimes)_1$--$(\boxtimes)_{10}$.  
  
Below, whenever we say {\em Generic chooses $x$ such that\/}
we mean {\em Generic chooses the $<^*_\chi$--first $x$ such that\/}, etc. 

First, Generic uses her favorite bookkeeping device to determine $w_\delta$
so that the demands of $(\boxtimes)_1$ are satisfied (and that at the end we
will have $\bigcup\limits_{\alpha<\lambda}\Dom(r_\alpha)=
\bigcup\limits_{\alpha<\lambda} w_\alpha$). If $\beta<\delta$ and $\zeta\in
w_\beta$, then we already have $\name{\bar{p}}_{\alpha,\zeta},
\name{\bar{q}}_{\alpha,\zeta}$ for $\alpha<\delta$ (see $(\boxtimes)_6$),
but we have not yet defined those objects when $\delta=\delta_0+1$ and
$\zeta\in w_\delta\setminus w_{\delta_0}$.  So if $\delta=\delta_0+1$ and
$\zeta\in w_\delta\setminus w_{\delta_0}$ then let
$\name{\bar{p}}_{\alpha,\zeta}=\langle\name{p}^\gamma_{\alpha,\zeta}:
\gamma<\mu_\alpha\cdot\xi^*_\alpha\rangle$ and $\name{\bar{q}}_{\alpha,
\zeta}=\langle\name{q}^\gamma_{\alpha,\zeta}:\gamma<\mu_\alpha \cdot 
\xi^*_\alpha\rangle$ (for $\alpha<\delta$) be such that 
\[\begin{array}{r}
\forces_{\bbP_\zeta}\mbox{`` }\langle\xi^*_\alpha,\langle
\name{p}^\gamma_{\alpha,\zeta},\name{q}^\gamma_{\alpha,\zeta}: \gamma<
\mu_\alpha\cdot \xi^*_\alpha\rangle:\alpha<\delta\rangle\mbox{ is a
partial play of }\ \\
\tagame(r_{\delta_0}(\zeta), \name{\bbQ}_\zeta)\mbox{ in which Generic uses
  $\name{\st}_\zeta$ and }\ \\
\name{p}^\gamma_{\alpha,\zeta}=\name{q}^\gamma_{\alpha,\zeta}\mbox{ for all
$\alpha<\delta$, $\gamma<\mu_\alpha\cdot\xi^*_\alpha$ ''.}
  \end{array}\]
Condition $(\boxtimes)_4$ and our rule of taking ``the $<^*_\chi$--first''
determine the enumeration $\bar{t}^\delta=\langle t^\delta_j:j<\mu_\delta
\rangle$ of $\prod\limits_{\zeta\in w_\delta}\mu_\delta$. Now Antigeneric
picks $\xi_\delta$ and the two players start a subgame of length
$\mu_\delta\cdot\xi_\delta$. During the subgame Generic will simulate
subgames of level $\delta$ at coordinates $\zeta\in w_\delta$ pretending
that Antigeneric played $\xi^*_\delta=\mu_\delta\cdot \xi_\delta$
there. Each step in the subgame of $\tagame(p,\bbP_{\zeta^*})$ will
correspond to $\mu_\delta$ steps in the subgames of $\tagame(r_\beta(\zeta),
\name{\bbQ}_\zeta)$ (when $\zeta\in w_{\beta+1}\setminus w_\beta$,
$\beta<\delta$). So suppose that the two opponents have arrived to a stage
$\vare=\mu_\delta\cdot i+j$ of the subgame, $i<\xi_\delta$, $j<\mu_\delta$,
and assume also that Generic (playing according to $\st$) has already
defined $\name{p}^\gamma_{\delta,\zeta},\name{q}^\gamma_{\delta,\zeta}$ for
$\zeta\in w_\delta$, $\gamma<\mu_\delta\cdot\vare$ and
$p^{\delta,*}_{\vare'}$ for $\vare'<\vare$, so that the requirements
of $(\boxtimes)_6$--$(\boxtimes)_9$ are satisfied. Note that (by
$(\boxtimes)_7$--$(\boxtimes)_9$) 
\begin{enumerate}
\item[$(\circledast)$] if $\vare>\vare'=\mu_\delta\cdot i'+j'> \vare''=
  \mu_\delta\cdot i''+j''$, $\zeta\in w_\delta\cup\{\zeta^*\}$ and
  $t^\delta_{j'}\rest\zeta = t^\delta_{j''}\rest \zeta$, then
  $p^\delta_{\vare''} \rest\zeta\leq q^\delta_{\vare''}\rest \zeta\leq
  p^{\delta,*}_{\vare'} \rest\zeta\leq p^\delta_\vare\rest\zeta$. 
\end{enumerate}
For each $\zeta\in w_\delta$ and $\beta<(t^\delta_j)_\zeta$ let
$\name{p}^{\mu_\delta \cdot \vare+\beta}_{\delta,\zeta}
=\name{q}^{\mu_\delta\cdot\vare+\beta}_{\delta,\zeta}$ be
$\bbP_\zeta$--names for conditions in $\name{\bbQ}_\zeta$ such that (the
relevant part of) $(\boxtimes)_6$ holds. The same clause determines also
$\name{p}^{\mu_\delta\cdot \vare+ \beta}_{\delta,\zeta}$ for $\beta=
(t^\delta_j)_\zeta$, $\zeta\in w_\delta$. Then the requirements in
$(\boxtimes)_7+(\boxtimes)_8$ essentially describe what $p^{\delta,*}_\vare$
is. Note that the ``upper bound demands'' in $(\boxtimes)_7$ can be
satisfied because of $(\boxtimes)_9+(\boxtimes)_3$ and
$(\circledast)$ above. Next, Generic's inning $p^\delta_\vare$ in
$\tagame(p,\bbP_{\zeta^*})$ is chosen so that $\Dom(p^\delta_\vare)=
\Dom(p^{\delta,*}_\vare)$ and clauses $(\boxtimes)_8+(\boxtimes)_9$
hold. After this Antigeneric answers with a condition $q^\delta_\vare\geq
p^\delta_\vare$, and Generic picks for the construction on the side names
$\name{q}^{\mu_\delta\cdot\vare+ \beta}_{\delta,\zeta}$ for $\zeta\in
w_\delta$ and $\beta=(t^\delta_j)_\zeta$ by the demand in $(\boxtimes)_8$.
She also picks $\name{p}^{\mu_\delta\cdot\vare+ \beta}_{\delta,\zeta}=
\name{q}^{\mu_\delta\cdot\vare+ \beta}_{\delta,\zeta}$ for $\zeta\in
w_\delta$ and $(t^\delta_j)_\zeta< \beta<\mu_\delta$ so that $(\boxtimes)_6$
holds.   
\smallskip

This completes the description of what happens during the $\mu_\delta\cdot
\xi_\delta$ steps of the subgame. After the subgame is over and the
sequence $\langle p^\delta_\gamma,q^\delta_\gamma:\gamma<\mu_\delta\cdot
\xi_\delta\rangle$ is constructed, Generic chooses conditions
$r^-_\delta,r_\delta\in \bbP_{\zeta^*}$ by $(\boxtimes)_1$--$(\boxtimes)_3$
and $(\boxtimes)_{10}$. (Note: since $\name{\st}_\zeta$ are names
for nice strategies, if $\zeta\in\zeta^*\setminus w_\delta$, $i_0,i_1<
\xi_\delta$, $j_0,j_1<\mu_\delta$, $\vare_0=\mu_\delta\cdot i_0+j_0$,
$\vare_1=\mu_\delta\cdot i_1+j_1$,  $t_0,t_1\in\prod\{\mu_\delta:\xi\in
w_\delta\cap\zeta\}$, $t_0\trianglelefteq t^\delta_{j_0}$,
$t_1\trianglelefteq t^\delta_{j_1}$ and $t_0\neq t_1$, then the conditions
$q^\delta_{\vare_0}\rest\zeta,q^\delta_{\vare_1}\rest\zeta$ are
incompatible.)   

This finishes the description of the strategy $\st$.
\medskip

Let us argue that $\st$ is a winning strategy for Generic. Suppose that 
\[\big\langle\xi_\delta,\langle p^\delta_\gamma,q^\delta_\gamma:\gamma< 
\mu_\delta\cdot \xi_\delta\rangle:\delta<\lambda\big\rangle\] 
is a play of $\tagame(p,\bbP_{\zeta^*})$ in which Generic followed $\st$ and
she constructed the side objects listed in $(\otimes)_\delta$ (for $\delta<
\lambda$) so that demands $(\boxtimes)_1$--$(\boxtimes)_{10}$ are
satisfied. We define a condition $r\in\bbP_{\zeta^*}$ as follows. Let
$\Dom(r)=\bigcup\limits_{\delta<\lambda}\Dom(r_\delta)$. For $\zeta\in
\Dom(r)$ let $r(\zeta)$ be a $\bbP_\zeta$--name for a condition in
$\name{\bbQ}_\zeta$ such that 
\begin{enumerate}
\item[$(\boxtimes)_{11}$] if $\zeta\in w_{\alpha+1}\setminus w_\alpha$,
  $\alpha<\lambda$ (or $\zeta=\alpha=0$), then 
\[\forces_{\bbP_\zeta}\mbox{`` }r(\zeta)\geq r_\alpha(\zeta)\mbox{ and }
r(\zeta)\forces_{\name{\bbQ}_\zeta} \big(\forall\delta{<}\lambda\big) \big( 
\exists j{<}\mu_\delta\big)\big(\forall\vare{<}\xi^*_\delta\big)\big( 
\name{q}^{\mu_\delta\cdot \vare+j}_{\delta,\zeta}\in
\name{G}_{\name{\bbQ}_\zeta}\big)\mbox{ ''}.\]
\end{enumerate}
Clearly $r$ is well defined (remember $(\boxtimes)_6$) and $(\forall\delta<
\lambda)(r_\delta\leq r)$ and $p\leq r$. 
\smallskip

Suppose now that $\delta<\lambda$ and $r'\geq r$. We are going to find
$j<\mu_\delta$ and a condition $r''\geq r'$ such that $(\forall
i<\xi_\delta)(q^\delta_{\mu_\delta\cdot i+j}\leq r'')$. To this end let
$\langle \zeta_\alpha:\alpha\leq\alpha^*\rangle$ be the increasing
enumeration of $w_\delta\cup\{\zeta^*\}$. For $\zeta\leq\zeta^*$ and
$q\in\bbP_\zeta$, let $\st(\zeta,q)$ be a winning strategy of Complete in
$\Game_0^\lambda(\bbP_\zeta,q)$ with the coherence properties given in
\ref{pA.6}.  

By induction on $\alpha\leq\alpha^*$ we will choose conditions $r^*_\alpha,
r^{**}_\alpha\in\bbP_{\zeta_\alpha}$ and $(t)_{\zeta_\alpha}<\mu_\delta$
such that 
\begin{enumerate}
\item[$(\boxtimes)_{12}$] $r'\rest\zeta_\alpha\leq r^*_\alpha$, 
\item[$(\boxtimes)_{13}$] if $i<\xi_\delta$, $j<\mu_\delta$ and 
  $(t^\delta_j)_{\zeta_\beta}=(t)_{\zeta_\beta}$ for $\beta<\alpha$, then
  $q^\delta_{\mu_\delta\cdot i+j}\rest\zeta_\alpha\leq r^*_\alpha$, 
\item[$(\boxtimes)_{14}$] $\langle r^*_\beta\conc r'\rest [\zeta_\beta,
  \zeta^*), r^{**}_\beta\conc r'\rest [\zeta_\beta,\zeta^*):\beta<\alpha
  \rangle$ is a partial legal play of $\Game^\lambda_0(\bbP_{\zeta^*},r')$
  in which Complete uses her winning strategy $\st(\zeta^*,r')$. 
\end{enumerate}
Suppose that $\alpha\leq\alpha^*$ is a limit ordinal and we have already
defined $(t)_{\zeta_\beta}<\mu_\delta$ and $r^*_\beta,r^{**}_\beta\in
\bbP_{\zeta_\beta}$ for $\beta<\alpha$. Let $\zeta=\sup(\zeta_\beta:
\beta<\alpha)$. It follows from $(\boxtimes)_{14}$ that we may pick a
condition $s\in\bbP_\zeta$ stronger than all $r^{**}_\beta$ for
$\beta<\alpha$. Put $r^*_\alpha=s\conc r'\rest [\zeta,\zeta_\alpha)\in
\bbP_{\zeta_\alpha}$. Then plainly $r'\rest \zeta_\alpha\leq r^*_\alpha$ and
$q^\delta_{\mu_\delta\cdot i+j}\rest \zeta\leq r^*_\alpha\rest \zeta$
whenever
\begin{enumerate}
\item[$(\boxtimes)_{15}^{i,j,\alpha}$] $i<\xi_\delta$, $j<\mu_\delta$ and
  $(t^\delta_j)_{\zeta_\beta}=(t)_{\zeta_\beta}$ for all $\beta<\alpha$. 
\end{enumerate}
Now by induction on $\xi\leq\zeta_\alpha$ we show that
$q^\delta_{\mu_\delta\cdot i+j}\rest\xi\leq r^*_\alpha\rest\xi$ whenever
$(\boxtimes)_{15}^{i,j,\alpha}$ holds. For $\xi\leq\zeta$ we are already
done, so assume $\xi\in [\zeta,\zeta_\alpha)$ and we have shown that  
$q^\delta_{\mu_\delta\cdot i+j}\rest\xi\leq r^*_\alpha\rest\xi$ whenever
$(\boxtimes)_{15}^{i,j,\alpha}$ holds. It follows from
$(\boxtimes)_7+(\boxtimes)_9$ that the condition $r^*_\alpha\rest\xi$ forces
in $\bbP_\xi$ that
\[\begin{array}{l}
\mbox{`` the set }\\
\quad\{p(\xi)\}\cup \{r_\alpha(\xi):\alpha<\delta\}\cup\\
\quad\big\{q^\delta_\vare(\xi):\vare=\mu_\delta\cdot i+j\ \&\
i<\xi_\delta\ \&\ j<\mu_\delta\ \&\ \big(\forall\beta<\alpha\big)\big((
t^\delta_j)_{\zeta_\beta}=(t)_{\zeta_\beta}\big)\}\\  
\mbox{ \ has an upper bound in } \name{\bbQ}_\xi\mbox{ ''.} 
  \end{array}\]
and therefore we may use $(\boxtimes)_{10}$ to conclude that 
\[r^*_\alpha\rest\xi\forces\mbox{`` if $(\boxtimes)_{15}^{i,j,\alpha}$ holds,
  then } q^\delta_{\mu_\delta\cdot i+j}(\xi)\leq r_\delta(\xi)\leq r'(\xi)=
r^*_\alpha(\xi)\mbox{ ''.}\] 
The limit stages are trivial and we may claim that
$q^\delta_{\mu_\delta\cdot i+j}\rest \zeta_\alpha\leq r^*_\alpha$ whenever
$(\boxtimes)_{15}^{i,j,\alpha}$ holds. Next, $r^{**}_\alpha$ is determined
by $(\boxtimes)_{14}$. 

Now suppose that $\alpha=\beta+1\leq\alpha^*$ and we have already defined
$r^*_\beta,r^{**}_\beta\in\bbP_{\zeta_\beta}$ and $\langle
(t)_{\zeta_\gamma}: \gamma<\beta\rangle$. It follows from $(\boxtimes)_{11}$
that
\[r^{**}_\beta\forces_{\bbP_{\zeta_\beta}}\mbox{`` }r(\zeta_\beta)
\forces_{\name{\bbQ}_{\zeta_\beta}}\big(\exists\rho<\mu_\delta\big) \big( 
\forall\vare<\xi^*_\delta\big)\big(\name{q}^{\mu_\delta \cdot\vare+
  \rho}_{\delta,\zeta_\beta}\in\name{G}_{\name{\bbQ}_{\zeta_\beta}}\big)\mbox{
  '',}\] 
so we may pick $\rho=(t)_{\zeta_\beta}$ and a condition $s\in
\bbP_{\zeta_\beta+1}$ such that $r^{**}_\beta\leq s\rest\zeta_\beta$ and 
\[s\rest\zeta_\beta\forces_{\bbP_{\zeta_\beta}}\big(\forall \vare< 
\xi^*_\delta\big)\big(\name{q}^{\mu_\delta\cdot\vare+ \rho}_{\delta, 
\zeta_\beta}\leq s(\zeta_\beta)\big).\] 
It follows from  $(\boxtimes)_{13}$+$(\boxtimes)_8$ that then also
$q^\delta_{\mu_\delta\cdot i+j}\rest (\zeta_\beta+1)\leq s$ whenever
$i<\xi_\delta$, $j<\mu_\delta$ and $(t^\delta_j)_{\zeta_\gamma}=
(t)_{\zeta_\gamma}$ for $\gamma\leq\beta$. We let $r^*_\alpha=s\conc r'\rest
(\zeta_\beta,\zeta_\alpha)$ and exactly like in the limit case we argue that
$r'\rest\zeta_\alpha\leq r^*_\alpha$ and $q^\delta_{\mu_\delta\cdot
  i+j}\rest \zeta_\alpha\leq r^*_\alpha$ whenever $i<\xi_\delta$,
$j<\mu_\delta$ and $(t^\delta_j)_{\zeta_\gamma}=(t)_{\zeta_\gamma}$ for
$\gamma\leq\beta$. Again, $r^{**}_\alpha$ is determined by
$(\boxtimes)_{14}$. 

After the induction is completed look at $r''=r^*_{\alpha^*}$ and
$j<\mu_\delta$ such that $t^\delta_j=\langle
(t)_{\zeta_\alpha}:\alpha<\alpha^*\rangle$.   
\end{proof}

\begin{theorem}
  \label{seconddouble}
Assume (a), (b) of \ref{presdouble}. Suppose that $\cU$ is a normal filter
on $\lambda$ and 
\begin{enumerate}
\item[(c)] $\bar{\bbQ}=\langle\bbP_\zeta,\name{\bbQ}_\zeta:\zeta<\zeta^*
\rangle$ is a $\lambda$--support iteration such that for every $\zeta
<\gamma$,   
\[\forces_{\bbP_\zeta}\mbox{`` $\name{\bbQ}_\zeta$ is nicely double {\bf
    b}--bounding over $\bar{\mu},\cU^{\bbP_\zeta}$ ''.}\] 
\end{enumerate}
Then $\bbP_{\zeta^*}=\lim(\bar{\bbQ})$ is nicely double {\bf b}--bounding
over $\bar{\mu},\cU$.
\end{theorem}

\begin{proof}
The proof essentially repeats that of \ref{presdouble} with the following
modifications in the arguments that $\st$ is a winning strategy for Generic
in $\bagame(p,\bbP_{\zeta^*})$. 

We assume that $\big\langle\xi_\delta,\langle p^\delta_\gamma,
q^\delta_\gamma:\gamma<\mu_\delta\cdot\xi_\delta\rangle: \delta<\lambda
\big\rangle$ is a play in which Generic follows $\st$ and the objects listed
in $(\otimes)_\delta$ were constructed on a side. A condition
$r\in\bbP_{\zeta^*}$ is chosen so that $\Dom(r)=\bigcup\limits_{\delta<
  \lambda} \Dom(r_\delta)= \bigcup\limits_{\delta<\lambda} w_\delta$ and for
each $\zeta\in w_{\alpha+1}\setminus w_\alpha$, $\alpha<\lambda$, we have 
\[\begin{array}{l}
\forces_{\bbP_\zeta}\mbox{`` }r(\zeta)\geq r_\alpha(\zeta)\mbox{ and}\\
\qquad\ \ r(\zeta)\forces_{\name{\bbQ}_\zeta} \{\delta<\lambda: (\exists
j<\mu_\delta)(\forall\vare<\xi^*_\delta)(\name{q}^{\mu_\delta\cdot\vare+j}_{\delta,\zeta}\in
\name{G}_{\name{\bbQ}_\zeta})\}\in\cU^{\bbP_{\zeta+1}}\mbox{ ''.}
\end{array}\]
Then, for each $\zeta\in\Dom(r)$, we choose $\bbP_{\zeta+1}$--names
$\name{A}^\zeta_i$ for elements of $\cU$ such that 
\[\forces_{\bbP_\zeta}\mbox{`` }r(\zeta)\forces_{\name{\bbQ}_\zeta}
(\forall\delta\in \mathop{\triangle}\limits_{\delta<\lambda} \name{A}^\xi_i)
(\exists j< \mu_\delta) (\forall\vare<\xi^*_\delta)(\name{q}^{\mu_\delta
\cdot \vare+j}_{\delta,\zeta}\in \name{G}_{\name{\bbQ}_\zeta})\mbox{ ''.}\] 
Finally, we show that for each limit ordinal $\delta<\lambda$,
\[r\forces_{\bbP_{\zeta^*}}\mbox{`` }(\forall\xi\in w_\delta)(\delta\in
\mathop{\triangle}\limits_{\delta<\lambda} \name{A}^\xi_i)\ \Rightarrow\ 
(\exists j< \mu_\delta) (\forall i<\xi_\delta)(q^\delta_{\mu_\delta
\cdot i+j}\in \name{G}_{\bbP_{\zeta^*}})\mbox{ ''.}\] 
For this we start with arbitrary condition $r'\geq r$ such that
\[r\forces_{\bbP_{\zeta^*}}\mbox{`` }(\forall\xi\in w_\delta)(\delta\in
\mathop{\triangle}\limits_{\delta<\lambda} \name{A}^\xi_i)\mbox{ ''}\] 
and we repeat the arguments from the end of the proof of \ref{presdouble} to
find $j<\mu_\delta$ and $r''\geq r'$ such that $(\forall i<\xi_\delta)(
q^\delta_{\mu_\delta\cdot i+j}\leq r'')$. 
\end{proof}

\section{Reasonable ultrafilters with small generating systems} 
Our aim here is to show that, consistently, there may exist a very
reasonable ultrafilter on an inaccessible cardinal $\lambda$ with generating
system of size less than $2^\lambda$. 

\begin{lemma}
\label{lemgetone}
Assume that $G^*\subseteq\bqz$ is directed (with respect to $\leq^0$) and
$\fil(G^*)$ is an ultrafilter on $\lambda$, $r\in G^*$. Let $\bbP$ be a
forcing notion not adding bounded subsets of $\lambda$, $p\in\bbP$ and
let $\name{A}$ be a $\bbP$--name for a subset of $\lambda$ such that
$p\forces_{\bbP}\name{A}\in \big(\fil(G^*)\big)^+$. Then    
\[Y\stackrel{\rm def}{=}\bigcup\big\{Z^r_\delta:\delta\in C^r\mbox{ and }p 
\nVdash_{\bbP}\mbox{`` }\name{A}\cap Z^r_\delta\notin d^r_\delta\mbox{
  ''}\big\}\in\fil(G^*).\]  
\end{lemma}

\begin{proof}
Assume towards contradiction that $Y\notin\fil(G^*)$. Then we may find $s\in 
G^*$ such that $r\leq^0 s$ and $\lambda\setminus Y\in\fil(s)$. Take
$\vare<\lambda$ such that 
\begin{quotation} 
\noindent if $\alpha\in C^s\setminus\vare$, 

\noindent then $Z^s_\alpha\setminus Y\in d^s_\alpha$ and $(\forall A\in 
d^s_\alpha)(\exists \beta\in C^r)(A\cap Z^r_\beta\in d^r_\beta)$. 
\end{quotation}
(Remember \ref{1.3C}.) Now take a generic filter $G\subseteq \bbP$ over
$\bV$ such that $p\in G$ and work in $\bV[G]$. Since $\name{A}^G\in
\fil(s)^+$, we may pick $\alpha\in C^s$ such that $\vare<\alpha$ and
$\name{A}^G\cap Z^s_\alpha\in d^s_\alpha$. Then also $Z^s_\alpha\cap
\name{A}^G\setminus Y\in d^s_\alpha$ and thus we may find $\beta\in C^r$
such that $Z^s_\alpha\cap\name{A}^G\cap Z^r_\beta \setminus Y\in
d^r_\beta$. In particular, $Z^r_\beta\setminus Y\neq\emptyset$, so
$p\forces\name{A}\cap Z^r_\beta\notin d^r_\beta$, and thus $\name{A}^G\cap
Z^r_\beta\notin d^r_\beta$. Consequently $Z^s_\alpha\cap\name{A}^G\cap
Z^r_\beta\setminus Y\notin d^r_\beta$ giving a contradiction.
\end{proof}

\begin{theorem}
\label{lemult}
Assume that 
\begin{enumerate}
\item[(i)] $\lambda$ is strongly inaccessible, $\bar{\mu}=\langle
  \mu_\alpha:\alpha<\lambda\rangle$, each $\mu_\alpha$ is a regular
  cardinal, $\aleph_0\leq\mu_\alpha\leq\lambda$ and $\big(\forall f\in
  {}^\alpha \mu_\alpha\big)\big(\big|\prod\limits_{\xi<\alpha} f(\xi)\big|<
  \mu_\alpha\big)$ for $\alpha<\lambda$;
\item[(ii)] $\bar{\bbQ}=\langle\bbP_\xi,\name{\bbQ}_\xi:\xi<\gamma\rangle$
  is a $\lambda$--support iteration such that for every $\xi<\gamma$,
\[\forces_{\bbP_\xi}\mbox{`` $\name{\bbQ}_\xi$ is reasonably A--bounding
  over $\bar{\mu}$ '';}\]
\item[(iii)] $G^*\subseteq\bqz$ is a $\leq^0$--downward closed
  $\bar{\mu}$--super reasonable family such that $\fil(G^*)$ is an
  ultrafilter on $\lambda$. 
\end{enumerate}
Then   
\[\forces_{\bbP_\gamma}\mbox{`` }\fil(G^*)\mbox{  is an ultrafilter on
  $\lambda$  ''.}\] 
\end{theorem}

\begin{proof}  
The proof is by induction on the length $\gamma$ of the iteration
$\bar{\bbQ}$. So we assume that (i)--(iii) hold and for each $\xi<\gamma$ 
\begin{enumerate}
\item[$(\odot)_\xi$] $\forces_{\bbP_\xi}$`` $\fil(G^*)$ is an ultrafilter on
  $\lambda$  ''. 
\end{enumerate}
Note that (by the strategic $({<}\lambda)$--completeness of $\bbP_\gamma$)
forcing with $\bbP_\gamma$ does not add bounded subsets of $\lambda$, and
therefore $\big(\bqz\big)^\bV\subseteq\big(\bqz\big)^{\bV^{\bbP_\gamma}}$.  

\begin{clx}
\label{cl4}
Assume that 
\begin{enumerate}
\item[(a)]  $\name{A}$ is a $\bbP_\gamma$--name for a subset of $\lambda$
  such that $\forces_{\bbP_\gamma}\name{A}\in \big(\fil(G^*)\big)^+$, 
\item[(b)]  $w\in [\gamma]^{<\omega}$ and $\cT$ is a finite standard
  $(w,1)^\gamma$--tree, and  
\item[(c)]  $\bar{p}=\langle p_t:t\in T\rangle$ is a (finite) tree of
  conditions in $\bar{\bbQ}$, and 
\item[(d)] $r\in G^*$ and $X$ is the set of all $\alpha\in C^r$ for which
  there is a tree of conditions $\bar{q}=\langle q_t:t\in T\rangle$ such
  that $\bar{q}\geq\bar{p}$ and 
\[(\forall t\in T)(\rk(t)=\gamma\ \Rightarrow\ q_t\forces\name{A}\cap
  Z^r_\alpha\in d^r_\alpha).\]
\end{enumerate}
Then $\bigcup\{Z^r_\alpha:\alpha\in X\}\in\fil(G^*)$. 
\end{clx}

\begin{proof}[Proof of the Claim]
Induction on $|w|$. 

If $w=\emptyset$ and so $T=\{\langle\rangle\}$, then the assertion follows
directly from Lemma \ref{lemgetone} (with $p,\bbP$ there standing for
$p_{\langle\rangle},\bbP_\gamma$ here). 

Assume that $|w|=n+1$, $\xi^*=\max(w)$, $w'=w\setminus\{\xi^*\}$ and the
claim is true for $w'$ (in place of $w$) and any $\name{A},\bar{p}$. Let
$\name{\bbP}_{\xi^*\gamma}$ be a $\bbP_{\xi^*}$--name for a forcing notion
with the universe $P_{\xi^*\gamma}=\{p\rest [\xi^*,\gamma):p\in
\bbP_\gamma\}$ and the order relation $\leq_{\name{\bbP}_{\xi^*\gamma}}$
such that    
\begin{quotation}
if $G\subseteq\bbP_{\xi^*}$ is generic over $\bV$ and $f,g\in
P_{\xi^*\gamma}$,\\
then $\bV[G]\models f\leq_{\name{\bbP}_{\xi^*\gamma}[G]} g$ if and only if $(\exists
p\in G)(p\cup f\leq_{\bbP_\gamma} p\cup g)$.
\end{quotation}
Note that $P_{\xi^*\gamma}$ is from $\bV$ but the relation
$\leq_{\name{\bbP}_{\xi^*\gamma}[G]}$ is defined in $\bV[G]$ only. Also
$\bbP_\gamma$ is isomorphic with a dense subset of the composition
$\bbP_{\xi^*}*\name{\bbP}_{\xi^*\gamma}$. 

We are going to define a $\bbP_{\xi^*}$--name $\name{Y}$ for a subset of
$\lambda$. Suppose that $G\subseteq\bbP_{\xi^*}$ is generic over $\bV$ and
work in $\bV[G]$. For $t\in T$ such that $\rk(t)=\gamma$ let $X_t$ consist
of all $\alpha\in C^r$ for which there is $f\in P_{\xi^*\gamma}$ such that 
\[p_t\rest [\xi^*,\gamma)\leq_{\name{\bbP}_{\xi^*\gamma}[G]} f\quad \mbox{ 
  and } \quad f\forces_{\name{\bbP}_{\xi^*\gamma}[G]} \name{A}\cap
Z^r_\alpha\in d^r_\alpha.\]
Let $Y_t=\bigcup\big\{Z^r_\alpha:\alpha\in X_t\big\}$ (for $t\in T$ such
that $\rk(t)=\gamma$). It follows from Lemma \ref{lemgetone} that each $Y_t$
belongs to $\fil(G^*)$ (remember that $\forces_{\bbP_{\xi^*}}$`` $\fil(G^*)$
is an ultrafilter'' by $(\odot)_{\xi^*}$). Hence
\[Y^*\stackrel{\rm def}{=}\bigcap\big\{Y_t:t\in T\ \&\ \rk(t)=\gamma\big\}\in
\fil(G^*).\] 
Note that for each $\alpha\in C^r$, either $Z^r_\alpha\cap Y^*=\emptyset$ or
$Z^r_\alpha\subseteq Y^*$. 

Going back to $\bV$, let $\name{Y}^*,\name{Y}_t,\name{X}_t$ be
$\bbP_{\xi^*}$--names for the objects described as $Y^*,Y_t,X_t$ above. Thus
$\forces_{\bbP_{\xi^*}}\name{Y}^*\in\fil(G^*)$ and we may apply the
inductive hypothesis to $w'$, $T'=\{t\rest\xi^*:t\in T\}$ and $\bar{p}'=
\langle p_{t'}: t'\in T'\rangle\subseteq\bbP_{\xi^*}$. Thus, if $X^*$ is the
set of all $\alpha\in C^r$ for which there is a tree of conditions
$\bar{q}'=\langle q_{t'}':t'\in T'\rangle\subseteq\bbP_{\xi^*} $ such that
$\bar{q}'\geq \bar{p}'$ and  
\[\big(\forall t'\in T\big)\big(\rk(t')=\xi^*\ \Rightarrow\ q_{t'}'
\forces_{\bbP_{\xi^*}} \name{Y}^*\cap Z^r_\alpha\in d^r_\alpha\big),\] 
then $\bigcup\big\{Z^r_\alpha:\alpha\in X^*\}\in\fil(G^*)$. 

Now suppose that $\alpha\in X^*$ is witnessed by $\bar{q}'$ and let $t'\in
T$ be such that $\rk(t')=\xi^*$. Then $q_{t'}' \forces_{\bbP_{\xi^*}}
Z^r_\alpha \subseteq\name{Y}^*$ and hence $q_{t'}'\forces_{\bbP_{\xi^*}}
\alpha \in \name{X}_t$ for all $t\in T$ with $\rk(t)=\gamma$, so we have
$\bbP_{\xi^*}$--names $\name{f}^{t'}_t$ for elements of
$\name{\bbP}_{\xi^*\gamma}$ such that  
\[q_{t'}'\forces_{\bbP_{\xi^*}}\mbox{`` } p_t\rest[\xi^*,\gamma)
\leq_{\name{\bbP}_{\xi^*\gamma}} \name{f}^{t'}_t\ \&\ \name{f}^{t'}_t 
\forces_{\name{\bbP}_{\xi^*\gamma}} \name{A}\cap Z^r_\alpha\in
d^r_\alpha\mbox{ ''}.\]
Now use \ref{dectree} (or just finite induction) to get a 
tree of conditions 
\[\bar{q}''=\langle q_{t'}'':t'\in T'\rangle\subseteq \bbP_{\xi^*}\] 
and objects $g^t_{t'}$ (for $t'\in T'$, $t\in T$, $\rk(t')=\xi^*$,
$\rk(t)=\gamma$) such that $\bar{q}'\leq\bar{q}''$ and  $q_{t'}''
\forces_{\bbP_{\xi^*}}\name{f}^{t'}_t=g^t_{t'}$. Now, for $t\in T$ put  
\begin{itemize}
\item $q_t=q_t''$ if $\rk(t)\leq\xi^*$, and 
\item $q_t=q_{t\rest\xi^*}''\conc g^{t\rest\xi^*}_t$ if $\rk(t)=\gamma$. 
\end{itemize}
It should be clear that $\bar{q}=\langle q_t:t\in T\rangle$ is a tree of
conditions in $\bar{\bbQ}$, $\bar{p}\leq\bar{q}$ and for every $t\in T$ with
$\rk(t)=\gamma$ we have $q_t\forces_{\bbP_\gamma}\name{A}\cap Z^r_\alpha\in
d^r_\alpha$. This shows that $X^*$ is included in the set $X$ defined in the
assumption (d), and hence $\bigcup\big\{Z^r_\alpha:\alpha\in X\big\}\in
\fil(G^*)$.  
\end{proof}

Let $\name{A}$ be a $\bbP_\gamma$--name for a subset of $\lambda$ such that
$\forces_{\bbP_\gamma}\name{A}\in \big(\fil(G^*)\big)^+$ and let $p\in
\bbP_\gamma$. We will find a condition $p^*\geq p$ such that $p^*
\forces_{\bbP_\gamma}\name{A}\in\fil(G^*)$. It will be provided by the
winning criterion $(\circledast)^{\rm tree}_{\bf A}$ of the game
$\agame(p,\bar{\bbQ})$ (see Definition \ref{betterA}; remember $\bbP_\gamma$
is reasonably$^*$ $A(\bar{\bbQ})$--bounding over $\bar{\mu}$ by Theorem 
\ref{verB}). 

Let $\st$ be a winning strategy of Generic in $\agame(p,\bar{\bbQ})$, and
for $\vare\leq\gamma$ and $q\in\bbP_\vare$ let us fix a winning strategy
$\st(\vare,q)$ of Complete in $\Game^\lambda_0(\bbP_\vare,q)$ so that the
coherence demands (i)--(iii) of Proposition \ref{pA.6} are satisfied.
\medskip

We are going to describe a strategy $\st^\boxplus$ of INC in the game
$\Game^\boxplus_{\bar{\mu}}(G^*)$. In the course of a play of
$\Game^\boxplus_{\bar{\mu}}(G^*)$, INC will construct on the side a play of
$\agame(p,\bar{\bbQ})$ in which Generic plays according to $\st$. So suppose
that INC and COM arrived to a stage $\alpha<\lambda$ of a play of
$\Game^\boxplus_{\bar{\mu}}(G^*)$, and they have constructed
\begin{enumerate}
\item[$(\circledast)^\alpha_1$]\qquad $\big\langle I_\gamma,i_\gamma,
\bar{u}_\gamma,\langle r_{\gamma,i},r_{\gamma,i}',(\beta_{\gamma,i},
Z_{\gamma,i},d_{\gamma,i}):i<i_\gamma\rangle:\gamma<\alpha\big\rangle$.  
\end{enumerate}
Also, let us assume that INC (playing according to $\st^\boxplus$) has
written on the side a partial play 
\begin{enumerate}
\item[$(\circledast)^\alpha_2$]\qquad $\big\langle\cT_\gamma,
  \bar{p}^\gamma, \bar{q}^\gamma: \gamma<\alpha\big\rangle$  
\end{enumerate}
of $\agame(p,\bar{\bbQ})$ (in which Generic plays according to $\st$). Let a  
standard tree $\cT_\alpha$ and a tree of conditions $\bar{p}^\alpha=\langle
p^\alpha_t: t\in T_\alpha\rangle$ be given to Generic by the strategy
$\st$ in answer to $(\circledast)^\alpha_2$ (so $|T_\alpha|<\mu_\alpha$). 
\smallskip

On the board of $\Game^\boxplus_{\bar{\mu}}(G^*)$, the strategy
$\st^\boxplus$ instructs INC to play the set
\[I_\alpha\stackrel{\rm def}{=}\{t\in T_\alpha: \rk_\alpha(t)= \gamma\}\] 
and the $<^*_\chi$--first enumeration $\bar{u}_\alpha=\langle
u_{\alpha,i}:i<i_\alpha\rangle$ of $[I_\alpha]^{<\omega}$ (so
$i_\alpha<\mu_\alpha$). Now the two players start playing a subgame of
length $i_\alpha$ to determine a sequence $\langle r_{\alpha,i},
r_{\alpha,i}', (\beta_{\alpha,i},Z_{\alpha,i},d_{\alpha,i}): i<i_\alpha
\rangle$. During the subgame INC will construct on the side a sequence
$\langle \bar{q}^0_i,\bar{q}^1_i:i<i_\alpha\rangle$ of trees of conditions
in $\bbP_\gamma$ so that   
\begin{enumerate}
\item[$(\circledast)_3$] $\bar{q}^\ell_i=\langle q^\ell_{t,i}:t\in
  T_\alpha\rangle$ (for $\ell<2$, $i<i_\alpha$) and for each $t\in
  T_\alpha$, the sequence $\langle q^0_{t,i}, q^1_{t,i}:i<i_\alpha\rangle$
  is a legal play of  $\Game^\lambda_0(\bbP_{\rk_\alpha(t)},p^\alpha_t)$ in which
  Complete uses her winning strategy $\st(\rk_\alpha(t),p^\alpha_t)$. 
\end{enumerate}
Suppose that COM and INC arrive at level $i<i_\alpha$ of the subgame (of
$\Game^\boxplus_{\bar{\mu}}(G^*)$) and   
\begin{enumerate}
\item[$(\circledast)^i_4$]\qquad $\langle r_{\alpha,j},r_{\alpha,j}',
  (\beta_{\alpha,j},Z_{\alpha,j},d_{\alpha,j}):j<i\rangle$ and  $\langle
  \bar{q}^0_j, \bar{q}^1_j:j<i\rangle$   
\end{enumerate}
have been determined and COM has chosen $r_{\alpha,i}\in G^*$. INC's answer  
is given by $\st^\boxplus$ as follows. First, INC takes the
$<^*_\chi$--first tree of conditions $\bar{q}^\diamond$ in $\bar{\bbQ}$ such that 
\begin{enumerate}
\item[$(\circledast)_5^{\rm a}$] $\bar{q}^\diamond=\langle q^\diamond_t:t\in T_\alpha
\rangle$ and $q_t^\diamond\in\bbP_{\rk(t)}$ is an upper bound to the set
$\{p^\alpha_t\}\cup\{q^1_{t,j}:j<i\}$ (for each $t\in T_\alpha$)
\end{enumerate}
(remember $(\circledast)_3$). Then INC lets $X\subseteq C^{r_{\alpha,i}}$ to
be the set of all $\beta\in C^{r_{\alpha,i}}$ greater than $\sup\Big(
\bigcup\limits_{\gamma<\alpha}\bigcup\limits_{j<i_\gamma} Z_{\gamma,j}\cup
\bigcup\limits_{j<i} Z_{\alpha,j} \Big)+890$ and such that 
\begin{enumerate}
\item[$(\circledast)_5^{\rm b}$] there is a tree of conditions $\bar{q}'$ in
  $\bar{\bbQ}$ such that $q^\diamond\leq\bar{q}'$ and  
\begin{center}
if $t\in u_{\alpha,i}$, then $q'_t\forces_{\bbP_\gamma}\name{A}\cap
Z^{r_{\alpha,i}}_\beta\in d^{r_{\alpha,i}}_\beta$.   
\end{center}
\end{enumerate}
Since $u_{\alpha,i}$ is finite, it follows from \ref{cl4} that $\bigcup
\big\{Z^{r_{\alpha,i}}_\beta:\beta\in X\big\}\in\fil(G^*)$.  Then INC picks
also the club $C$ of $\lambda$ such that $C\subseteq C^{r_{\alpha,i}}$ and
$r_{\alpha,i}$ is restrictable to $\langle X,C\rangle$ (see Definition
\ref{restrictable}) and $\min(C)=\min(X)$, and his inning at the stage $i$
of the subgame of $\Game^\boxplus_{\bar{\mu}}(G^*)$ is $r_{\alpha,i}'=
r_{\alpha,i}\rest\langle X,C\rangle$ (again, see Definition
\ref{restrictable}; note that $r_{\alpha,i}'\in G^*$ by \ref{2.10}).  

After this COM answers with $(\beta_{\alpha,i},Z_{\alpha,i},d_{\alpha,i})\in
\#\big(r_{\alpha,i}'\big)$, and then INC chooses (for the construction
on the side) the $<^*_\chi$--first tree of conditions $\bar{q}^0_i$ in 
$\bar{\bbQ}$ such that  $\bar{q}^\diamond\leq\bar{q}^0_i$ and 
\begin{enumerate}
\item[$(\circledast)_6$]  if $t\in u_{\alpha,i}$, then
  $q^0_{t,i}\forces_{\bbP_\gamma}\name{A}\cap Z_{\alpha,i}\in d_{\alpha,i}$.     
\end{enumerate}
Then $\bar{q}^1_i=\langle q^1_{t,i}:t\in T_\alpha\rangle$ is a tree of
conditions determined by the demand in $(\circledast)_3$ and the strategies
$\st(\rk_\alpha(t),p^\alpha_t)$ (for $t\in T_\alpha$); remember the
coherence conditions of \ref{pA.6}. 
\smallskip

This completes the description of how INC plays in the subgame of stage
$\alpha$. After the subgame is finished, INC determines the move
$\bar{q}^\alpha$  of Antigeneric in the play of $\agame(p,\bar{\bbQ})$ which
he is constructing on the side: 
\begin{enumerate}
\item[$(\circledast)_7$] $\bar{q}^\alpha$ is the $<^*_\chi$--first tree of 
conditions $\langle q^\alpha_t:t\in T_\alpha\rangle$ such that $\bar{q}^0_i
\leq\bar{q}^1_i\leq\bar{q}^\alpha$ for all $i<i_\alpha$. 
\end{enumerate}
(There is such a tree of conditions by $(\circledast)_3$; remember $i_\alpha<
\mu_\alpha\leq\lambda$.) 
\medskip
 
This completes the description of the strategy $\st^\boxplus$. Since $G^*$
is $\bar{\mu}$--super reasonable, $\st^\boxplus$ cannot be a winning  
strategy, so there is a play  
\begin{enumerate}
\item[$(\circledast)_8$]\qquad $\big\langle I_\alpha,i_\alpha,
\bar{u}_\alpha,\langle r_{\alpha,i},r_{\alpha,i}',(\beta_{\alpha,i},
Z_{\alpha,i},d_{\alpha,i}):i<i_\alpha\rangle:\alpha<\lambda\big\rangle$
\end{enumerate}
of $\Game^\boxplus_{\bar{\mu}}(G^*)$ in which INC follows $\st^\boxplus$,
but  
\begin{enumerate}
\item[$(\circledast)_9$]\qquad for some $r\in G^*$, for every $\langle
j_\alpha:\alpha<\lambda\rangle\in\prod\limits_{\alpha<\lambda} I_\alpha$ we
have   
\[\{(\beta_{\alpha,i},Z_{\alpha,i},d_{\alpha,i}):\alpha<\lambda\ \&\
i<i_\alpha\ \&\ j_\alpha\in u_{\alpha,i}\}\leq^* \#(r).\] 
\end{enumerate}
Let $\langle \cT_\alpha,\bar{p}^\alpha,\bar{q}^\alpha:\alpha<\lambda
\rangle$ be the play of $\agame(p,\bar{\bbQ})$ constructed on the side by
INC (so this is a play in which Generic uses her winning strategy
$\st$). Since Generic won that play, there is a condition
$p^*\in\bbP_\gamma$ stronger than $p$ and such that for each
$\alpha<\lambda$ the set $\{q^\alpha_t: t\in T_\alpha\ \&\
\rk_\alpha(t)=\gamma\}$ is pre-dense above $p^*$. Note that if we show that    
\begin{enumerate}
\item[$(\circledast)_{10}$] it is forced in $\bbP_\gamma$ that for every
  $\bar{j}=\langle j_\alpha:\alpha<\lambda\rangle\in \prod\limits_{\alpha<
    \lambda} I_\alpha$ we have
\[\big\{(\beta_{\alpha,i},Z_{\alpha,i},d_{\alpha,i}):\alpha<\lambda\ \&\
i<i_\alpha\ \&\ j_\alpha \in u_{\alpha,i}\big\}\leq^* \#(r),\]
\end{enumerate}
then we will be able to conclude that $p^*\forces\name{A}\in\fil(r)$
(remember $(\circledast)_6+(\circledast)_7$ and \ref{qstar}), finishing the
proof of the Theorem. So let us argue that $(\circledast)_{10}$ holds true.  

It follows from the description of $\st^\boxplus$ (see the description of $X$
after $(\circledast)_5^{\rm a}$) that we may choose a continuous increasing
sequence $\langle\delta_\alpha:\alpha<\lambda\rangle\subseteq\lambda$ such
that  
\[\big(\forall\alpha<\lambda\big)\big(\delta_\alpha\leq\beta_{\alpha,0}
\leq\sup\big(\bigcup\limits_{i<i_\alpha}Z_{\alpha,i}\big)<\delta_{\alpha+1}
\big).\] 
Now, we will say that {\em $\beta\in C^r$ is a sick case\/} whenever 
there are $\alpha_0<\alpha_1<\lambda$ and $B\in d^r_\beta$ such that
$Z^r_\beta\subseteq [\delta_{\alpha_0},\delta_{\alpha_1})$ and 
\[\big(\forall\alpha\in [\alpha_0,\alpha_1)\big)\big(\exists t\in I_\alpha
\big)\big(\forall i<i_\alpha\big)\big(t\notin u_{\alpha,i}\ \mbox{ or }\ B\cap
Z_{\alpha,i}\notin d_{\alpha,i}\big).\] 
Using \ref{qstar}(2) one can easily verify that the following two conditions
are equivalent: 
\begin{enumerate}
\item[$(\circledast)_{11}^{\rm one}$] there is $\langle
j_\alpha:\alpha<\lambda\rangle\in \prod\limits_{\alpha<\lambda} I_\alpha$
such that 
\[\big\{(\beta_{\alpha,i},Z_{\alpha,i},d_{\alpha,i}):\alpha<\lambda\ \&\
i<i_\alpha\ \&\ j_\alpha\in u_{\alpha,i}\big\}\nleq^* \#(r),\]
\item[$(\circledast)_{11}^{\rm two}$] there are $\lambda$ many sick cases of
  $\beta\in C^r$. 
\end{enumerate}
Since the forcing with $\bbP_\gamma$ does not add bounded subsets of
$\lambda$, being a sick case is absolute between $\bV$ and
$\bV^{\bbP_\gamma}$. So we may conclude (from $(\circledast)_9$) that
$(\circledast)_{10}$ is true and thus the proof of Theorem \ref{lemult}
is complete.   
\end{proof}

\begin{theorem}
\label{lemsup}
Assume (i) and (ii) of \ref{lemult} and

\begin{enumerate}
\item[$(\alpha)$] $\bar{\kappa}=\langle\kappa_\alpha:\alpha<\lambda\rangle$
  is a sequence of regular cardinals such that for each $\alpha<\lambda$:
\[\mu_\alpha\leq\kappa_\alpha\leq\lambda\quad\mbox{ and }\quad
(\forall\mu<\mu_\alpha)(2^\mu<\kappa_\alpha),\]  
\item[$(\beta)$]  $G^*\subseteq\bqz$ is $\bar{\kappa}$--super reasonable.
\end{enumerate}
Then  $\forces_{\bbP_\gamma}$`` $G^*$ is $\bar{\mu}$--strongly reasonable ''.
\end{theorem}

\begin{proof}
First of all note that the forcing notion $\bbP_\gamma$ is reasonably$^*$
$A(\bar{\bbQ})$--bounding over $\bar{\mu}$ and $\lambda$--proper (see 
\ref{verB}).  Therefore $\forces_{\bbP_\gamma}$`` $\big([G^*]^{\leq\lambda}
\big)^\bV$ is cofinal in  $[G^*]^{\leq\lambda}$'', and consequently
$\forces_{\bbP_\gamma}$`` $G^*$ is $({<}\lambda^+)$--directed (with respect
to $\leq^0$)''. 
\smallskip

Suppose that $\name{\st}^\oplus$ is a $\bbP_\gamma$--name, $p\in\bbP_\gamma$
and    
\[\begin{array}{r}
\forces_{\bbP_\gamma}\mbox{`` }\name{\st}^\oplus\mbox{ is a strategy of
  INC in } \Game^\oplus_{\bar{\mu}}(G^*)\mbox{ such that }\quad \\
\mbox{ all values given by it are from $\bV$ ''}.
\end{array}\] 
We are going to find a condition $p^*\geq p$ and a $\bbP_\gamma$--name 
$\name{g}_\lambda$  such that 
\[\begin{array}{r}
p^*\forces_{\bbP_\gamma}\mbox{`` }\name{g}_\lambda\mbox{ is a play of 
$\Game^\oplus_{\bar{\mu}}(G^*)$ in which INC uses $\name{\st}^\oplus$ but }\
\ \\
\mbox{ COM wins the play ''.}
\end{array}\]
The condition $p^*$ will be provided by the winning criterion
$(\circledast)^{\rm tree}_{\bf A}$ of the game $\agame(p,\bar{\bbQ})$ (see
Definition \ref{betterA}). 
\smallskip

In the rest of the proof whenever we say ``INC chooses/picks $x$ such that''
we mean ``INC chooses/picks the $<^*_\chi$--first $x$ such that''. Let us
fix   
\begin{enumerate}
\item[(i)]  a winning strategy $\st$ of Generic in $\agame(p,\bar{\bbQ})$, 
\item[(ii)] winning strategies $\st(\vare,q)$ of Complete in
  $\Game^\lambda_0(\bbP_\vare,q)$ (for $\vare\leq\gamma$, $q\in\bbP_\vare$)
  such that the coherence conditions of \ref{pA.6} are satisfied. 
\end{enumerate}
\smallskip

We are going to describe a strategy $\st^\boxplus$ of INC in the game
$\Game^\boxplus_{\bar{\kappa}}(G^*)$. In the course of a play of
$\Game^\boxplus_{\bar{\kappa}}(G^*)$, INC will simulate a play of $\agame(p,
\bar{\bbQ})$ and he will consider names for partial plays of
$\Game^\oplus_{\bar{\mu}}(G^*)$ in which INC uses $\name{\st}^\oplus$. Thus
players INC/COM will appear in the play of $\Game^\boxplus_{\bar{\kappa}}
(G^*)$ in $\bV$ and in the play of  $\Game^\oplus_{\bar{\mu}}(G^*)$ in
$\bV^{\bbP_\gamma}$. To avoid confusion we will refer to them as
$\vcom,\vinc$ for $\Game^\boxplus_{\bar{\kappa}}(G^*)$ (in $\bV$) and
$\fcom,\finc$ for $\Game^\oplus_{\bar{\mu}}(G^*)$ (in $\bV^{\bbP_\gamma}$).  
\medskip

So suppose that $\vinc$ and $\vcom$ arrived at a stage $\alpha<\lambda$ of
the play of $\Game^\boxplus_{\bar{\kappa}}(G^*)$ (in $\bV$), and $\vinc$
(playing according to $\st^\boxplus$) has written on the side: 
\begin{enumerate}
\item[$(\oplus)_1^\alpha$] a partial play $\langle \cT_\beta,\bar{p}^\beta,
  \bar{q}^\beta:\beta<\alpha\rangle$ of $\agame(p,\bar{\bbQ})$ in which
  Generic plays according to $\st$, and   
\item[$(\oplus)_2^\alpha$] a $\bbP_\gamma$--name $\name{g}_\alpha=
\langle\name{I}_\beta, \name{i}_\beta,\name{\bar{u}}_\beta,
\name{\bar{x}}_\beta:\beta<\alpha\rangle$  of a partial play of
$\Game^\oplus_{\bar{\mu}}(G^*)$ (in $\bV^{\bbP_\gamma}$) in which $\finc$
uses the strategy $\name{\st}^\oplus$, 
\item[$(\oplus)_3^\alpha$] ordinals $i_\beta<\mu_\beta$ such that $q^\beta_t
\forces\name{i}_\beta=i_\beta$ for every $t\in T_\beta$ with $\rk_\beta(t)=
\gamma$ (for $\beta<\alpha$).  
\end{enumerate}
Note that $\name{I}_\beta$ is a $\bbP_\gamma$--name for a set of size
$<\mu_\beta$ from $\bV$, $\name{\bar{u}}_\beta$ is a $\bbP_\gamma$--name for
an $\name{i}_\beta$--sequence of finite subsets of $\name{I}_\beta$ and
$\name{\bar{x}}_\beta$ is a $\bbP_\gamma$--name for the result of the
subgame of length $\name{i}_\beta$ of level $\beta$.     
\smallskip

Let $\name{I}_\alpha$ be a $\bbP_\gamma$--name for the answer by
$\name{\st}^\oplus$ to the play $\name{g}_\alpha$ of
$\Game^\oplus_{\bar{\mu}}(G^*)$ (in $\bV^{\bbP_\gamma}$).  

Let $\cT_\alpha$ and $\bar{p}^\alpha=\langle p^\alpha_t:t\in
T_\alpha\rangle$ be given to Generic by the strategy $\st$ as an answer to
$(\oplus)^\alpha_1$. Let $\bar{q}^\diamond=\langle q^\diamond_t:t\in
T_\alpha\rangle$ be a tree of conditions in
$\bar{\bbQ}$ such that 
 \begin{enumerate}
\item[$(\oplus)_4^{\rm a}$] $\bar{p}^\alpha\leq \bar{q}^\diamond$ and
  $q^\diamond_{t_0}, q^\diamond_{t_1}$ are incompatible for distinct $t_0,t_1\in
  T_\alpha$ with $\rk_\alpha(t_0)=\rk_\alpha(t_1)$, 
\item[$(\oplus)_4^{\rm b}$] for every $t\in T_\alpha$ with
  $\rk_\alpha(t)=\gamma$ the condition $q^\diamond_t$ decides the value of
  $\name{I}_\alpha$, say $q^\diamond_t\forces_{\bbP_\gamma}$`` $\name{I}_\alpha=I_\alpha^t$ ''. 
\end{enumerate} 
(Note that $\forces_{\bbP_\gamma}\name{I}_\alpha\in\bV$ by the choice of
$\name{\st}^\oplus$; remember \ref{dectree}.)

In the play of $\Game^\boxplus_{\bar{\kappa}}(G^*)$, the strategy $\st^\boxplus$
instructs $\vinc$ to choose the set 
\[I_\alpha=\prod\{I_\alpha^t:t\in T_\alpha\ \&\
\rk_\alpha(t)=\gamma\}\]
and an enumeration $\bar{u}_\alpha=\langle
u_{\alpha,i}: i<i_\alpha\rangle$ of $[I_\alpha]^{<\omega}$. Note that
$|I^t_\alpha|<\mu_\alpha$ for all relevant $t\in T_\alpha$ and
$|T_\alpha|<\mu_\alpha$, so by our assumptions on $\mu_\alpha$ and
$\kappa_\alpha$ we know that $|I_\alpha|<\kappa_\alpha$ (so also
$i_\alpha< \kappa_\alpha$).  

Then, in the play of $\Game^\oplus_{\bar{\mu}}(G^*)$, $\finc$ pretends that
$\fcom$ played an ordinal $\name{i}_\alpha\in [i_\alpha,\lambda)$ and
$\name{\bar{u}}_\alpha=\langle \name{u}_{\alpha,i}:i<\name{i}_\alpha\rangle$
such that 
\[\forces_{\bbP_\gamma}\mbox{`` } \name{\bar{u}}_\alpha\subseteq
[\name{I}_\alpha]^{<\omega}\ \mbox{ and }\ \bigcup\{\name{u}_{\alpha,i}:
i<\name{i}_\alpha\} =\name{I}_\alpha\mbox{ ''}\]
and for each $t\in T_\alpha$ with $\rk_\alpha(t)=\gamma$ we have 
\[q^\diamond_t\forces_{\bbP_\gamma}\mbox{`` }\name{i}_\alpha=i_\alpha \mbox{
and }\name{u}_{\alpha,i}=\{c(t):c\in u_{\alpha,i}\}\mbox{ for }i<i_\alpha
\mbox{ ''.}\]  

Now, both in $\Game^\oplus_{\bar{\mu}}(G^*)$ of $\bV^{\bbP_\gamma}$ and in
$\Game^\boxplus_{\bar{\kappa}}(G^*)$ of $\bV$ the two players start a
subgame. The length of the subgame in $\bV^{\bbP_\gamma}$ may be longer than
$i_\alpha$, but we will restrict our attention to the first $i_\alpha$ steps
of that subgame. In our {\em active case\/} we will have $\name{i}_\alpha
=i_\alpha$, see the choice of $\name{i}_\alpha$ above. When playing the
subgame, $\vinc$ will build a sequence $\langle\bar{q}^0_i,\bar{q}^1_i:
i<i_\alpha\rangle$ of trees of conditions in $\bar{\bbQ}$ such that (in
addition to demands stated later): 
\begin{enumerate}
\item[$(\oplus)_5^{\rm a}$]  $\bar{q}^\ell_j=\langle q^\ell_{t,j}:t\in 
T_\alpha\rangle$, $\bar{q}^\diamond\leq\bar{q}^0_j\leq\bar{q}^1_j\leq
\bar{q}^0_i$ for $\ell<2$, $j<i<i_\alpha$, and  
\item[$(\oplus)_5^{\rm b}$] for each $t\in T_\alpha$,  the sequence
$\langle q^0_{t,i},q^1_{t,i}:i<i_\alpha\rangle$ is a legal play of the game 
$\Game^\lambda_0(\bbP_{\rk_\alpha(t)},q^\diamond_t)$ in which Complete uses
her winning strategy $\st(\rk_\alpha(t),q^\diamond_t)$. 
\end{enumerate}
He (as $\finc$) will also construct a name for a play of a subgame of
$\Game^\oplus_{\bar{\mu}}(G^*)$ of $\bV^{\bbP_\gamma}$ for this stage.  
\smallskip

Suppose that $\vinc$ and $\vcom$ have arrived to a stage $i<i_\alpha$ of the 
subgame and $\vinc$ has determined on the side  $\bar{q}^\ell_j$ for $j<i$, 
$\ell<2$ and a $\bbP_\gamma$--name $\langle\name{z}^\alpha_j:j<i\rangle$ for
a partial play of the subgame of $\Game^\oplus_{\bar{\mu}}(G^*)$ of
$\bV^{\bbP_\gamma}$. Now $\vcom$ chooses $r_{\alpha,i}\in G^*$ which $\vinc$
passes to $\finc$ as an inning of $\fcom$ at the $i$th step of the subgame
of level $\alpha$ of $\Game^\oplus_{\bar{\mu}}(G^*)$ in $\bV^{\bbP_\gamma}$.
There the strategy $\name{\st}^\oplus$ gives $\finc$ an answer
$\name{\delta}_{\alpha,i}<\lambda$.

Next, $\vinc$ picks a tree of conditions $\bar{q}^0_i=\langle q^0_{t,i}:t\in
T_\alpha\rangle$ in $\bar{\bbQ}$ such that 
\begin{enumerate}
\item[$(\oplus)_6^{\rm a}$] $(\forall j<i)(\bar{q}^1_j\leq
  \bar{q}^0_i)$ and $\bar{q}^\diamond\leq\bar{q}^0_i$, and 
\item[$(\oplus)_6^{\rm b}$] for every $t\in T_\alpha$ with
  $\rk_\alpha(t)=\gamma$, the condition $q^0_{t,i}$ decides the value of
  $\name{\delta}_{\alpha,i}$, say $q^0_{t,i}\forces_{\bbP_\gamma}
  \name{\delta}_{\alpha,i}=\delta^t_{\alpha,i}$. 
\end{enumerate}
Then $\vinc$ lets
\[\delta^*_{\alpha,i}=\sup\Big(\big\{\delta^t_{\alpha,i}:t\in T_\alpha\ \&\
\rk_\alpha(t)=\gamma\big\}\cup\bigcup\limits_{\beta<\alpha}
\bigcup\limits_{j<i_\beta} Z_{\beta,j}\cup\bigcup\limits_{j<i}Z_{\alpha,j}
\Big) +890\] 
and in the subgame of $\Game^\boxplus_{\bar{\kappa}}(G^*)$ (in $\bV$) he is 
instructed to put $r'_{\alpha,i}$ such that 
\[C^{r'_{\alpha,i}}=C^{r_{\alpha,i}}\setminus\delta^*_{\alpha,i}\quad\mbox{
  and }\quad d^{r'_{\alpha,i}}_\beta=d^{r_{\alpha,i}}_\beta\mbox{ for }
\beta\in  C^{r'_{\alpha,i}}.\]
(Note that $r'_{\alpha,i}\in G^*$ by \ref{superplus}(2)(ii).)

After this $\vcom$ chooses $(\beta_{\alpha,i},Z_{\alpha,i},d_{\alpha,i})\in
\#(r'_{\alpha,i})$, so $\beta_{\alpha,i}\in C^{r_{\alpha,i}}$,
$\beta_{\alpha,i}\geq \delta^*_{\alpha,i}$ and
$d_{\alpha,i}=d^{r_{\alpha,i}}_{\beta_{\alpha,i}}$. Next $\vinc$ lets
\begin{itemize}
\item $\bar{q}^1_i$ be the tree of conditions in $\bar{\bbQ}$ fully determined by
demand $(\oplus)^{\rm b}_5$ and 
\item $\name{z}^\alpha_i$ be a $\bbP_\gamma$--name for a legal result of
stage $i$ of the subgame of level $\alpha$ of $\Game^\oplus_{\bar{\mu}}
(G^*)$ in $\bV^{\bbP_\gamma}$ such that for each $t\in T_\alpha$ with
$\rk_\alpha(t)=\gamma$ we have  
\[q^0_{t,i}\forces_{\bbP_\gamma} \name{z}^\alpha_i=\big(r_{\alpha,i},
\name{\delta}_{\alpha,i}, (\beta_{\alpha,i},Z_{\alpha,i},d_{\alpha,i})\big).\]
\end{itemize} 
Then the subgame continues.

After all $i_\alpha$ steps of the subgame are completed, $\vinc$ chooses a
tree of conditions $\bar{q}^\alpha=\langle q^\alpha_t:t\in T_\alpha\rangle$
in $\bar{\bbQ}$ such that $(\forall i<i_\alpha)(\bar{q}^1_i\leq
\bar{q}^\alpha)$ and he also lets $\name{\bar{x}}_\alpha$ be a
$\bbP_\gamma$--name for the result of the subgame of level $\alpha$ of
$\Game^\oplus_{\bar{\mu}}(G^*)$  in $\bV^{\bbP_\gamma}$ such that
$\name{\bar{x}}_\alpha\rest i_\alpha=\langle\name{z}^\alpha_i:i<i_\alpha 
\rangle$. Note that all the objects described by
$(\oplus)^{\alpha+1}_1$--$(\oplus)^{\alpha+1}_3$ are determined now.  
\medskip

This completes the description of the strategy $\st^\boxplus$ of INC (i.e.,
$\vinc$) in $\Game^\boxplus_{\bar{\kappa}}(G^*)$. Since $G^*$ is
$\bar{\kappa}$--super reasonable, this strategy cannot be a winning one, so 
there is a play   
\begin{enumerate}
\item[$(\oplus)_7$] $\big\langle I_\alpha,i_\alpha,
\bar{u}_\alpha,\langle r_{\alpha,i},r_{\alpha,i}',(\beta_{\alpha,i},
Z_{\alpha,i},d_{\alpha,i}):i<i_\alpha\rangle:\alpha<\lambda\big\rangle$
\end{enumerate}
of $\Game^\boxplus_{\bar{\kappa}}(G^*)$ in which INC follows $\st^\boxplus$, 
but  
\begin{enumerate}
\item[$(\oplus)_8$] for some $r\in G^*$, for every $\langle
j_\alpha:\alpha<\lambda\rangle\in\prod\limits_{\alpha<\lambda} I_\alpha$ we
have   
\[\{(\beta_{\alpha,i},Z_{\alpha,i},d_{\alpha,i}):\alpha<\lambda\ \&\
i<i_\alpha\ \&\ j_\alpha\in u_{\alpha,i}\}\leq^* \#(r).\] 
\end{enumerate}
Exactly as in the proof of Theorem \ref{lemult} we may argue that then also 
\begin{enumerate}
\item[$(\oplus)_9$] it is forced in $\bbP_\gamma$ that 
\[\big(\forall
\bar{j}\!\in\!\prod\limits_{\alpha<\lambda} I_\alpha\big)\big(\big\{(
\beta_{\alpha,i},Z_{\alpha,i},d_{\alpha,i})\!:\alpha<\lambda\ \&\
i<i_\alpha\ \&\ j_\alpha \in u_{\alpha,i}\big\}\leq^* \#(r)\big).\]   
\end{enumerate}
(See $(\circledast)_{10}$ in the proof of \ref{lemult}.)

Let $\langle \cT_\alpha,\bar{p}^\alpha,\bar{q}^\alpha:\alpha<\lambda
\rangle$ be the play of $\agame(p,\bar{\bbQ})$ constructed on the side by
INC. Generic won that play, so there is a condition $p^*\in\bbP_\gamma$
stronger than $p$ and such that for each $\alpha<\lambda$ the set
$\{q^\alpha_t: t\in T_\alpha\ \&\ \rk_\alpha(t)=\gamma\}$ is pre-dense above
$p^*$.  Also, let $\name{g}_\lambda$ be the $\bbP_\gamma$--name of a play of
$\Game^\oplus_{\bar{\mu}}(G^*)$ (in $\bV^{\bbP_\gamma}$) constructed on the
side in the same run of $\Game^\boxplus_{\bar{\kappa}}(G^*)$ (see
$(\oplus)_2$). We are going to argue that
\begin{enumerate}
\item[$(\oplus)_{10}$] the condition $p^*$ forces (in $\bbP_\gamma$) that
\[\big(\forall\bar{j}\in\prod\limits_{\alpha<\lambda} \name{I}_\alpha\big)
\big(\{(\beta_{\alpha,i},Z_{\alpha,i},d_{\alpha,i}):\alpha<\lambda\ \&\
i<\name{i}_\alpha\ \&\ j_\alpha\in \name{u}_{\alpha,i}\}\leq^* \#(r)\big),\] 
\end{enumerate}
that is 
\[p^*\forces_{\bbP_\gamma}\mbox{`` }\fcom \mbox{ wins the play
}\name{g}_\lambda\mbox{ as witnessed by $r$ ''}.\]
Suppose that $G\subseteq\bbP_\gamma$ is generic over $\bV$, $p^*\in
G$ and let us work in $\bV[G]$.  For every $\alpha<\lambda$ there is a
unique $t=t(\alpha)\in T_\alpha$ such that $\rk_\alpha(t)=\gamma$ and
$q^\alpha_t\in G$, and thus $\big(\name{I}_\alpha\big)^G=I^t_\alpha$,
$\big(\name{i}_\alpha\big)^G=i_\alpha$ and $\big(\name{\bar{u}}_\alpha
\big)^G=\langle (\name{u}_{\alpha,i})^G:i<i_\alpha\rangle$, where $\big(
\name{u}_{\alpha,i}\big)^G=\{c(t):c\in u_{\alpha,i}\}\subseteq
I^t_\alpha$. Suppose that $\bar{j}=\langle j_\alpha:\alpha< \lambda\rangle 
\in\prod\limits_{\alpha<\lambda} I^{t(\alpha)}_\alpha$. For each
$\alpha<\lambda$ fix $j^*_\alpha\in I_\alpha=\prod\{I_\alpha^t:t\in
T_\alpha\ \&\ \rk_\alpha(t)=\gamma\}$ such that $j^*_\alpha(t(\alpha))= 
j_\alpha$. Note that if $j^*_\alpha\in u_{\alpha,i}$, $i<i_\alpha$, then
$j_\alpha\in \big(\name{u}_{\alpha,i}\big)^G$ and therefore
\[\begin{array}{rcl}
\{(\beta_{\alpha,i},Z_{\alpha,i},d_{\alpha,i}):\alpha<\lambda\ \&\
i<i_\alpha\ \&\ j_\alpha\in\big(\name{u}_{\alpha,i}\big)^G\}&\leq^*&\\
\{(\beta_{\alpha,i},Z_{\alpha,i},d_{\alpha,i}):\alpha<\lambda\ \&\
i<i_\alpha\ \&\ j^*_\alpha\in u_{\alpha,i}\}&\leq^*& \#(r)
\end{array}\] 
(remember $(\oplus)_9$). Now $(\oplus)_{10}$ follows and the proof of the
theorem is complete. 
\end{proof}

\begin{corollary}
\label{p.4A}
Assume that $\lambda$ is a strongly inaccessible cardinal. Then there is a
forcing notion $\bbP$ such that
\[\begin{array}{ll}
\forces_\bbP\mbox{``}&\lambda\mbox{ is strongly inaccessible and }
2^\lambda=\lambda^{++}\mbox{ and}\\
&\mbox{there is a strongly reasonable family $G^*\subseteq\bqz$ such that}\\
&\fil(G^*)\mbox{ is an ultrafilter on $\lambda$ and $|G^*|=\lambda^+$, in
  particular }\\
&\mbox{there is a very reasonable ultrafilter on }\lambda\\
&\mbox{with a generating system of size $<2^\lambda$ ''}  
\end{array}\]
\end{corollary}

\begin{proof}
We may start with a universe $\bV$ in which
$\diamondsuit_{S^{\lambda^+}_\lambda}$ holds (and $\lambda$ is strongly
inaccessible). It follows from \ref{1.5x} that (in $\bV$) there is a 
$\leq^0$--increasing sequence $\langle r_\alpha:\alpha<\lambda^+\rangle
\subseteq\bqz$ such that $G^*\stackrel{\rm def}{=}\{r\in\bqz:(\exists
\alpha<\lambda^+)(r\leq^0 r_\alpha)\}$ is super reasonable and
$\fil(G^*)$ is an ultrafilter on $\lambda$. 

Let $\bar{\bbQ}=\langle\bbP_\alpha,\name{\bbQ}_\alpha:\alpha<\lambda^{++}
\rangle$ be a $\lambda$--support iteration of the forcing notion $\bbQ^{\rm
  tree}_{\cD_\lambda}(K_1,\Sigma_1)$ defined in the proof of
\cite[Prop. B.8.5]{RoSh:777}. This forcing is reasonably A--bounding
(by \cite[Prop. 4.1, p. 221]{RoSh:860} and \cite[Thm B.6.5]{RoSh:777}), 
so we may use Theorems \ref{lemult} and \ref{lemsup} to conclude that  
\[\begin{array}{ll}
\forces_{\bbP_{\lambda^{++}}}&\mbox{`` }G^*\mbox{ is strongly reasonable,
  $|G^*|=\lambda^+<2^\lambda$ and }\\
&\quad\fil(G^*)\mbox{ is  ultrafilter on $\lambda$ ''.}
\end{array}\]   
If one analyzes the proof of Theorem \ref{lemsup}, one may notice that even 
\[\forces_{\bbP_{\lambda^{++}}}\mbox{`` }\{r_\alpha:\alpha<\lambda^+\}
  \mbox{ is strongly reasonable ''.}\]   
\end{proof}

\section{A feature, not a bug}
One may wonder if Theorems \ref{lemult}, \ref{lemsup} could be improved by
replacing the assumption that we are working with the iteration of
reasonably A--bounding forcings by, say, just dealing with a nicely double
{\bf a}--bounding forcing. A result of that sort would be more natural and
the fact that we had to refer to an iteration-specific property could be
seen as some lack of knowledge. However, this is {\em a feature, not a
bug\/} as nicely double {\bf a}--bounding forcing notions may cause that
$\fil(G^*)$ is not an ultrafilter anymore.

In this section we assume that $\lambda$ is a strongly inaccessible
cardinal. 

\begin{definition}
  \label{goshfor}
  \begin{enumerate}
\item Let $\bbP^*$ consist of all pairs $p=(\eta^p,C^p)$ such that $\eta^p:
  \lambda\longrightarrow\{-1,1\}$ and $C^p$ is a club of $\lambda$. A binary 
  relation $\leq=\leq_{\bbP^*}$ on $\bbP^*$ is defined by letting $p\leq
  q$\quad if and only if 
\begin{enumerate}
\item[$(\alpha)$] $C^q\subseteq C^p$, $\eta^q\rest \min(C^p)=\eta^p\rest 
  \min(C^p)$,  and 
\item[$(\beta)$]   for every successive members $\alpha<\beta$ of $C^p$ we
  have 
\[\big(\forall\gamma\in [\alpha,\beta)\big)\big(\eta^q(\gamma)=
\frac{\eta^p(\alpha)}{\eta^q(\alpha)}\cdot \eta^p(\gamma)\big).\]
\end{enumerate}
\item For $p\in \bbP^*$ and $\alpha\in C^p$ let 
\[\pos(p,\alpha)\stackrel{\rm def}{=}\big\{\eta^q\rest\alpha:q\in\bbP^*\ \&\
  p\leq q\big\}.\] 
\item For $p\in \bbP^*$, $\alpha<\lambda$ and
  $\nu:\alpha\longrightarrow\{-1,1\}$ we define 
\[\nu*_\alpha p=( \nu\conc\eta^p\rest [\alpha,\lambda),
C^p\setminus\alpha).\]  
(Plainly, $\nu*_\alpha p\in\bbP^*$.)  
  \end{enumerate}
\end{definition}

\begin{remark}
  $\bbP^*$ is a natural generalization of the forcing notion used by 
  Goldstern and Shelah \cite{GoSh:388} to the context of uncountable
  cardinals. 
\end{remark}

\begin{proposition}
  \label{onP}
Let $\bar{\mu}=\langle\mu_\alpha:\alpha<\lambda\rangle$, $\mu_\alpha=
2^{|\alpha|+\aleph_0}$ (for $\alpha<\lambda$). Then $\bbP^*$ is a nicely
double {\bf a}--bounding over $\bar{\mu}$ forcing notion. Also
$|\bbP^*|=2^\lambda$.  
\end{proposition}

\begin{proof}
One easily verifies that the relation $\leq_{\bbP^*}$ is transitive and
reflexive, also plainly $|\bbP^*|=2^\lambda$.

\begin{clx}
\label{cl3}
$\bbP^*$ is $({<}\lambda)$--complete. 
\end{clx}

\begin{proof}[Proof of the Claim]
Suppose that $\delta<\lambda$ and $\langle p_\xi:\xi<\delta\rangle$ is a
$\leq_{\bbP^*}$--increasing sequence of conditions in $\bbP^*$. Let
$C=\bigcap\limits_{\xi<\delta} C^{p_\xi}$ (it is a club of $\lambda$) and
let $\eta:\lambda\longrightarrow\{-1,1\}$  be defined by
\begin{itemize}
\item if $\gamma<\min(C)$ and $\zeta=\min\big(\vare<\delta:\gamma< 
  \min(C^{p_\vare})\big)$,\\
then $\eta(\gamma)=\eta^{p_\zeta}(\gamma)$,  
\item if $\alpha<\beta$ are successive members of the club $C$,
  $\alpha\leq\gamma< \beta$ and $\zeta=\min\big(\vare<\delta: \gamma<
  \min\big(C^{p_\vare}\setminus (\alpha+1)\big)\big)$, then $\eta(\gamma)
  =\eta^{p_\zeta}(\alpha)\cdot \eta^{p_\zeta}(\gamma)$. 
\end{itemize}
Plainly, $\eta$ is well defined and $q\stackrel{\rm def}{=} (\eta,C)\in
\bbP^*$. We claim that $(\forall\xi<\delta)(p_\xi\leq q)$. To this end
suppose $\xi<\delta$. Clearly $C\subseteq C^{p_\xi}$. Now, if $\gamma<
\min(C^{p_\xi})$, then $\eta(\gamma)=\eta^{p_\zeta}(\gamma)$ for some
$\zeta\leq\xi$ such that $\gamma<\min(C^{p_\zeta})$. Since $p_\zeta\leq
p_\xi$, we have $\eta^{p_\zeta}(\gamma)=\eta^{p_\xi}(\gamma)$ and thus 
$\eta^{p_\xi}(\gamma)=\eta(\gamma)$.

Next,  suppose that $\alpha<\beta$ are successive members of $C^{p_\xi}$ and
$\alpha\leq\gamma<\beta$. If $\gamma<\min(C)$ and $\zeta=\min\big(\vare<
\delta:\gamma<\min(C^{p_\vare})\big)$, then $\zeta>\xi$, $\eta(\alpha)=
\eta^{p_\zeta}(\alpha)$ and 
\begin{enumerate}
\item[$(*)^1$] $\eta(\gamma)=\eta^{p_\zeta}(\gamma)=
\frac{\eta^{p_\xi}(\alpha)}{\eta^{p_\zeta}(\alpha)}\cdot
\eta^{p_\xi}(\gamma)=\frac{\eta^{p_\xi}(\alpha)}{\eta(\alpha)}\cdot  
  \eta^{p_\xi}(\gamma)$. 
\end{enumerate}
So assume $C\cap\beta\neq\emptyset$ and let $\alpha'<\beta'$ be successive
members of $C$ such that $\alpha'\leq\alpha\leq\gamma<\beta\leq\beta'$. Let
$\zeta=\min\big(\vare<\delta:\gamma<\min\big(C^{p_\vare}\setminus
(\alpha'+1)\big)\big)$. If $\alpha=\alpha'$, then $\zeta\leq\xi$ and
\begin{enumerate}
\item[$(*)^2$] $\eta(\gamma)=\eta^{p_\zeta}(\alpha)\cdot\eta^{p_\zeta}
(\gamma)=\eta^{p_\zeta}(\alpha)\cdot\frac{\eta^{p_\xi}(\alpha)}{
\eta^{p_\zeta}(\alpha)}\cdot\eta^{p_\xi}(\gamma)=\\
\eta^{p_\xi}(\alpha)\cdot    
\eta^{p_\xi}(\gamma)=\frac{\eta^{p_\xi}(\alpha)}{\eta(\alpha)} \cdot 
\eta^{p_\xi}(\gamma)$ 
\end{enumerate}
(as $\eta(\alpha)=\eta(\alpha')=1$). If $\alpha'<\alpha$, then $\xi<\zeta$
and $\eta(\alpha)=\eta^{p_\zeta}(\alpha')\cdot \eta^{p_\zeta}(\alpha)$, and
hence  
\begin{enumerate}
\item[$(*)^3$] $\eta(\gamma)=\eta^{p_\zeta}(\alpha')\cdot \eta^{p_\zeta}  
(\gamma)=\frac{\eta(\alpha)}{\eta^{p_\zeta}(\alpha)}\cdot \eta^{p_\zeta}
(\gamma)=\\
\frac{\eta(\alpha)}{\eta^{p_\zeta}(\alpha)}\cdot
\frac{\eta^{p_\zeta}(\alpha)}{\eta^{p_\xi}(\alpha)}\cdot \eta^{p_\xi}(
\gamma)=\frac{\eta^{p_\xi}(\alpha)}{\eta(\alpha)}\cdot
\eta^{p_\xi}(\gamma)$.   
\end{enumerate}
Clearly $(*)^1$--$(*)^3$ are what we need to justify \ref{goshfor}(1$\beta$) 
and conclude $p_\xi\leq q$. 
\end{proof}

\begin{clx}
\label{cl2}
Let $p\in \bbP^*$. Then Generic has a nice winning strategy in the game 
$\tagame(p,\bbP^*)$ (see Definition \ref{po27}).  
\end{clx}

\begin{proof}[Proof of the Claim]
We will describe a strategy $\st$ for Generic in
$\tagame(p,\bbP^*)$. Whenever we say {\em Generic chooses $x$ such that\/} 
we mean {\em Generic chooses the $<^*_\chi$--first $x$ such that\/} (and
likewise for other variants). 

During a play of $\tagame(p,\bbP^*)$ Generic constructs on the side
sequences $\langle p_\alpha:\alpha<\lambda\rangle$ and $\bar{\delta}=\langle
\delta_\alpha: \alpha<\lambda\rangle$ so that for each $\alpha<\lambda$:
\begin{enumerate}
\item[(a)] $\bar{\delta}$ is a strictly increasing continuous sequence of
  ordinals below $\lambda$, $p_\alpha\in\bbP^*$ and $\{\delta_\xi:
  \xi\leq\omega+\alpha\}=C^{p_\alpha}\cap (\delta_{\omega+\alpha}+1)$,   
\item[(b)] if $\beta<\alpha$, then $p_\beta\leq p_\alpha$ and
  $\eta^{p_\alpha}\rest \delta_{\omega+\beta}=\eta^{p_\beta} \rest
  \delta_{\omega+\beta}$,     
\item[(c)] $\{\delta_\xi:\xi\leq\omega\}=\{\delta\in C^p:\otp(\delta\cap
  C^p)\leq\omega\}$ and $p_0=p$, 
\item[(d)] $\delta_{\omega+\alpha+1}$ and $p_{\alpha+1}$ are
  determined right after stage $\alpha$ of $\tagame(p,\bbP^*)$.  
\end{enumerate}
So suppose that the two players have arrived to a stage $\alpha<\lambda$ of
a play of $\tagame(p,\bbP^*)$, and Generic has constructed on the side
$\delta_{\omega+\beta+1}$ and $p_{\beta+1}$ for $\beta<\alpha$. If
$\alpha=0$ or $\alpha$ is a limit ordinal, then conditions (a)--(c) and our
rule of taking ``the $<^*_\chi$--first'' fully determine $\{\delta_\xi:\xi\leq
\omega+\alpha\}$ and $p_\alpha$ (the suitable bounds exists essentially by
\ref{cl3}).  

Now Generic chooses an enumeration (without repetition) $\bar{\rho}=
\langle\rho^\alpha_j:j<\mu_\alpha\rangle$ of $\pos(p_\alpha,\delta_{\omega+
  \alpha})$ such that $\rho^\alpha_0=\eta^{p_\alpha}\rest \delta_{\omega+
  \alpha}$. Antigeneric picks a non-zero ordinal $\xi_\alpha<\lambda$ and
the two players start a subgame of length $\mu_\alpha\cdot\xi_\alpha$. In
the course of the subgame, in addition to her innings $p^\alpha_\gamma$,
Generic will also choose ordinals $\vare^\alpha_\gamma=\vare_\gamma<\lambda$
and sequences $\varphi^\alpha_\gamma=\varphi_\gamma:\vare_\gamma
\longrightarrow\{-1,1\}$. These objects will satisfy the following demands
(letting $q^\alpha_\gamma$ be the innings of Antigeneric): 
\begin{enumerate}
\item[(e)] $\delta_{\omega+\alpha}<\vare_{\gamma'}<\vare_\gamma\in
C^{q^\alpha_\gamma}$ and $\varphi_{\gamma'}\rest [\delta_{\omega+\alpha},
\vare_{\gamma'})=\varphi_\gamma \rest [\delta_{\omega+\alpha},
\vare_{\gamma'})$ for $\gamma'<\gamma<\mu_\alpha\cdot\xi_\alpha$,
\item[(f)] if $\gamma=\mu_\alpha\cdot i+2j$, $i<\xi_\alpha$ and
  $j<\mu_\alpha$, then 
\begin{enumerate}
\item[(i)] $\rho^\alpha_j\vtl\varphi_\gamma\vtl
  \varphi_{\gamma+1}$, $\varphi_\gamma=\eta^{q^\alpha_\gamma}\rest
  \vare_\gamma$, and $\varphi_{\gamma+1}(\delta)=
  -\eta^{q^\alpha_{\gamma+1}}(\delta)$ for $\delta\in [\delta_{\omega+
  \alpha},\vare_{\gamma+1})$,  
\item[(ii)] $p^\alpha_0\geq \rho_0^\alpha *_{\delta_{\omega+\alpha}}
  p_\alpha$, $\min(C^{p^\alpha_0})>\delta_{\omega+\alpha}$, and
  $(\varphi_\gamma\rest \vare_{\gamma'})*_{\vare_{\gamma'}}
  q^\alpha_{\gamma'} \leq p^\alpha_\gamma$ for $\gamma'<\gamma$, and  
\item[(iii)] $q^\alpha_\gamma\leq\varphi_\gamma*_{\vare_\gamma}
  p^\alpha_{\gamma+1}\leq \varphi_{\gamma+1}*_{\vare_{\gamma+1}}
  q^\alpha_{\gamma+1}$.  
  \end{enumerate}
\end{enumerate}

So suppose that the two players have arrived to a stage
$\gamma=\mu_\alpha\cdot i+2j$ ($i<\xi_\alpha$, $j<\mu_\alpha$) of the
subgame and $p^\alpha_{\gamma'},q^\alpha_{\gamma'},\varphi_{\gamma'},
\vare_{\gamma'}$ have been determined for $\gamma'<\gamma$. Let
$\varphi=\rho^\alpha_j\conc \bigcup\limits_{\gamma'<\gamma}
\varphi_{\gamma'}\rest [\delta_{\omega+\alpha},\vare_{\gamma'})$. It follows
from (f) that the sequence $\langle(\varphi\rest \vare_{\gamma'}
*_{\vare_{\gamma'}} q^\alpha_{\gamma'}:\gamma'<\gamma\rangle$ is
$\leq_{\bbP^*}$--increasing, so Generic may choose an upper bound
$p^\alpha_\gamma\in\bbP^*$ to it. (Note that necessarily $\varphi\vtl
\eta^{p^\alpha_\gamma}$, $\sup(\vare_{\gamma'}:\gamma'<\gamma)\leq
\min(C^{p^\alpha_\gamma})$.) She plays $p^\alpha_\gamma$ in the subgame and
Antigeneric answers with $q^\alpha_\gamma\geq p^\alpha_\gamma$. Now Generic
lets $\vare_\gamma\in C^{q^\alpha_\gamma}$ be such that
$|C^{q^\alpha_\gamma}\cap \vare_\gamma|=1$ and she puts $\varphi_\gamma=
\eta^{q^\alpha_\gamma}\rest\vare_\gamma$ and she lets $\psi:\vare_\gamma
\longrightarrow \{-1,1\}$ be defined by $\psi\rest \delta_{\omega+\alpha}=
\rho_j^\alpha$ and $\psi(\delta)=-\varphi_\gamma(\delta)$ for $\delta\in 
[\delta_{\omega+\alpha},\vare_\gamma)$. Then Generic plays
$p^\alpha_{\gamma+1}=\psi*_{\vare_\gamma} q^\alpha_\gamma$ as her inning at
stage $\gamma+1$ of the subgame and Antigeneric answers with
$q^\alpha_{\gamma+1}\geq p^\alpha_{\gamma+1}$. Finally, Generic picks
$\vare_{\gamma+1}\in C^{q^\alpha_{\gamma+1}}$ such that
$|C^{q^\alpha_{\gamma+1}}\cap\vare_{\gamma+1}|=1$ and she takes
$\varphi_{\gamma+1}:\vare_{\gamma+1}\longrightarrow\{-1,1\}$ such that
$\varphi_\gamma\vtl\varphi_{\gamma+1}$ and $\varphi_{\gamma+1}(\delta)=
-\eta^{q^\alpha_{\gamma+1}}(\delta)$ for $\delta\in [\vare_\gamma,
\vare_{\gamma+1})$. Plainly, if $i'<i$, $\gamma'=\mu_\alpha\cdot i'+2j$ then
$q^\alpha_{\gamma'}\leq p^\alpha_\gamma$ and $q^\alpha_{\gamma'+1}\leq
p^\alpha_{\gamma+1}$ so both $p^\alpha_\gamma$ and $p^\alpha_{\gamma+1}$ are
legal innings in $\tagame(p,\bbP^*)$. Also easily the demands in (e)+(f) are
satisfied. Moreover, if $j'<j<\mu_\alpha$ then the conditions
$p^\alpha_{\mu_\alpha+j'}$ and $p^\alpha_{\mu_\alpha+j}$ are incompatible.     

After the subgame is over, Generic lets 
\[\varphi=\rho^\alpha_0\conc\bigcup\{\varphi_\gamma\rest
[\delta_{\omega+\alpha},\vare_\gamma):\gamma< \mu_\alpha\cdot
\xi_\alpha\},\] 
and she picks a $\leq_{\bbP^*}$--upper bound
$p_{\alpha+1}'$ to the increasing sequence 
\[\langle (\varphi\rest\vare_\gamma) *_{\vare_\gamma}q^\alpha_\gamma:
\gamma<\mu_\alpha\cdot \xi_\alpha\rangle.\]
Note that $\vare_\gamma\leq\min\big(C^{p_{\alpha+1}'}
\big)$ and $\varphi\rest \vare_\gamma\vtl\eta^{p_{\alpha+1}'}$ for all
$\gamma<\mu_\alpha\cdot \xi_\alpha$, so also $\rho^\alpha_0=\eta^{p_\alpha}
\rest\delta_{\omega+\alpha} \vtl\eta^{p_{\alpha+1}'}$. Also
\begin{enumerate}
\item[(g)] $(\eta^{p_{\alpha+1}'}\rest\vare_\gamma)*_{\vare_\gamma}
  q^\alpha_\gamma \leq p_{\alpha+1}'$ for all $\gamma<\mu_\alpha\cdot
  \xi_\alpha$.  
\end{enumerate}
Let $p_{\alpha+1}\in\bbP^*$ be such that $C^{p_{\alpha+1}}=\{\delta_\xi:
\xi\leq\omega+\alpha\}\cup C^{p_{\alpha+1}'}$ and $\eta^{p_{\alpha+1}}=
\eta^{p_{\alpha+1}'}$ (plainly $p_\alpha\leq p_{\alpha+1}$) and let
$\delta_{\omega+\alpha+1}=\min\big(C^{p_{\alpha+1}'}\big)$. 

This finishes the description of the strategy $\st$. Let us argue that $\st$
is a winning strategy for Generic. To this end suppose that 
\begin{enumerate}
\item[$(\boxplus)$] $\big\langle\xi_\alpha,\langle p^\alpha_\gamma,
  q^\alpha_\gamma:\gamma<\mu_\alpha\cdot\xi_\alpha\rangle:\alpha< \lambda
  \big\rangle$  
\end{enumerate}
is a result of a play of $\tagame(p,\bbP^*)$ in which Generic follows $\st$
and the objects constructed on the side are
\begin{enumerate}
\item[$(\boxplus)^*_\alpha$] $p_\alpha',p_\alpha,\delta_\xi,\langle
  \vare^\alpha_\gamma, \varphi^\alpha_\gamma:\gamma<\mu_\alpha\cdot
  \xi_\alpha \rangle,\langle\rho^\alpha_j:j<\mu_\alpha\rangle$
\end{enumerate}
(and the demands in (a)--(g) are satisfied).  Let $C=\{\delta_\xi:
\xi<\lambda\}$ (so it is a club of $\lambda$) and $\eta=
\bigcup\limits_{\alpha<\lambda} \eta^{p_\alpha}\rest\delta_{\omega+\alpha}$
(clearly $\eta:\lambda\longrightarrow\{-1,1\}$; remember (b)), and let 
$p^*=(\eta,C)$. It is a condition in $\bbP^*$ and it is stronger than all
$p_\alpha$ (for $\alpha<\lambda$) so also $p^*\geq p$. Suppose that
$\alpha<\lambda$ and $p'\geq p^*$. We will show that there is $p''\geq p'$ 
such that for some $j<\mu_\alpha$, the condition $p''$ is stronger than all
$q^\alpha_{\mu_\alpha\cdot i+j}$ for all $i<\xi_\alpha$. Without loss of
generality, $\min(C^{p'})\geq\delta_{\omega+\alpha+1}$. Let $j'<\mu_\alpha$
be such that $\eta^{p'}\rest\delta_{\omega+\alpha}=\rho^\alpha_{j'}$. We
consider two cases now.
\smallskip

\noindent {\sc Case 1:}\qquad $\eta^{p'}(\delta_{\omega+\alpha})=
\eta^{p^*}(\delta_{\omega+\alpha})=\eta^{p_{\alpha+1}}(\delta_{\omega+\alpha})$.\\
Then $\eta^{p'}\rest [\delta_{\omega+\alpha},\delta_{\omega+\alpha+1})=
\eta^{p_{\alpha+1}}\rest[\delta_{\omega+\alpha},\delta_{\omega+\alpha+1})$. Let
$j=2\cdot j'<\mu_\alpha$, and we will argue that $q^\alpha_{\mu_\alpha\cdot
  i+j}\leq p'$ for all $i<\xi_\alpha$. So let $i<\xi_\alpha$, $\gamma=
\mu_\alpha \cdot i+j$. By the choice of $j'$ we know that $\eta^{p'}\rest
\delta_{\omega+\alpha}=\rho^\alpha_{j'}=\eta^{q^\alpha_\gamma} \rest 
\delta_{\omega+\alpha}$ and also 
\[\varphi^\alpha_\gamma\rest [\delta_{\omega+\alpha}, \vare^\alpha_\gamma) = 
\eta^{p_{\alpha+1}} \rest  [\delta_{\omega+\alpha}, \vare^\alpha_\gamma) =
\eta^{p'} \rest  [\delta_{\omega+\alpha}, \vare^\alpha_\gamma).\]
Hence (by (f)(i)) $\eta^{p'}\rest\vare^\alpha_\gamma=\eta^{q^\alpha_\gamma}
\rest \vare^\alpha_\gamma$ and now
\[\begin{array}{l}
q^\alpha_\gamma\leq (\eta^{p'}\rest\vare_\gamma^\alpha)*_{\vare_\gamma} 
q^\alpha_\gamma \leq (\eta^{p'}\rest\delta_{\omega+\alpha})
*_{\delta_{\omega+\alpha}} p_{\alpha+1}'=\\
(\eta^{p'}\rest\delta_{\omega+\alpha+1}) *_{\delta_{\omega+ \alpha+1}}
p_{\alpha+1}'=(\eta^{p'}\rest \delta_{\omega+\alpha+1}) *_{\delta_{\omega+\alpha+1}}
p_{\alpha+1} \leq\\
(\eta^{p'}\rest\delta_{\omega+\alpha+1})*_{\delta_{\omega+\alpha+1}} p^* \leq
 (\eta^{p'}\rest \delta_{\omega+\alpha+1}) *_{\delta_{\omega+\alpha+1}} p'
 =p' 
\end{array}\] 
(for the second inequality remember (g)).
\smallskip

\noindent {\sc Case 2:}\qquad $\eta^{p'}(\delta_{\omega+\alpha})=
-\eta^{p^*}(\delta_{\omega+\alpha})=-\eta^{p_{\alpha+1}}(\delta_{\omega
  +\alpha})$.\\ 
Then $\eta^{p'}(\delta)=-\eta^{p^*}(\delta)=-\eta^{p_{\alpha+1}}(\delta)$
for all $\delta\in [\delta_{\omega+\alpha}, \delta_{\omega+\alpha+1})$. Let
  $j=2\cdot j'+1$ and let us argue that $q^\alpha_{\mu_\alpha\cdot i+j}\leq
  p'$ for all $i<\xi_\alpha$. So let $i<\xi_\alpha$, $\gamma=\mu_\alpha\cdot
  i+j$. Like in the previous case we show that $\eta^{p'}\rest
  \vare^\alpha_\gamma=\eta^{q^\alpha_\gamma} \rest \vare^\alpha_\gamma$ and
  then easily  
\[\begin{array}{l}
q^\alpha_\gamma\leq (\eta^{p'}\rest \vare_\gamma)*_{\vare_\gamma^\alpha} 
q^\alpha_\gamma\leq (\eta^{p'}\rest\delta_{\omega+\alpha+1})
*_{\delta_{\omega+\alpha+1}} p_{\alpha+1}'\leq\\
(\eta^{p'}\rest\delta_{\omega+\alpha+1}) *_{\delta_{\omega+\alpha+1}}
 p'=p'.
\end{array}\]
\end{proof}
\end{proof}

\begin{proposition}
  \label{killG}
Let $\name{\eta}$ be a $\bbP^*$--name such that
\[\forces_{\bbP^*}\name{\eta}=\bigcup\{\eta^p\rest\min(C^p):p\in
\name{G}_{\bbP^*}\}.\]
Then $\forces_{\bbP^*}$`` $\name{\eta}:\lambda\longrightarrow\{-1,1\}$ ''
and for every $s\in\bqz\cap\bV$,   
\[\forces_{\bbP^*}\mbox{`` }\{\alpha<  
\lambda: \name{\eta}(\alpha)=-1\}\in \fil(s)^+\ \mbox{ and }\ \{\alpha<
\lambda:\name{\eta}(\alpha)=1\}\in \fil(s)^+\mbox{ ''.}\]
\end{proposition}

\begin{proof}
It should be clear that $\forces_{\bbP^*}$`` $\name{\eta}:\lambda 
\longrightarrow\{-1,1\}$ '', so let us show the second statement.  Assume
$p\in \bbP^*$, $s\in\bqz$.  Choose  a continuous increasing sequence
$\langle\delta_\xi: \xi<\lambda\rangle\subseteq C^p$ such that for every
$\xi<\lambda$ there is $\alpha=\alpha(\xi)\in C^s$ such that $Z^s_\alpha
\subseteq [\delta_\xi,\delta_{\xi+1})$. Then let
$C=\{\delta_\xi:\xi<\lambda$ is even $\}$ (it is a club of $\lambda$) and
let $\eta:\lambda\longrightarrow\{-1,1\}$ be such that 
\begin{itemize}
\item $\eta\rest [\delta_\xi,\delta_{\xi+1})\in\big\{\eta^p\rest
  [\delta_\xi,\delta_{\xi+1}), -\eta^p\rest [\delta_\xi,
  \delta_{\xi+1})\big\}$,  
\item if $\xi<\lambda$ is even, then $\big\{\delta\in Z^s_{\alpha(\xi)}:
  \eta(\delta)=1\big\}\in d^s_{\alpha(\xi)}$,
\item if $\xi<\lambda$ is odd, then $\big\{\delta\in Z^s_{\alpha(\xi)}:
  \eta(\delta)=-1\big\}\in d^s_{\alpha(\xi)}$.
\end{itemize}
Now note that $(\eta,C)\in\bbP^*$ is a condition stronger than $p$ and it
forces in $\bbP^*$ that 
\[\mbox{`` }\big\{\alpha<\lambda:\name{\eta}(\alpha)=1\big\}\in\fil(s)^+\
\mbox{ and }\ \big\{\alpha<\lambda:\name{\eta}(\alpha) =-1\big\} \in
\fil(s)^+\mbox{ ''.}\] 
\end{proof}

\begin{corollary}
\label{cor4.5}
Assume $\lambda$ is a strongly inaccessible cardinal. Then there is a
forcing notion $\bbP$ such that 
\[\begin{array}{ll}
\forces_\bbP\mbox{``}&\lambda\mbox{ is strongly inaccessible and }
2^\lambda=\lambda^{++}\mbox{ and there is no very}\\
&\mbox{reasonable ultrafilter on $\lambda$ with a
  generating system of size $<2^\lambda$ ''}  
\end{array}\]
\end{corollary}

\begin{proof}
We may start with the universe $\bV$ in which $2^\lambda=\lambda^+$. 

Let $\bar{\bbQ}=\langle\bbP_\alpha,\name{\bbQ}_\alpha:\alpha<\lambda^{++}
\rangle$ be a $\lambda$--support iteration of the forcing notion $\bbP^*$
(see Definition \ref{goshfor}). This forcing is nicely double {\bf
  a}--bounding over $\bar{\mu}$ (where $\mu_\alpha=2^{|\alpha|+\aleph_0}$;
remember Proposition \ref{onP}) and hence $\bbP_{\lambda^{++}}$ is nicely
double {\bf a}--bounding over $\bar{\mu}$ (by Theorem
\ref{presdouble}). Using Theorem \ref{lppcc} we conclude that
$\bbP_{\lambda^{++}}$ does not collapse any cardinals and forces that
$2^\lambda=\lambda^{++}$. Proposition \ref{killG} implies that 
\[\forces_{\bbP_{\lambda^{++}}}\mbox{`` for no family $G^*\subseteq\bqz$ of
  size $<2^\lambda$, $\fil(G^*)$ is an ultrafilter on $\lambda$ ''.}\]
\end{proof}

\begin{problem}
\begin{enumerate}
\item Is it consistent that for some uncountable regular cardinal $\lambda$
  we have that there is no super-reasonable ultrafilter on $\lambda$? Or
  even no very reasonable one? 
\item In particular, are there  super-reasonable ultrafilters on $\lambda$
  in the model constructed for Corollary \ref{cor4.5}? 
\item Do we need the inaccessibility of $\lambda$ for the assertions of
  Corollaries \ref{p.4A}, \ref{cor4.5} (concerning ultrafilters on
  $\lambda$)?   
\end{enumerate}
\end{problem}

%\bibliographystyle{hplain}
%\bibliography{lista,listb,listx,listf,liste,listy}

\begin{thebibliography}{10}

\bibitem{Ab}
Uri Abraham.
\newblock {Lectures on proper forcing}.
\newblock In M.~Foreman A.~Kanamori and M.~Magidor, editors, {\em Handbook of
  Set Theory}.

\bibitem{Ei03}
Todd Eisworth.
\newblock {On iterated forcing for successors of regular cardinals}.
\newblock {\em Fundamenta Mathematicae}, 179:249--266, 2003, math.LO/0210162.

\bibitem{Gi81}
Moti Gitik.
\newblock {On nonminimal $p$-points over a measurable cardinal}.
\newblock {\em Annals of Mathematical Logic}, 20:269--288, 1981.

\bibitem{GoSh:388}
Martin Goldstern and Saharon Shelah.
\newblock {Ramsey ultrafilters and the reaping number---${\rm Con}({\mathfrak 
  r}<{\mathfrak u})$}.
\newblock {\em {Annals of Pure and Applied Logic}}, 49:121--142, 1990.

\bibitem{J}
Thomas Jech.
\newblock {\em {Set theory}}.
\newblock Springer Monographs in Mathematics. Springer-Verlag, Berlin, 2003.
\newblock The third millennium edition, revised and expanded.

\bibitem{RoSh:888}
Andrzej Roslanowski and Saharon Shelah.
\newblock {Lords of the iteration}.
\newblock In {\em Proceedings of the Conference on Boise Extravaganza in Set
  Theory (BEST 2009)}, volume accepted of {\em Contemporary Mathematics
  (CONM)}.
\newblock math.LO/0611131.

\bibitem{RoSh:470}
Andrzej Roslanowski and Saharon Shelah.
\newblock {Norms on possibilities I: forcing with trees and creatures}.
\newblock {\em {Memoirs of the American Mathematical Society}}, 141(671):xii +
  167, 1999.
\newblock math.LO/9807172.

\bibitem{RoSh:860}
Andrzej Roslanowski and Saharon Shelah.
\newblock {Reasonably complete forcing notions}.
\newblock {\em Quaderni di Matematica}, 17, 2005.
\newblock math.LO/0508272.

\bibitem{RoSh:777}
Andrzej Roslanowski and Saharon Shelah.
\newblock {Sheva-Sheva-Sheva: Large Creatures}.
\newblock {\em Israel Journal of Mathematics}, 159:109--174, 2007.
\newblock math.LO/0210205.

\bibitem{RoSh:889}
Andrzej Roslanowski and Saharon Shelah.
\newblock {Generating ultrafilters in a reasonable way}.
\newblock {\em Mathematical Logic Quarterly}, 54:202--220, 2008.
\newblock math.LO/0607218.

\bibitem{Sh:f}
Saharon Shelah.
\newblock {\em {Proper and improper forcing}}.
\newblock {Perspectives in Mathematical Logic}. {Springer}, 1998.

\bibitem{Sh:830}
Saharon Shelah.
\newblock {The combinatorics of reasonable ultrafilters}.
\newblock {\em Fundamenta Mathematicae}, 192:1--23, 2006.
\newblock math.LO/0407498.

\end{thebibliography}

\end{document}